\documentclass{article}

\usepackage{arxiv}

\usepackage[utf8]{inputenc} % allow utf-8 input
\usepackage[T1]{fontenc}    % use 8-bit T1 fonts
\usepackage{hyperref}       % hyperlinks
\usepackage{url}            % simple URL typesetting
\usepackage{booktabs}       % professional-quality tables
\usepackage{amsfonts}       % blackboard math symbols
\usepackage{nicefrac}       % compact symbols for 1/2, etc.
\usepackage{microtype}      % microtypography
\usepackage{lipsum}		% Can be removed after putting your text content
\usepackage{graphicx}
\usepackage{doi}

\usepackage{graphicx}
\usepackage[export]{adjustbox}% for the second solution
\usepackage{amsmath}
\usepackage{breqn}
\usepackage{dutchcal}
\usepackage{amssymb}
\usepackage{microtype}  
\usepackage{hyperref}
\usepackage{caption}
\usepackage{subcaption}
\usepackage[titletoc,title]{appendix}

\title{Response Behavior of Bi-stable Point Wave Energy Absorbers under Harmonic Wave Excitations}

\date{} 	
\author{ {\hspace{1mm}Mohammad A.~Khasawneh} \thanks{Corresponding author.}\\
Department of Mechanical Engineering, Tandon School of Engineering, New York University New York, 11201, USA.\\
Division of Engineering, NYU Abu Dhabi, Abu Dhabi, UAE. \\
	\texttt{mak1011@nyu.edu} \\
	\And
	{\hspace{1mm}Mohammed F.~Daqaq}\\ 
	Department of Mechanical Engineering, Tandon School of Engineering, New York University New York, 11201, USA.\\
Division of Engineering, NYU Abu Dhabi, Abu Dhabi, UAE. \\
	\texttt{mfd6@nyu.edu} \\
}

\renewcommand{\shorttitle}{Response Behavior of Bi-stable PWA under Harmonic Wave Excitations}

\hypersetup{
pdftitle={Response Behavior of Bi-stable Point Wave Energy Absorbers under Harmonic Wave Excitations},
pdfsubject={math},
pdfauthor={Mohammad A.~Khasawneh, Mohammed F. D.~Daqaq},
pdfkeywords={Wave energy, Point wave energy absorber, Bi-stability,  Nonlinearity},
}

\begin{document}
\maketitle
\begin{abstract}
To expand the narrow response bandwidth of linear point wave energy absorbers (PWAs), a few research studies have recently proposed incorporating a bi-stable restoring force in the design of the absorber. Such studies have relied on numerical simulations to demonstrate the improved bandwidth of the bi-stable absorbers. In this work, we aim to understand how the shape of the bi-stable restoring force influences the effective bandwidth of the absorber. To this end, we use perturbation methods to obtain an approximate analytical solution of the nonlinear differential equations governing the complex motion of the absorber under harmonic wave excitations. The approximate solution is validated against a numerical solution obtained via direct integration of the equations of motion. Using a local stability analysis of the equations governing the slow modulation of the amplitude and phase of the response, the loci of the different bifurcation points are determined as function of the wave frequency and amplitude. Those bifurcation points are then used to define an effective bandwidth of the absorber. 
The influence of the shape of the restoring force on the effective bandwidth is also characterized by generating design maps that can be used to predict the kind of response behavior (small amplitude periodic, large amplitude periodic, or aperiodic) for any given combination of wave amplitude and frequency.  Such maps are critical towards designing efficient bi-stable PWAs for known wave conditions.
\end{abstract}


%\pagebreak
\section{Introduction}
Wave energy constitutes one of the most promising and dense renewable energy sources that is yet to be fully exploited. Since the beginning of human civilization, several devices have been devised to harness energy from ocean waves both at small and large scales. Today, methods used to exploit wave energy can be categorized based on their working principle into three different categories; namely,  oscillating water columns,  overtopping devices, and point wave energy absorbers (PWAs) \cite{al2019point}.  Among such approaches, PWAs received the most attention due to their simple design and working principle. In its simplest form, a PWA is composed of a partially-submerged body (buoy) connected through a mooring mechanism to a linear electromagnetic generator attached to the seabed. When waves set the buoy into motion, it pulls a cable connecting it to a linear generator, which creates relative motions between the translating part of the generator (translator) and stationary magnets (stator). As per Faraday's law of induction, this motion induces a current in the generator coils.

Due to their fundamental principle of operation,  traditional PWAs which employ a linear restoring force can work efficiently only near resonance; i.e., when the buoy's velocity is in phase with the wave excitation force. Unfortunately, for typical energetic marine sites, this condition cannot be easily satisfied for reasonably sized systems. Because of the high stiffness of the hydrostatic restoring force emanating from buoyancy, the resonance frequency of the absorber is typically higher than the dominant frequencies in the spectrum of the incoming ocean waves \cite{Falnes2012}. Furthermore,  because linear PWAs has a narrow bell-shaped frequency response with a distinct peak occurring at the resonance frequency, they are incapable of efficiently extracting power from the wide frequency content of the ocean waves; thereby leaving most of the wave energy unexploited.

To overcome such issues, different ideas and solutions have been proposed \cite{drew2009review}. These include the use of active control strategies to bring the natural frequency of the absorber closer to the dominant frequency in the ocean wave spectrum, and the introduction of a bi-stable restoring force to broaden the frequency response bandwidth of the absorber \cite{younesian2017multi,schubert2020performance}.  The idea of utilizing a bi-stable restoring force in PWAs emanated from the field of vibration energy harvesting, where it was shown that vibratory energy harvesters whose potential energy function has two potential wells separated by a potential energy barrier have a broader frequency bandwidth, and are, therefore,  less sensitive to changes in the excitation parameters \cite{daqaq2014role}. 

\begin{figure} [h!]
\centering
\includegraphics[width=0.7\textwidth]{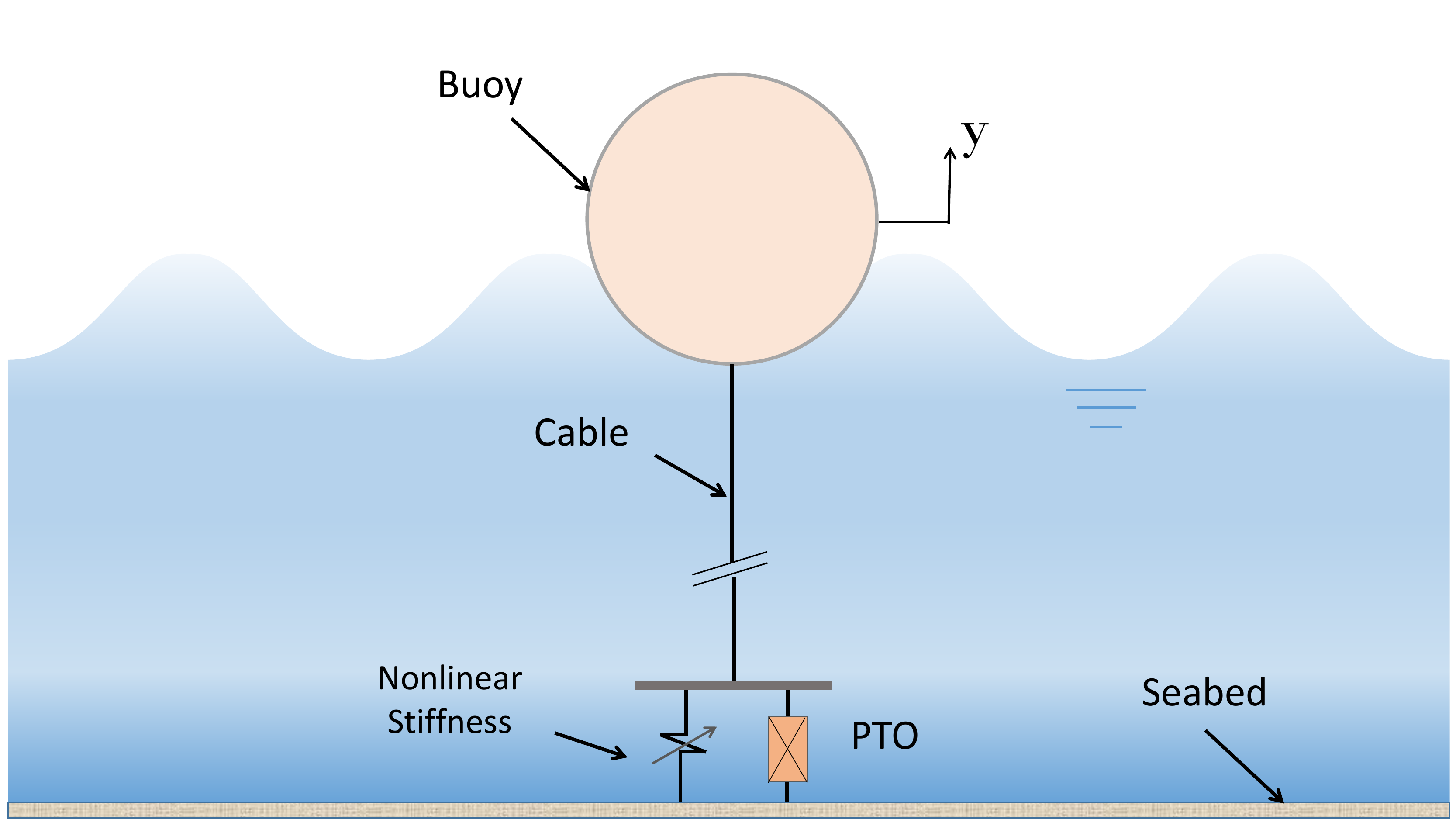}
\vspace{0.0cm}
\caption{Schematic diagram of a bi-stable PWA.}
\label{fig:schematic1}
\end{figure}

\begin{figure} [h!]
\centering
\includegraphics[width=0.5\textwidth]{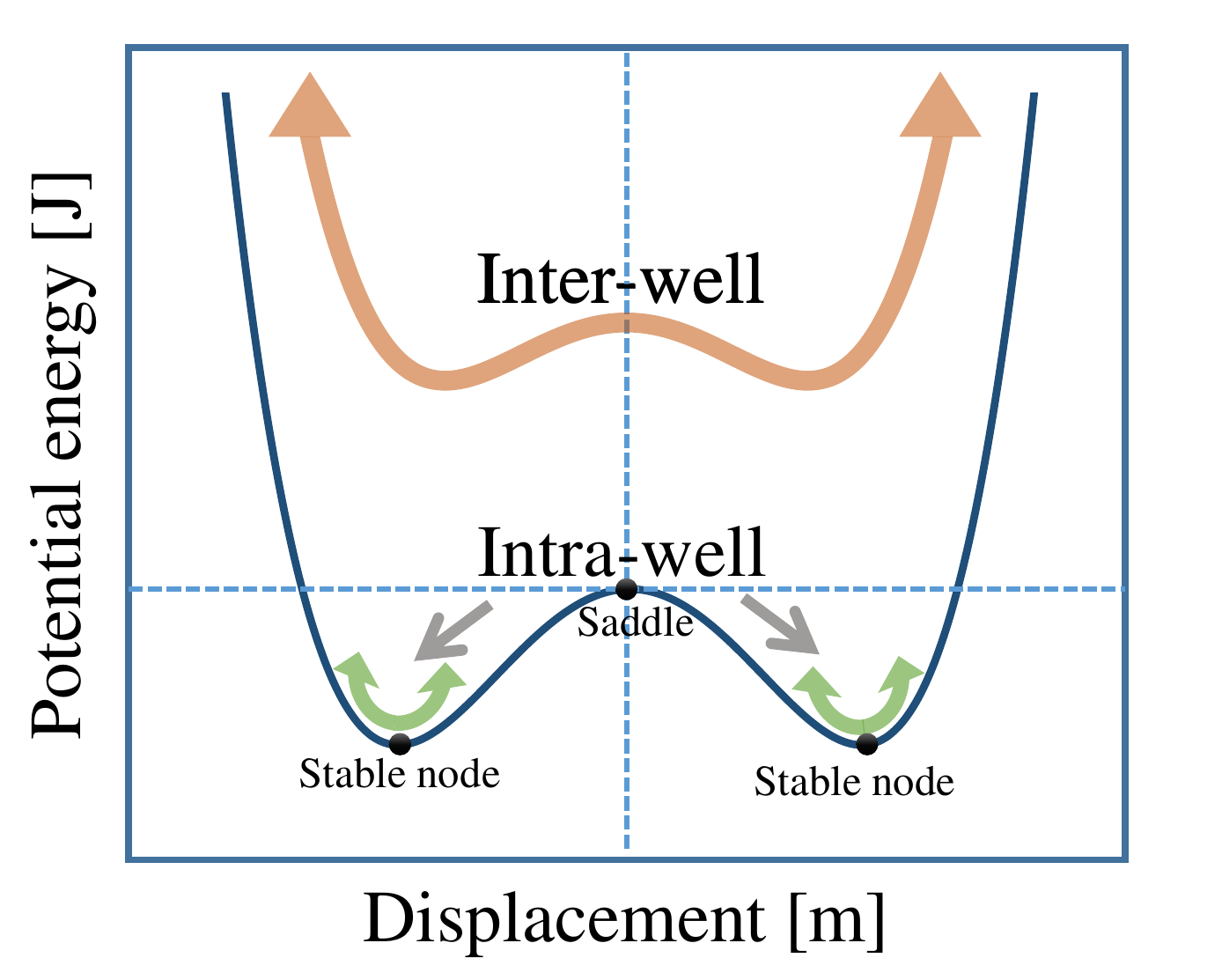}
\vspace{0.0cm}
\caption{Typical potential energy function of a symmetric bi-stable system showing the two types of possible motions (inter- and intra-well).}
\label{fig:Types of Motion}
\end{figure}

A schematic diagram of a bi-stable PWA is shown in Figure \ref{fig:schematic1}. The only difference between the linear and bi-stable PWAs is the addition of the bi-stable spring attachment in parallel with the power take-off unit (PTO). This attachment is specifically designed to create a bi-stable restoring force behavior and can be created by using a set of pre-stretched springs \cite{younesian2017multi} or by using magnetic interactions \cite{schubert2020performance,xi2021high,xiao2017comparative,zhang2019efficiency}. The shape of the potential energy function associated with the bi-stable PWA is shown in Figure \ref{fig:Types of Motion}. The system has two stable equilibria (nodes) separated by a potential barrier (saddle). For some combination of the wave frequency and amplitude, the response of the buoy remains confined to a single potential well (intra-well motion), while for others, the dynamic trajectories overcome the potential barrier causing the buoy to undergo large-amplitude inter-well motions that span the two stable equilibria. This type of large-amplitude motion can extend over a wide spectrum of frequencies which, depending on the shape of the restoring force, can even extend to very low frequencies. These characteristics are key to improving the energy capture from the lower frequency content of the ocean waves.

In terms of performance, a comparison between linear and bi-stable PWAs has revealed a superior bandwidth for the bi-stable absorbers under harmonic waves conditions  \cite{younesian2017multi}. In addition, when considering irregular random waves, results demonstrated superior robustness of the bi-stable absorber with less sensitivity to variations in the frequency content of the waves. A numerical analysis performed in Ref.  \cite{zhang2016oscillating} demonstrated that the performance of the bi-stable PWA is dependent on the shape of its potential energy function and the ability of the dynamic trajectories to escape the potential wells for any given combination of wave frequency and amplitude. Thus, in order to improve the ability of the absorber to perform large-amplitude inter-well motions for a wide range of wave conditions, an adaptive bi-stable absorber, which can adjust the depth of its potential barrier to match the waves excitation intensity was proposed first in Ref. \cite{zhang2018application}, followed by other studies \cite{ zhang2019mechanism,song2020performance}. 

We noticed that, despite the relatively large body of research focused on studying the behaviour of bi-stable PWAs, all of the previous studies relied on purely numerical means without attempting to investigate the complex underlying dynamics of the system via analytical or semi-analytical techniques. Utilizing approximate analytical solutions of the governing nonlinear equations can provide key additional insights into the long-time behavior and bandwidth characteristics of the PWA that cannot be otherwise inferred by relying on numerical simulations alone \cite{daqaq2014role}.

Aiming to bridge this gap, we derive in this paper an approximate analytical solution of the nonlinear differential equations governing the complex motion of the absorber under harmonic wave excitations. Using a local stability analysis of the equations governing the slow modulation of the amplitude and phase of the response, we determine the loci of the different bifurcation points as function of the waves' excitation frequency and amplitude. Those bifurcation points are then used to define an effective bandwidth of the absorber. 
We generate design maps that characterize the influence of the shape of the potential energy function of the PWA on its effective bandwidth. Those maps can be used to predict the type of response behavior of the absorber; e.g. small amplitude periodic, large amplitude periodic, or aperiodic, for any given combination of wave amplitude and frequency.  We believe that such maps are valuable towards designing efficient bi-stable PWAs for known wave conditions.

 The rest of the paper is organized as follow: in Section \ref{sec:mathematical formulation}, the mathematical model governing the motion of the bi-stable PWA is presented and discussed.  In Section \ref{sec:MOMS}, an asymptotic analytical solution of the governing equations is derived using the method of multiple scales for both intra-well and inter-well oscillations. In Section 4, the different bifurcations of the asymptotic solution are identified and analyzed using a stability analysis of the equations governing the slow modulation of the response. In Section \ref{sec: maps}, design maps that characterize the influence of the shape of the bi-stable restoring force of the PWA on its effective bandwidth are generated and discussed.  Finally, in Section \ref{sec:conclusion}, the main conclusions of this work are presented.
\section{Mathematical formulation} \label{sec:mathematical formulation}
\subsection{Governing equations}
\begin{figure} [h!]
\centering
\includegraphics[width=0.75\textwidth]{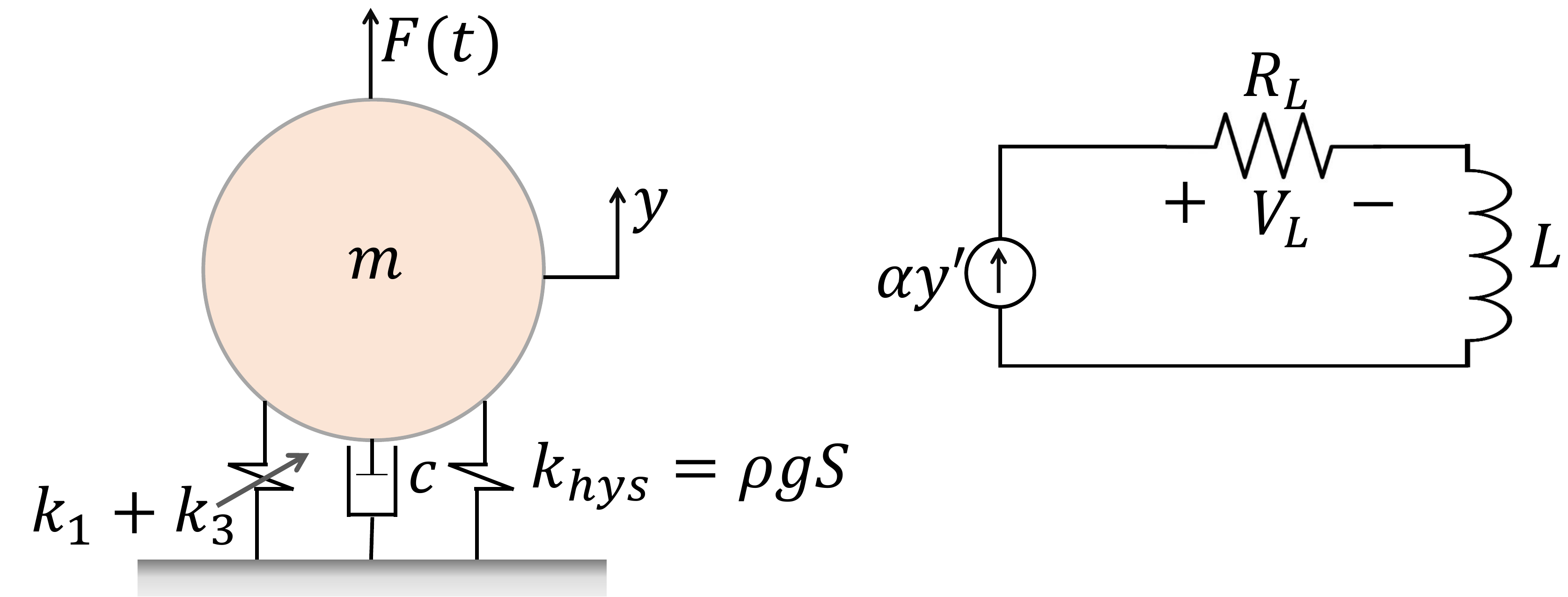}
\vspace{0.0cm}
\caption{A lumped-parameter model of the bi-stable PWA.}
\label{fig:mathmodel}
\end{figure}
Assuming that the buoy undergoes motions in the heave direction only, the equations governing the motion of the absorber can be obtained by applying Newton's second law on the buoy, and Kirchhoff's current law on the harvesting circuit to obtain the following governing equations for the equivalent lumped system shown in Figure \ref{fig:mathmodel}:
 %\begin{equation} \label{eq:Commins}
\begin{subequations} \label{eq:Commins}
\begin{align}
\begin{split} \label{eq:Subeq1}
      &(m+m_{\infty})y'' +\int_{0}^{t} h(t-\tau) y' d\tau +c y' +(k_{hys}-k_{1})y+k_{3}y^3=f_{wave} \cos(\omega t),
\end{split}\\
 &V_L^{'} + \frac{R_L}{L} V_L = \alpha y' .\label{eq:Subeq2}
\end{align}
\end{subequations}
 %\end{equation} 
Here, $y$ represents the displacement of the buoy in the heave direction and the overprime represents a derivative with respect to time, $t$. In Equation (\ref{eq:Subeq1}), $m$ and $m_\infty$ represent, respectively,  the mass of the buoy and the added mass of the fluid. The integral term is used to account for the radiation damping effect; $c$ is a linear viscous damping coefficient; $k_{hys}$ is the stiffness resulting from the hydrostatic buoancy force, which is equal to $\rho g S y$. Here, $\rho$ is the density of water, $g$ is the gravitational acceleration constant, and $S$ is the buoy's wet surface area. The coefficients $k_1>0$ and $k_3>0$ represent, respectively, the linear and cubic coefficients of the nonlinear restoring force added to introduce the bi-stable behavior, and $f_{wave}$ and $\omega$ represent, respectively, the wave amplitude and frequency. Note that to induce a bi-stable potential energy function, $k_1$ must be larger than $k_{hys}$.

In Equation (\ref{eq:Subeq2}), $V_L$ represents the voltage induced by the generator across a purely resistive load $R_L$; $L$ represents the inductance of the harvesting coil, and $\alpha$ is the electromechanical coupling coefficient.

\begin{figure} [h!]
\centering
\includegraphics[width=0.5\textwidth]{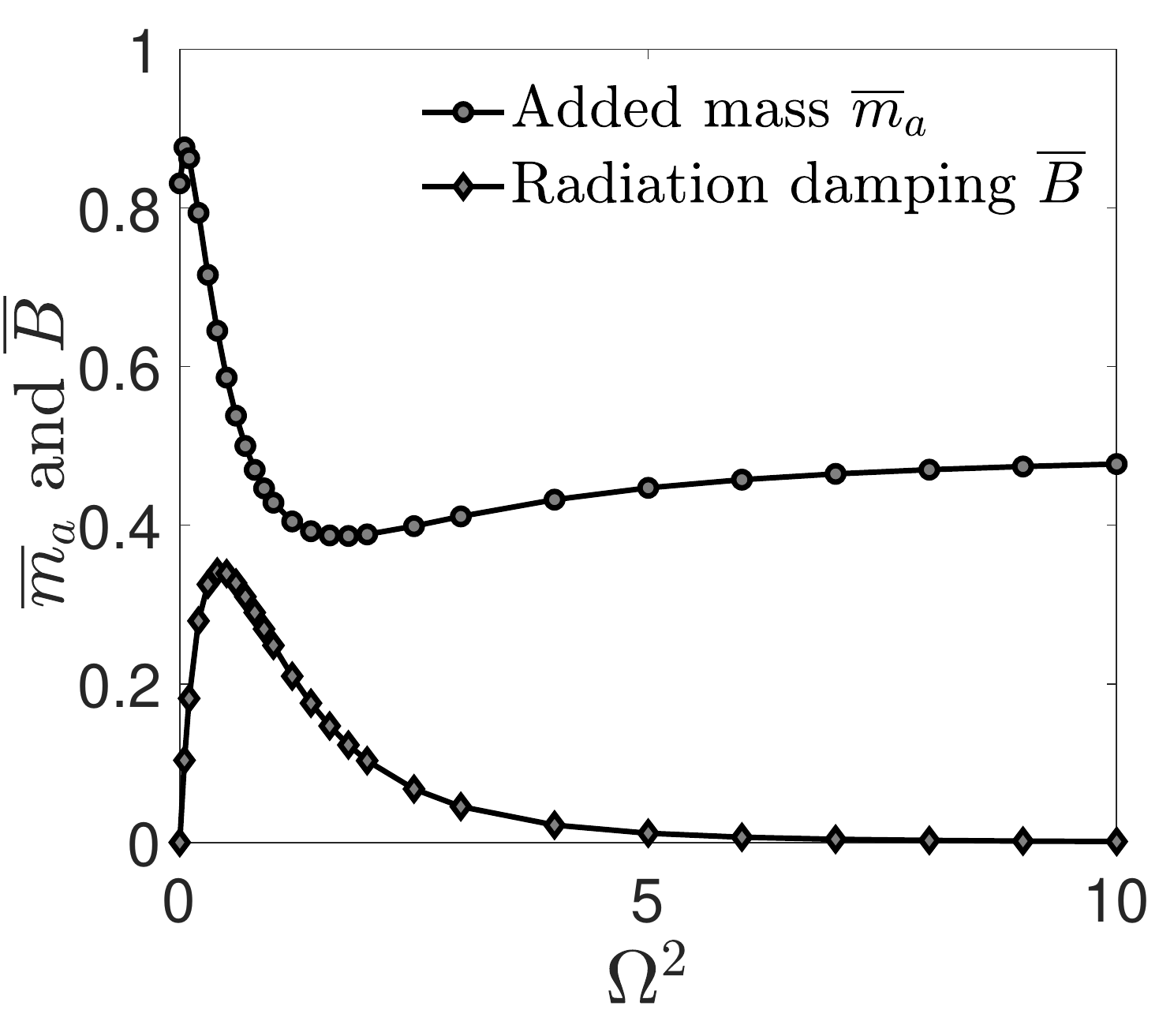}
\vspace{0.0cm}
\caption{Variation of the normalized added mass $\overline{m}_a$, and the normalized radiation damping coefficient $\overline{B}$ with the square of the normalized wave frequency.}
\label{fig:addedMass and Damping}
\end{figure}

The added mass and the radiation damping depend on the fluid velocity field around the buoy, and hence, are a function of the wave frequency.  Following the asymptotic analysis provided in the work of  \textit{Holme} \cite{hulme1982wave} for a spherical buoy of radius $R$, the curves governing the dependence of the normalized added mass denoted here as $\overline{m}_a=\frac{m_a}{M}$ $(M=\frac{2}{3}\pi R^3 \rho)$, and the normalized radiation damping coefficient $\overline{B}=\frac{B}{M\omega}$ on the normalized wave frequency $\Omega=\omega/\sqrt{g/R}$ are shown in Figure \ref{fig:addedMass and Damping}.

The radiation kernel can be further related to the radiation damping coefficient $B(\omega)$ via the following equation \cite{ogilvie1969rational}:
\begin{equation} \label{eq:Ogilive int}
    h(t) = \frac{2}{\pi} \int_0^\infty B(\omega)\cos(\omega t) d\omega,
\end{equation}
which can be discretized as following: 
\begin{equation} \label{eq:ogilvie}
    h(t)= \frac{2}{\pi}\lim_{\delta \omega\to 0}\sum_{i=1}^{\infty} B(\omega_i)\cos(\omega_i t)\delta\omega.
    \end{equation}
Here, the $B(\omega_i)$ are the values of the radiation damping coefficient calculated at the discrete value of $\omega_i$ shown in Figure \ref{fig:addedMass and Damping}. 

The amplitude of the force $f_{wave}$ acting on the buoy due to the incident waves can be related to the radiation damping coefficients by employing \textit{Haskinds's} relation  \cite{haskind2010exciting,newman1962exciting}, which states that
\begin{equation} \label{eq:Haskind}
    f_{wave}(\omega) = A_{wave} \sqrt{\frac{2 \rho g^3}{\omega^3} B(\omega)},
\end{equation}
where $A_{wave}$ is the regular wave amplitude. The physical interpretation of Equation (\ref{eq:Haskind}) is fairly intuitive, as it relates the tendency of the buoy to radiate waves in a certain direction to the excitation forces the buoy experiences from waves propagating in the same direction. For more insight, the reader can refer to the seminal works of \textit{Haskind} and \textit{Newman}\cite{haskind2010exciting,newman1962exciting}.

\subsection{Approximation of the convolution integral} \label{sec: convolution }
As aforestated, the goal of our work is to obtain approximate analytical solutions of the equations governing the motion of the buoy in order to gain deeper insights into the influence of the shape of the restoring force on the performance of the absorber. To achieve this goal, we will use the method of multiple scales \cite{nayfeh2008perturbation}. In order to facilitate the implementation of the method, we obtain in this section an approximation of the convolution integral governing the radiation damping in Equation (\ref{eq:Subeq1}).

To this end, we first use the fact that the output, $z(t)$, of a linear dynamical system can be expressed as a convolution integral between its input $u(t)$, and an impulse response function $h(t)$, in the form
\begin{equation} \label{eq:state-space output}
    \textbf{z}(t)=\int_0^t h(t-\tau) \textbf{u}(\tau) d\tau \approx \textbf{C}_r \textbf{x}(t).
\end{equation}
Here, $\textbf{z(t)} \in \mathbb{R}^q$ is the output vector,  $\textbf{x(t)} \in \mathbb{R}^n$ is the state vector and $\textbf{C}_r$ is a $q \times n $ output matrix governed by the following linear state-space equation
\begin{align} \label{state-space}
        \textbf{x}'(t) & =\textbf{A}_r\textbf{x}(t)+\textbf{B}_r\textbf{u}(t).
\end{align}
 Using the Eigensystem Realization Algorithm (ERA) detailed in Ref. \cite{brunton2019data}, the convolution integral can be expressed in terms of the realized state-space matrices $\textbf{A}_r, \textbf{B}_r$ and $\textbf{C}_r$ as:
 \begin{equation} \label{eq:h(t)=CexpAB}
     h(t)=\left(M \sqrt{\frac{g}{R}} \right)\textbf{C}_r e^{\textbf{A}_rt} \textbf{B}_r,
 \end{equation}
where the realized state-space matrices are :
\begin{equation*}
\textbf{A}_r = 0.8
\begin{pmatrix}
-1 & 1 & 1\\
-1 & 0 & 0\\
-1 & 0 & -2
\end{pmatrix},
\end{equation*}
\begin{equation*}
\textbf{B}_r =
\begin{pmatrix}
-0.48 & -0.02 & -0.22
\end{pmatrix}^T,
\end{equation*}
and
\begin{equation*}
\textbf{C}_r =
\begin{pmatrix}
-0.46 & 0 & 0.18
\end{pmatrix}.
\end{equation*}
Note that the numerical values appearing in the realized state-space matrices are general for any spherical buoy of radius, $\textit{R}$. For more details on the ERA procedure, the interested reader can refer to Appendix \ref{apndx: ERA}.  

The matrix exponential $e^{\textbf{A}_rt}$ in Equation (\ref{eq:h(t)=CexpAB}) can be further expressed in the following form
\begin{equation} \label{eq:state-transition}
    e^{\textbf{A}_rt}=\mathcal{L}^{-1}\left(s\textbf{I} - \textbf{A}_r\right)^{-1}.
\end{equation}
Here, $\mathcal{L}^{-1}$ is the inverse Laplace transform, and \textbf{I} is the identity matrix. Substituting Equation (\ref{eq:state-transition}) into Equation (\ref{eq:h(t)=CexpAB}), we obtain the following analytical expression for $h(t)$:
\begin{equation} \label{eq:h(t) analytical}
    h(t)=\left( M \sqrt{\frac{g}{R}} \right) e^{-\mu t}\left(\lambda_1 + \lambda_2 \cos(\mu t) + \lambda_3 \sin(\mu t)\right), 
\end{equation}
where $\mu , \lambda_1, \lambda_2$ and $\lambda_3$ are constants listed in Table \ref{Table:h(t) constants}, and are valid for any spherical buoy of radius, R.
\begin{table}
\caption{Numerical values of the constants appearing in Equation (\ref{eq:h(t) analytical}).}
\label{Table:h(t) constants}
\begin{center}
\begin{tabular}{c  c} 
\hline
 $Parameter$ & $Value$  \\ 
 \hline
 $\mu$       & 0.8    \\  
 $\lambda_1$ & -0.44  \\  
 $\lambda_2$ & 0.62   \\    
 $\lambda_3$ & 0.24   \\
 \hline
\end{tabular} 
\end{center} 
\end{table} 

\begin{figure} [h!]
\centering
\includegraphics[width=0.6\textwidth]{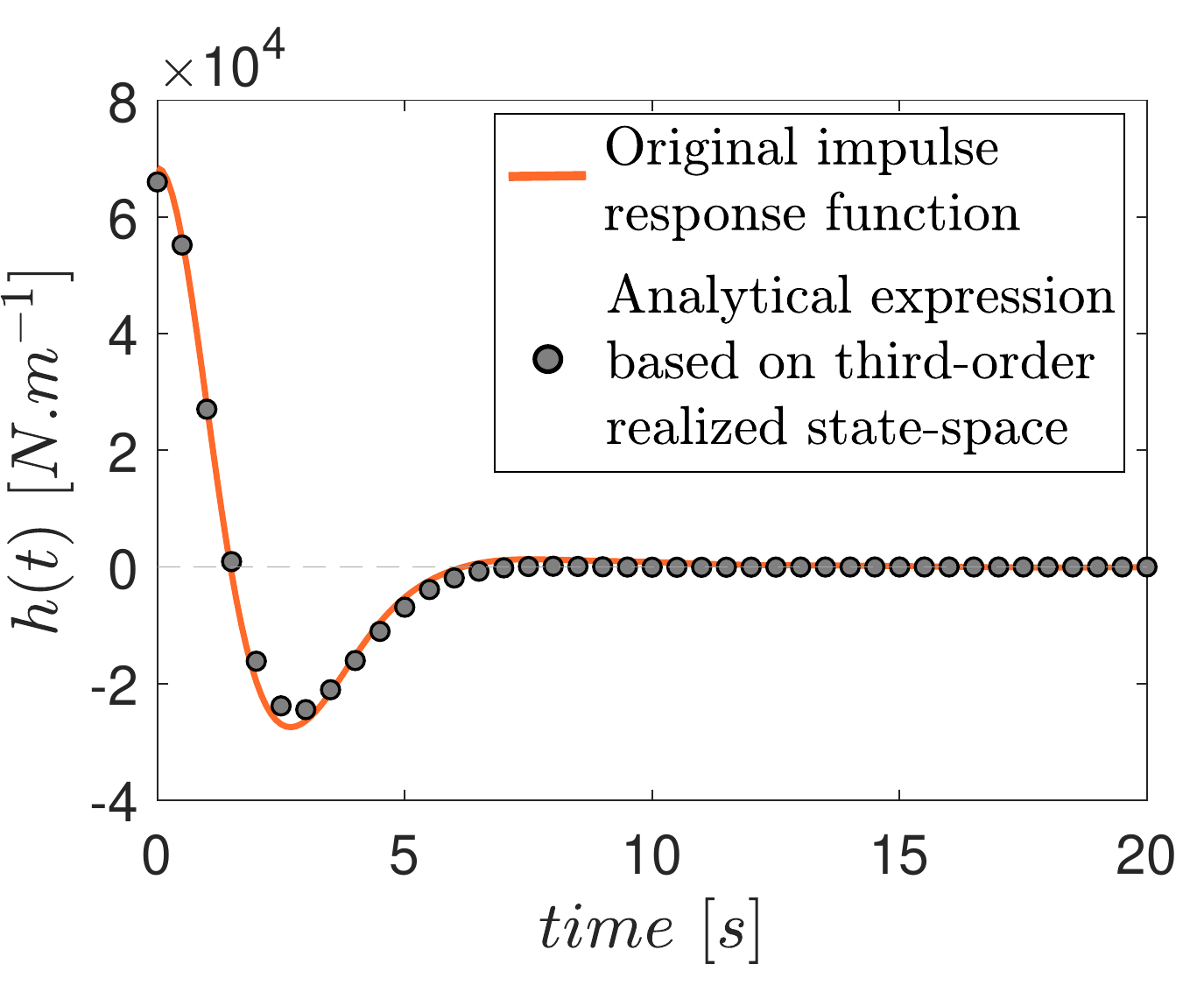}
\vspace{0.0cm}
\caption{Impulse response function $h(t)$ for a hemispherical buoy of radius 5 [m]. Solid line represents the original function obtained by applying Equation (\ref{eq:ogilvie}). Circles represent the values of $h(t)$ obtained from the analytical expression in Equation (\ref{eq:h(t) analytical}) which is based on the third-order realized state-space estimated through employing the eigensystem realization algorithm.}
\label{fig:h(t)}
\end{figure}

Figure \ref{fig:h(t)} depicts a comparsion between the analytical expression of $h(t)$ as obtained using Equation (\ref{eq:h(t) analytical}) and that obtained using the original expression of Equation (\ref{eq:ogilvie}) for a hemispherical buoy of radius $R=5$ [m]. It can be clearly seen that the analytical expression for $h(t)$ is in an excellent agreement with the original impulse response function obtained via Equation (\ref{eq:ogilvie}).
\\

\section{Approximate Analytical Solution} \label{sec:MOMS}
In this section, we employ the method of multiple scales \cite{nayfeh2008perturbation} to obtain an approximate analytical solution of Equation (\ref{eq:Commins}). As aforementioned, and shown in Figure \ref{fig:Types of Motion}, there are two possible steady-state motions. Those that are confined to the potential well and known as intra-well oscillations and those that span the two potential wells and known as the inter-well motions. We first obtain analytical approximations of the motion trajectories within one potential well, then we seek approximate solutions that govern the inter-well motions.

\subsection{Local intra-well response}
In this subsection, we obtain an approximate solution for Equation (\ref{eq:Commins}) when the buoy undergoes motions within a single potential well; i.e, the dynamics around one stable equilibrium node. Upon using the radius of the buoy, $R$, as a length scale and the $\sqrt{g/R}$ as a time scale, we obtain the following dimensionless equivalent equation of motion:
 \begin{subequations} \label{eq:ComminsWN}
  \begin{align}
      \ddot{Y}&+\delta_1 \int_{0}^{t^*} \overline{h}(t^*-\tau) \dot{Y}(\tau) d\tau +\delta_2 \dot{Y}-\omega^2_n Y+\gamma Y^3=g_{wave} \cos(\Omega t^*),\\
  &\dot{v} + \theta v=  \dot{Y},
 \end{align}
 \end{subequations}
where the overdot represents the derivative with respect to the nondimensional  time, $t^*$, and the other parameters and nondimensional groups are defined as:
 \begin{center}
 $\delta_1=\frac{m+m_\infty}{M}$,  $\delta_2=\frac{c}{(m+m_{\infty})}\sqrt{\frac{R}{g}}$,  $\omega_n=\sqrt{\frac{(k_1-\rho g S) R }{(m+m_{\infty})g}}$,  $\gamma=\frac{R^3 k_3}{(m+m_{\infty})g}$,   $g_{wave}=\frac{A M}{R(m+m_{\infty})} \Omega\sqrt{\frac{3\overline{B}}{\pi}}$ ,  $M=\frac{2}{3}\pi R^3 \rho$,  $Y=\frac{y}{R}$,  $v=\frac{V}{\alpha R}$,  $\overline{h}=\frac{h}{M}\sqrt{\frac{R}{g}}$,  $\theta=\frac{R_L}{L } \sqrt{\frac{R}{g}}$,  $t^*=t\sqrt{g/R}$.
 \end{center}
We expand the dynamics governed by Equation (\ref{eq:ComminsWN}) about the stable node, $Y_s= \sqrt{\omega^2_n/\gamma}$, by introducing the  transformation $Y(t^*)=z(t^*) - Y_s$ into Equation (\ref{eq:ComminsWN}). This yields
 \begin{subequations} \label{eq:Commins intra-well}
 \begin{align}
 \begin{split}
      \ddot{z} &+\delta_1 \int_{0}^{t^*} \overline{h}(t^*-\tau) \dot{z}(\tau) d\tau +\delta_2 \dot{z}+\omega^2_0 z+\eta z^2+\gamma z^3=g_{wave} \cos(\Omega t^*),
 \end{split}\\
  &\dot{v} + \theta  v = \dot{z}.
 \end{align}
 \end{subequations}
Here $z(t)$ represents the dynamic trajectories within the potential well, $\omega_o=\sqrt{2}\omega_n$ is the corresponding local frequency of oscillations.  It is worth noting that the expansion about the stable node introduces a new quadratic term which captures the asymmetric nature of the motion trajectories about the stable node, $Y_s$.

To implement the method of multiple scales on Equation (\ref{eq:Commins intra-well}), we introduce multiple time scales, $T_n= \epsilon^n t^*$, $n=0,1,2$. It follows that
\begin{equation} \label{eq:time scales_intra}
\begin{split}
    \dot{(.)} &= D_0 + \epsilon D_1 + \epsilon^2 D_2 + O (\epsilon^3), \\
    \ddot{(.)} &= D_0^2 + \epsilon^2 D_1 + 2\epsilon D_0 D_1 + 2 \epsilon^2 D_0 D_2 + O (\epsilon^3),
\end{split}
\end{equation}
where $\epsilon$ is a scaling parameter, and $D_n$ is the temporal derivative operator with respect to the time scales, $T_n$. Next, we seek an expansion for the response $z(t^*)$ and voltage $v(t^*)$ as
\begin{subequations} \label{eq:y expansion_intra}
\begin{align}
        z(t^*, \epsilon)&=z_0(T_0,T_1,T_2) + \epsilon z_1(T_0,T_1,T_2) +\epsilon^2 z_2(T_0,T_1,T_2) + O (\epsilon^3), \\
       v (t^*, \epsilon)&= v_0 (T_0,T_1,T_2) + \epsilon v_1(T_0,T_1,T_2) +\epsilon^2 v_2(T_0,T_1,T_2) + O (\epsilon^3).
\end{align}
\end{subequations}
Based on the typically small values of the damping coefficients, quadratic and cubic nonlinearities, electromechanical coupling, the time constant of the harvesting circuit, and the excitation amplitude, we scale them to be at order $\epsilon^2$; that is 
\begin{center}
    $\delta_1=\epsilon^2 \delta_1$,  $\delta_2=\epsilon^2 \delta_2$,  $\eta= \epsilon^2 \eta$,  $\gamma=\epsilon^2 \gamma$ ,  $g_{wave}=\epsilon^2 g_{wave}$. 
    \end{center}
Since large-amplitude intra-well motions occur near the primary resonance of the system; i.e. when $\Omega$ is near $\omega_o$,  we limit the intra-well analysis to wave frequencies that are close to $\omega_o$ by introducing the detuning parameter $\sigma$ such that
\begin{equation} \label{eq:detuning paraemter_intra}
    \Omega=\omega_o + \epsilon^2 \sigma.
\end{equation}
Upon substituting Equations (\ref{eq:time scales_intra} - \ref{eq:detuning paraemter_intra}) into Equation (\ref{eq:Commins intra-well}), then collecting terms of equal powers of $\epsilon$, we obtain the following perturbation problems at the different scales:

\noindent
\underline{$O(\epsilon^0)$}:\\
\begin{subequations} \label{eq:order 0 intra-well}
\begin{align}
    D_0^2 z_0 + \omega_o^2 z_0 &=0,\\
    D_0 v_{0} + \theta  v_0 &= D_0z_0,
\end{align}
\end{subequations}\\
\underline{$O(\epsilon^1)$}:\\
\begin{subequations} \label{eq:order 1 intra-well}
\begin{align}
    D_0^2 z_1 +\omega_o^2 z_1 &=-2D_0 D_1 z_0 - \eta z_0^2, \\
    D_0 v_1+\theta z_1 &= D_0z_1+  D_1z_0 - D_0v_0,
\end{align}
\end{subequations}\\
\underline{$O(\epsilon^2)$}:\\
\begin{subequations} \label{eq:order 2 intra-well}
\begin{align}
\begin{split}
    D_0^2 z_2 + \omega^2_o z_2 =&-2D_0 D_1 z_1 - 2D_0 D_2z_0 -D_1^2 z_0 -\delta_1 \int_0^{T_0} \overline{h} (T_0 - \tau) D_0 z_0 (\tau) d\tau \\
    &-\delta_2 D_0z_0 -2\eta z_0 z_1 - \gamma z^3_0+g_{wave} \cos\left((\omega_o+\epsilon^2 \sigma)T_0\right),
\end{split}\\
\begin{split}
D_0 v_2+ \theta  v_{2} =&  D_0 z_2 +  D_2 z_0 +  D_1z_1 - D_2 v_0 - D_1 v_1.
\end{split}
\end{align}
\end{subequations}\\
Upon solving Equations (\ref{eq:order 0 intra-well}) and (\ref{eq:order 1 intra-well}), we obtain the following expressions for $z_0$, $v_0$, $z_1$ and $v_1$:
\begin{subequations} \label{eq:y_0 intra-well}
\begin{align}
    z_0 &=A(T_1,T_2) e^{i \omega_o T_0} + cc,\\
    v_0 &=\Gamma_0 A(T_1,T_2) e^{i \omega_o T_0} + cc,
\end{align}
\end{subequations}
and
\begin{subequations} \label{eq:y_1 intra-well}
\begin{align}
   z_1 &=\frac{\eta}{\omega^2_o} \left( \frac{A^2(T_1,T_2)}{3} e^{2 i \omega_o T_0} - 2 A(T_1,T_2) \overline{A}(T_1,T_2) \right) + cc, \\
   v_1 &=\eta \Gamma_1 \frac{A^2(T_1,T_2)}{3 \omega_o^2} e^{2 i\omega_0 T_0} +cc,
\end{align}
\end{subequations}
where
\begin{equation*}
    \Gamma_0 = \frac{\omega^2_o + i \theta \omega_o }{\omega^2_o + \theta^2},  \qquad \Gamma_1 = \frac{4 \omega^2_o + 2i  \theta \omega_o }{4 \omega^2_o +\theta^2},
\end{equation*}
and $cc$ stands for the complex conjugate of the preceding terms. 

Elimination of the secular terms from the second-order perturbation problem, yields
 $D_1 A(T_1,T_2)=0$, which implies that the complex valued function $A$ is only dependent on the third time scale, $T_2$. 
 
Elimination of the secular terms from the third-order problem associated with Equation (\ref{eq:order 2 intra-well}) requires evaluating the convolution integral $ \int_0^{T_0} \overline{h}(T_0 - \tau) D_0 z_0 (\tau) d\tau$. To this end, we use Equation (\ref{eq:h(t) analytical}) to write:
\begin{equation} \label{eq:h(T0-tau)}
%\begin{split}
    \overline{h}(T_0 - \tau)= \frac{e^{-\mu T_0}}{M} \sqrt{\frac{R}{g}} \Big( \lambda_1 e^{\mu \tau} + \lambda_2 e^{\mu \tau} \cos(\mu \tau - \mu T_0) - \lambda_3 e^{\mu \tau} \sin(\mu \tau - \mu T_0) \Big).
%\end{split}
\end{equation}
Also, from Equation (\ref{eq:y_0 intra-well}), we have
\begin{equation} \label{eq:D_0 y_0}
    D_0 z_0 (\tau)= i \omega_o A e^{i \omega_o \tau } + cc.
\end{equation}
Thus, using Equations (\ref{eq:h(T0-tau)}) and (\ref{eq:D_0 y_0}), we can express the convolution integral as
\begin{equation} \label{eq:convolution expanded}
\begin{split}
    \int_0^{T_0} \overline{h}(T_0 - \tau) D_0z_0 (\tau) d\tau =i \omega_o A e^{-\mu T_0}\int_0^{T_0} \Big( &\lambda_1 e^{(\mu+i \omega_o) \tau}+ \lambda_2 e^{(\mu+i \omega_o) \tau} \cos(\mu \tau - \mu T_0) \\
    &- \lambda_3 e^{(\mu+i \omega_o) \tau} \sin(\mu \tau - \mu T_0) \Big) d\tau,    
\end{split}
\end{equation}
which upon integration by parts yields
\begin{equation} \label{eq:convolution solution -intrawell}
     \int_0^{T_0} \overline{h}(T_0 - \tau) D_0 z_0 (\tau) d\tau= A \omega_o (\xi_0 +i \bar{\xi}_0) e^{i \omega_o T_0} +NST+cc,
\end{equation}
where $NST$ stands for non-secular terms, and 
\begin{equation} \label{eq:C_1 -intrawell}
    \begin{split}
        \xi_0&=\left( \frac{\lambda_1\omega_o}{\mu^2+\omega_o^2} + \frac{2\lambda_2 \mu^2\omega_o}{4\mu^4+\omega_o^4} -\frac{\lambda_3\omega_o^3}{4\mu^4+\omega_o^4} \right),\\
        \bar{\xi}_0&=\left( \frac{\lambda_1 \mu}{\mu^2+\omega_o^2} + \frac{\lambda_2(2\mu^3-\mu\omega_o^2)}{4\mu^4+\omega_o^4}-\frac{\lambda_3(2\mu^3+\mu\omega_o^2)}{4\mu^4+\omega_o^4} \right).
    \end{split}
\end{equation}
Substituting Equations (\ref{eq:convolution solution -intrawell}) and (\ref{eq:C_1 -intrawell}) into the third perturbation problem presented in Equation (\ref{eq:order 2 intra-well}), and using the following polar transformation for the complex valued function $A(T_2)$
\begin{equation} \label{eq:polar transformation}
    \begin{split}
        A(T_2)&= \frac{a(T_2)}{2} e^{i \beta (T_2)}, \\
        \overline{A}(T_2)&= \frac{a(T_2)}{2} e^{-i \beta (T_2)},
    \end{split}
\end{equation}
then eliminating the secular terms, yields the following modulation equations:
\begin{equation} \label{eq:amplitude modulation-intra-well}
    D_2 a= -\left(\frac{\delta_1 \bar{\xi}_0}{2}  + \frac{\delta_2}{2}\right) a + \frac{g_{wave}}{2\omega_o} \sin \psi
\end{equation}
\begin{equation} \label{eq:phase modulation-intra-well}
        a D_2 \psi=\left(\sigma-\frac{\delta_1 \xi_0}{2}\right)a+\left(\frac{5\eta^2}{12\omega_o^3} - \frac{3\gamma}{8\omega_o}\right)a^3+ \frac{g_{wave}}{2\omega_o} \cos \psi,
\end{equation}
where $a$ and $\beta$ represent, respectively, the amplitude and phase of oscillations, and $\psi=\sigma T_2-\beta$.

 \subsection{Global inter-well response}
 In this section, we seek an approximate solution for the high energy orbits; i.e., when the system undergoes global inter-well oscillations. Since Equation (\ref{eq:ComminsWN}) has a negative linear stiffness,  it is difficult to employ the method of multiple scales in its straightforward fashion. To overcome this, we first expand the natural frequency of symmetric oscillations to be in the following form
  \begin{equation} \label{eq:Nat freq.}
     \omega_N^2=-\omega_n^2+ \sigma_1,
 \end{equation}
where $\sigma_1$ is a detuning parameter. Since we are seeking an approximate solution in the vicinity of the primary resonance of the inter-well motions, we express the nearness of the wave frequency to the natural frequency, $\omega_N$, by introducing the detuning parameter $\sigma_2$, such that
   \begin{equation} \label{eq:detuning interwell}
     \Omega^2=\omega_N^2+\epsilon^2 \sigma_2.
 \end{equation}
  Upon adding Equations (\ref{eq:Nat freq.} - \ref{eq:detuning interwell}), we arrive at
 \begin{equation} \label{eq:Nat freq 2}
     -\omega_n^2=\Omega^2 - (\omega_N^2+\omega_n^2) - \epsilon^2(\Omega^2-\omega_N^2).
 \end{equation}
Substituting Equation (\ref{eq:Nat freq 2}) into Equation (\ref{eq:ComminsWN}), we obtain 
 \begin{equation} \label{eq:Commins-interwell2}
 \begin{split}
      \ddot{Y} + &\epsilon^2 \delta_1 \int_{0}^{t^*} h(t^*-\tau) \dot{Y}(\tau) d\tau +  \epsilon^2 \delta_2 \dot{Y} +\Omega^2 Y \\
      + &\epsilon(-(\omega_N^2+\omega_n^2)Y+ \gamma Y^3)- \epsilon^2(\Omega^2-\omega_N^2)Y = \epsilon^2 g_{wave} \cos(\Omega t^*).
 \end{split}
 \end{equation}
 Next, we employ the method of multiple scales on Equation (\ref{eq:Commins-interwell2}) by introducing slow and fast time scales $T_n= \epsilon^n t^*$, $n=0,1,2$, and seek an expansion of the response in the form 

\begin{subequations} \label{eq:yexpansion}
\begin{align}
        Y(t^*, \epsilon)&=Y_0(T_0,T_1,T_2) + \epsilon Y_1(T_0,T_1,T_2) +\epsilon^2 Y_2(T_0,T_1,T_2) + O (\epsilon^3), \\
       v (t^*, \epsilon)&= v_0 (T_0,T_1,T_2) + \epsilon v_1(T_0,T_1,T_2) +\epsilon^2 v_2(T_0,T_1,T_2) + O (\epsilon^3).
\end{align}
\end{subequations}  
Implementing the method of multiple scales on the scaled equations as described in the previous subsection yields the following solution for $Y$ and $v$. 

\begin{subequations}
\begin{align}
     Y &= a \cos (\Omega t^* + \psi) + O (\epsilon), \\
     v &= \Gamma_2 a \cos (\Omega t^*+ \psi) + O (\epsilon),
\end{align}
\end{subequations}
where
\begin{equation*}
    \Gamma_2 = \frac{\Omega^2 + i \theta \Omega}{\Omega^2 + \theta^2}
\end{equation*}
and  the amplitude, $a$, and phase, $\psi$, of the response are governed by the following modulation equations:
\begin{equation} \label{eq:amplitude modulation-interwell}
     \omega_n D_2 a = -\left(\frac{\omega_n \delta_1  \bar{\xi_1} + \omega_n \delta_2}{2}\right)a - \frac{f_{wave}}{2} \sin \psi
\end{equation}
\begin{equation} \label{eq:phase modulation-interwell}
\omega_n a D_2 \psi=-\left(\frac{\Omega^2 + \omega_n^2}{2 } -\omega_n \delta_1 \xi_1 \right)a + \frac{3 \gamma}{8} a^3+\frac{3 \gamma^2 }{256 \omega_n^2}a^5-\frac{f_{wave}}{2} \cos \psi,
\end{equation}
where
\begin{equation} \label{eq:C_2 -interwell}
    \begin{split}
        \xi_1&=\left( \frac{\lambda_1 \omega_n}{\mu^2+\omega_n^2} + \frac{2\lambda_2 \mu^2 \omega_n}{4\mu^4+\omega_n^4} -\frac{\lambda_3 \omega_n^3}{4\mu^4+\omega_n^4} \right)\\
        \bar{\xi}_1&=\left( \frac{\lambda_1 \mu}{\mu^2+\omega_n^2} + \frac{\lambda_2(2\mu^3-\mu \omega_n^2)}{4\mu^4+\omega_n^4}-\frac{\lambda_3(2\mu^3+\mu \omega_n^2)}{4\mu^4+\omega_n^4} \right)
    \end{split}
\end{equation}

\subsection{Steady-state response}
The long-time steady-state behavior of the absorber is of particular interest to assess its performance. Thus, the  time derivatives in Equations (\ref{eq:amplitude modulation-intra-well} - \ref{eq:phase modulation-intra-well}) and  (\ref{eq:amplitude modulation-interwell} - \ref{eq:phase modulation-interwell}) are set to zero and the resulting algebraic equations are then solved numerically for the steady-state amplitude $a_o$, and phase $\psi_o$. For intra-well oscillations, we obtain: 
\begin{subequations} \label{eq: intrawell oscillations}
\begin{align}
\begin{split} \label{eq:Intra-Subeq1}
  Y=a_o \cos(\Omega t^*- \psi_o)+ \frac{\eta}{2 \omega_0^2} \left(-a_o^2+\frac{a_o^2}{3} \cos(2 \Omega t^*- 2 \psi_o)\right) +O(\epsilon^2),
\end{split}\\
 &v=\frac{\omega_0^2}{\omega_0^2+\theta^2}a_o \cos(\Omega t^* -\psi_o) - \frac{\omega_0 \theta}{\omega_0^2+\theta^2}a_o \sin(\Omega t^* -\psi_o) +O(\epsilon^2), 
 \label{eq:Intra-Subeq2}
\end{align}
\end{subequations}
and for the inter-well oscillations, we obtain
\begin{subequations} \label{eq: interwell oscillations}
\begin{align}
\begin{split} \label{eq:Inter-Subeq1}
Y&=a_o \cos(\Omega t^* - \psi_o) + \left( \frac{\gamma}{32 \omega_n^2} a_o^3 + \frac{3 \gamma^2}{1024 \omega_n^4} a_o^5 \right) \cos(3\Omega t^* - 3\psi_o) \\
&+ \frac{\gamma^2}{1024 \omega_n^4} a_o^5 \cos(5\Omega t^* - 5\psi_o) +O(\epsilon^2),
\end{split}\\
 v&=\frac{\omega_n^2}{\omega_n^2+\theta^2}a_o \cos(\Omega t^* -\psi_o) - \frac{\omega_n \theta}{\omega_n^2+\theta^2}a_o \sin(\Omega t^* -\psi_o) +O(\epsilon^2).
 \end{align}
\end{subequations}

Using the steady-state response,  the average power available at the buoy can be expressed as
\begin{equation} \label{eq: Normalized averaged power}
    P_{avg}=\frac{1}{T^*} \int_0^{T^*} \delta_2 \dot{Y}^2 dt^* ,
\end{equation}
where $T^*$ is the period of oscillations. For intra-well oscillations Equation (\ref{eq: Normalized averaged power}) reduces to
\begin{equation} \label{eq: POWER intrawell}
    P_{avg}= \delta_2 \left(\frac{\Omega^2 a_o^2}{2} + \frac{\eta^2 \Omega^2 a_o^4}{18 \omega_0^4}  \right) + O(\epsilon^2),
\end{equation}
while for inter-well oscillations, we get
\begin{equation} \label{eq: POWER interwell}
    P_{avg}=\delta_2 \left(\frac{\Omega^2 a_o^2}{2} + \frac{9 \gamma^2 \Omega^2 a_o^6}{2048 \omega_n^4}\right) +\ O(\epsilon^2).
\end{equation}

In addition to the averaged power, we are also interested in evaluating the capture width ratio (CWR), which is a common parameter used to evaluate a PWA's performance. CWR or absorption width as commonly coined in some literature is defined as the ratio between the average absorbed power, $P_{avg}$, and the power available at the wave front, $P_{wave}$. The latter can be obtained by multiplying the wave energy flux per unit crest length with the buoy's characteristic length, which is the diameter for a hemispherical buoy. This yields  \cite{wang2017modelling}:

\begin{equation} \label{eq: CWR}
    CWR=\frac{6 (m+m_\infty)  \Omega}{\rho R A^2} P_{avg}.
\end{equation}

\section{Stability Analysis and Numerical Simulations} \label{sec:bifurcations}
Bi-stable PWAs are known to produce large-amplitude responses over certain frequency ranges under harmonic waves excitations. Unfortunately, such desired motions can be incited only when the excitation is capable of channelling enough energy to the PWA to overcome the potential energy barrier and perform periodic inter-well motions.  Even when excited, those desired motions can only be uniquely realized over a specific bandwidth of the wave frequency, which we coin here as the effective bandwidth of the PWA. Outside the effective bandwidth, the large-amplitude motions are often accompanied with other, less desirable responses; e.g., small periodic or aperiodic motions. To define this effective bandwidth, it is important to study the stability of the steady-state periodic solutions $a_o$ and $\psi_o$ as the excitation frequency is varied. At the first level, this can be realized by finding the eigenvalues of the Jacobian matrices associated with Equations (\ref{eq:amplitude modulation-intra-well}) - (\ref{eq:phase modulation-intra-well}) and (\ref{eq:amplitude modulation-interwell})-(\ref{eq:phase modulation-interwell}). 

\begin{figure} [h!]
\centering
\includegraphics[width=0.75\textwidth]{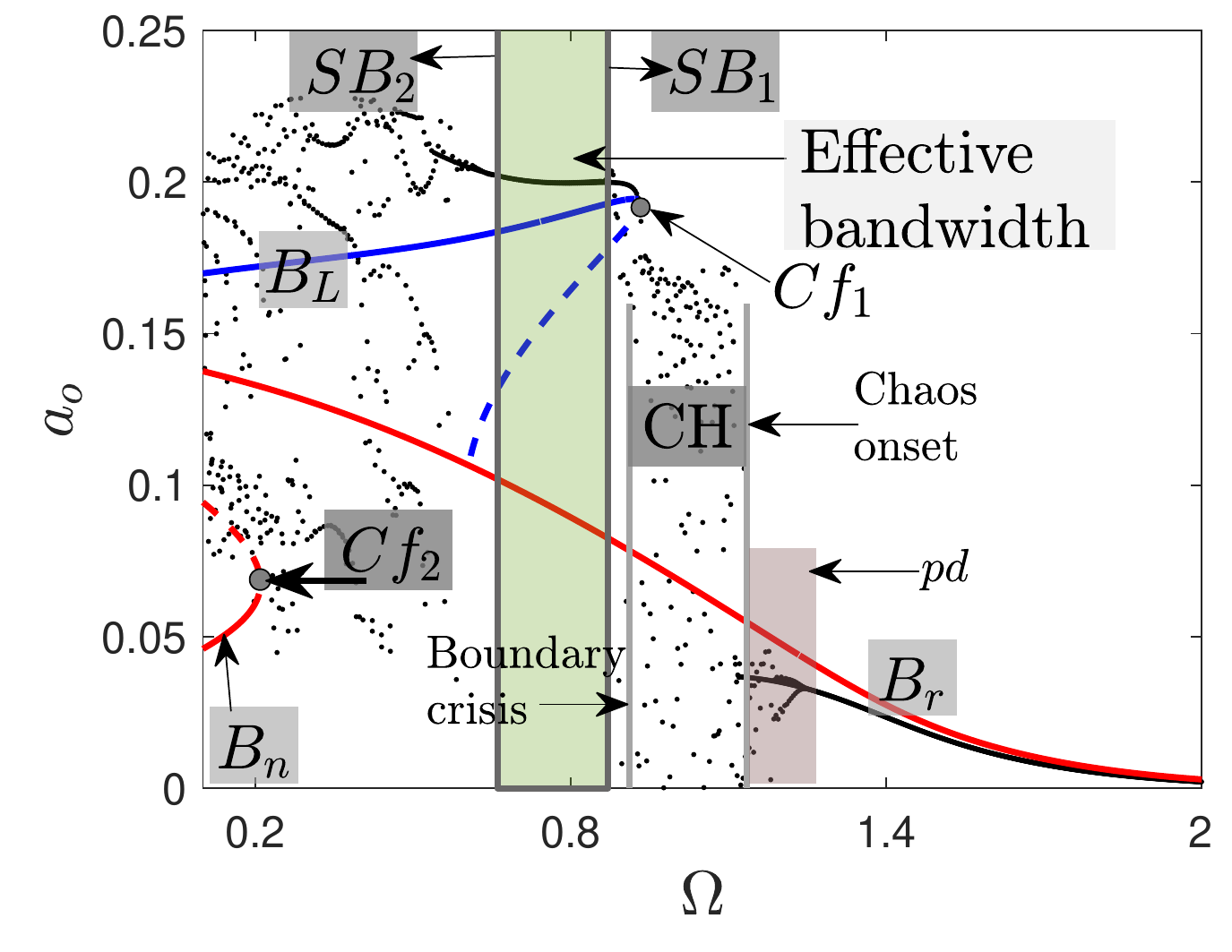}
\vspace{0.0cm}
\caption{Stroboscopic and analytical bifurcation diagrams for the bi-stable PWA under regular wave excitation at nondimensional wave
amplitude of 0.1. Solid lines: stable solution. Dashed lines: unstable solutions.}
\label{fig:JacobianStability}
\end{figure}

\begin{figure} [h!]
\centering
\includegraphics[width=0.75\textwidth]{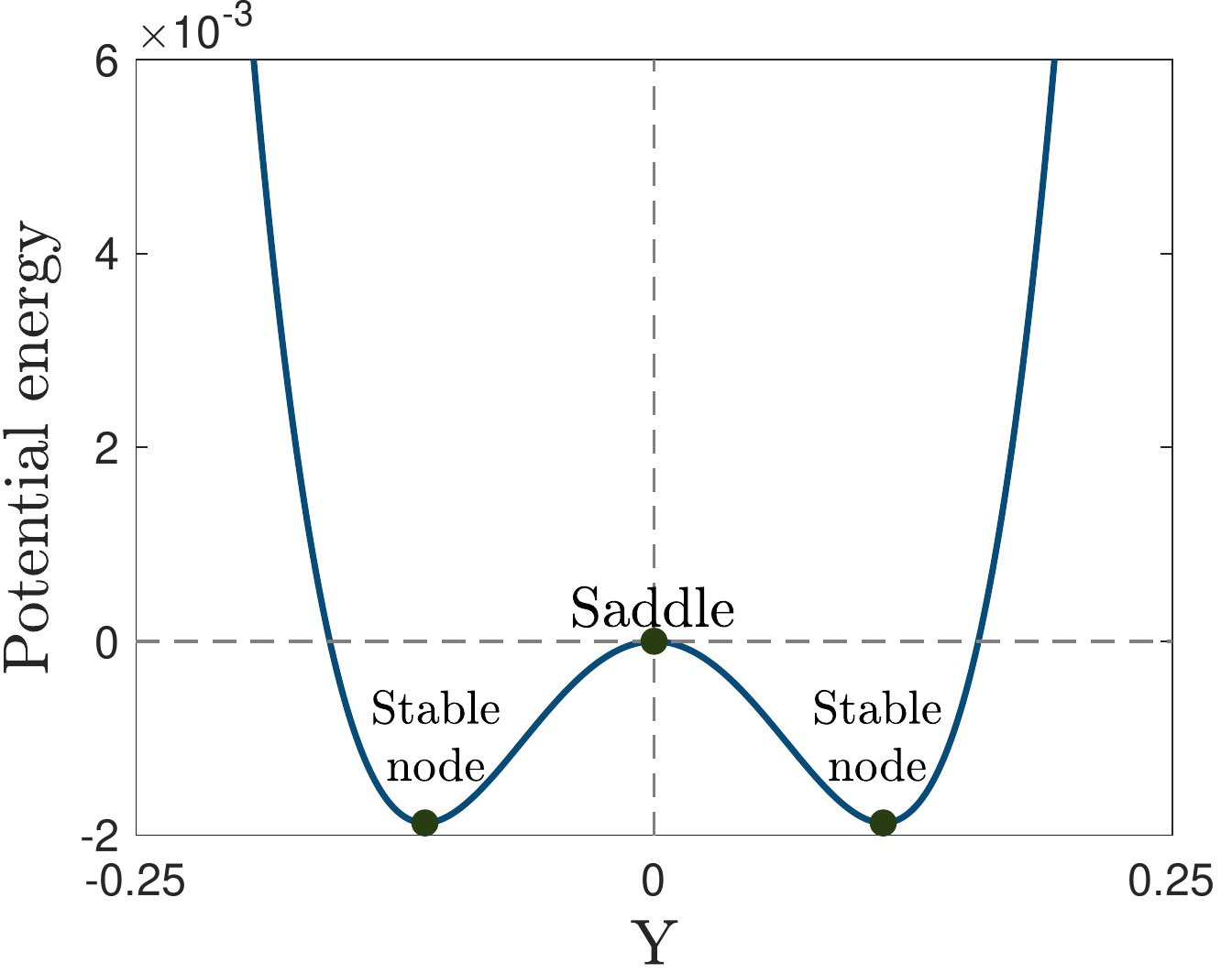}
\vspace{0.0cm}
\caption{Potential energy function associated with the system's parameters listed in Table \ref{Table:WEC parameters}.}
\label{fig:POT2}
\end{figure}

Figure \ref{fig:JacobianStability} depicts variation of the steady-state amplitude of the absorber with the excitation frequency as obtained using equations (\ref{eq:amplitude modulation-intra-well}) - (\ref{eq:phase modulation-intra-well}) and (\ref{eq:amplitude modulation-interwell})-(\ref{eq:phase modulation-interwell}). The parameters used in the simulation are listed in Table \ref{Table:WEC parameters} and the associated symmetric potential energy function of the PWA is depicted in Figure \ref{fig:POT2}.

\begin{table}
\caption{WEC parameters.}
\label{Table:WEC parameters}
\begin{center}
\begin{tabular}{c  c} 
\hline
 $Parameter$ & $Value$  \\ 
 \hline
 $\omega_n$ & 0.78\\  
 $\gamma$ & 50 \\    
 $\delta_2$ & 0.13  \\
 \hline
\end{tabular} 
\end{center} 
\end{table}

On the figure, the solid lines represent stable steady-state single-period periodic solutions of amplitude $a_o$, and dashed lines represent unstable unrealizable periodic solutions. It is evident that there are three branches of stable periodic solutions based on the Jacobian-based stability analysis. The branches $B_r$ and $B_n$ which represent, respectively, the resonant and non-resonant branches of the small-amplitude intra-well motions, and the branch $B_L$, which represents the large amplitude inter-well motions. The stability analysis reveals two cyclic-fold bifurcations, $Cf_1$, and $Cf_2$, which result from the stable and unstable periodic orbits colliding and destructing each other.

\begin{figure}  [h!]
\centering
  \begin{tabular}{@{}c|c|c@{}}
  $(a)$ &
  $(b)$&
  $(c)$\\
    \includegraphics[width=0.31\textwidth]{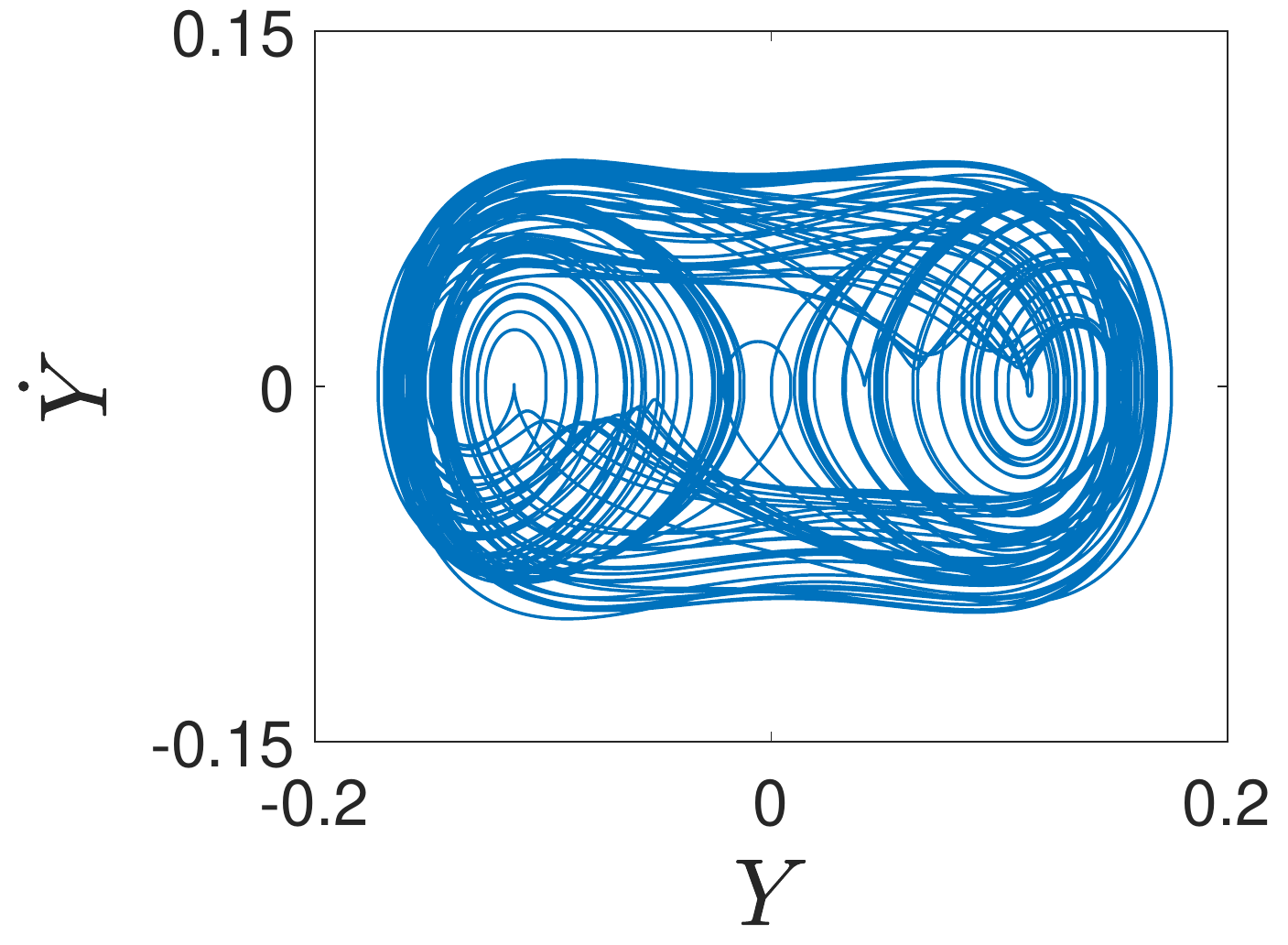} & 
    \includegraphics[width=0.31\textwidth]{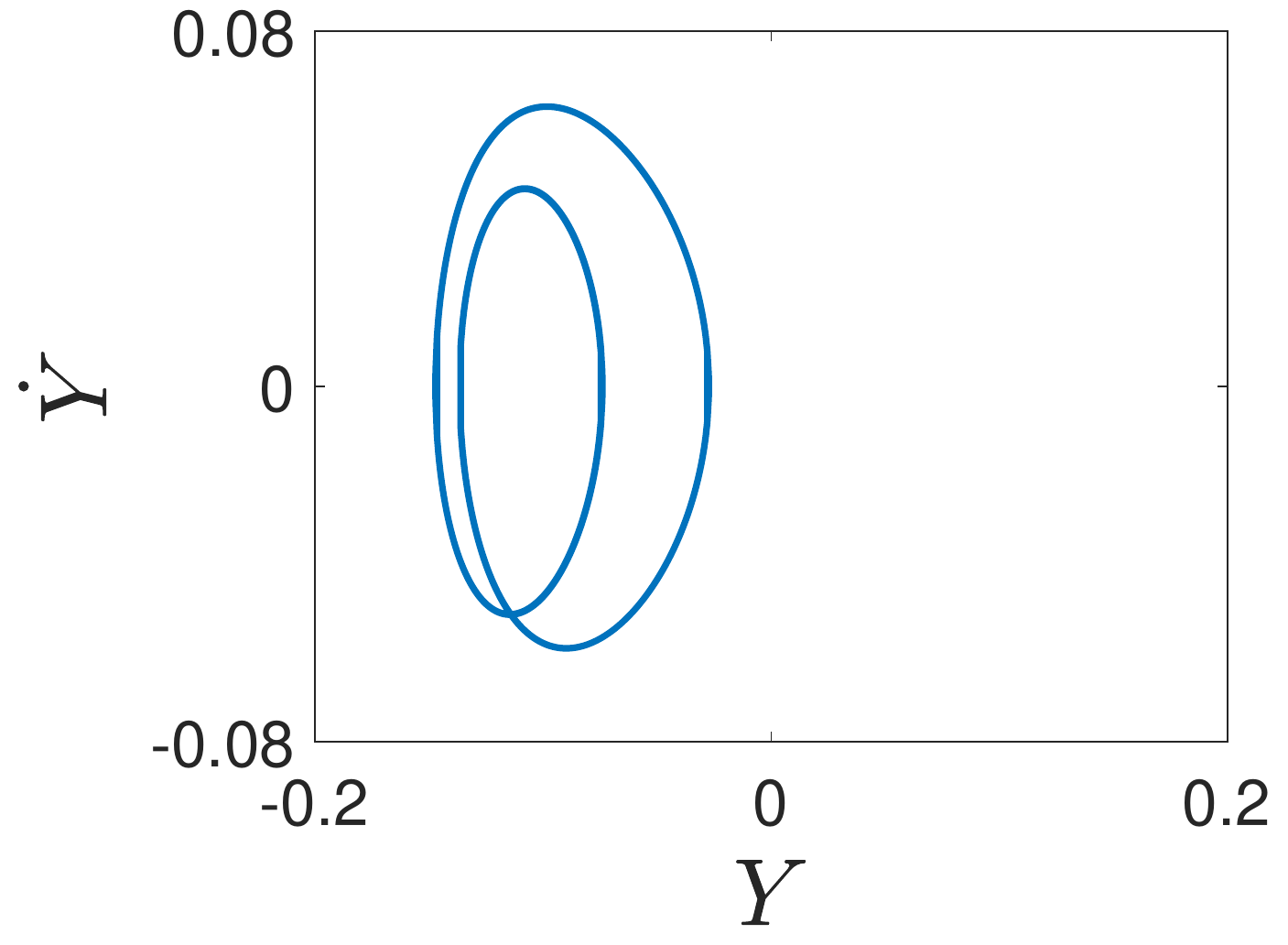}&
    \includegraphics[width=0.31\textwidth]{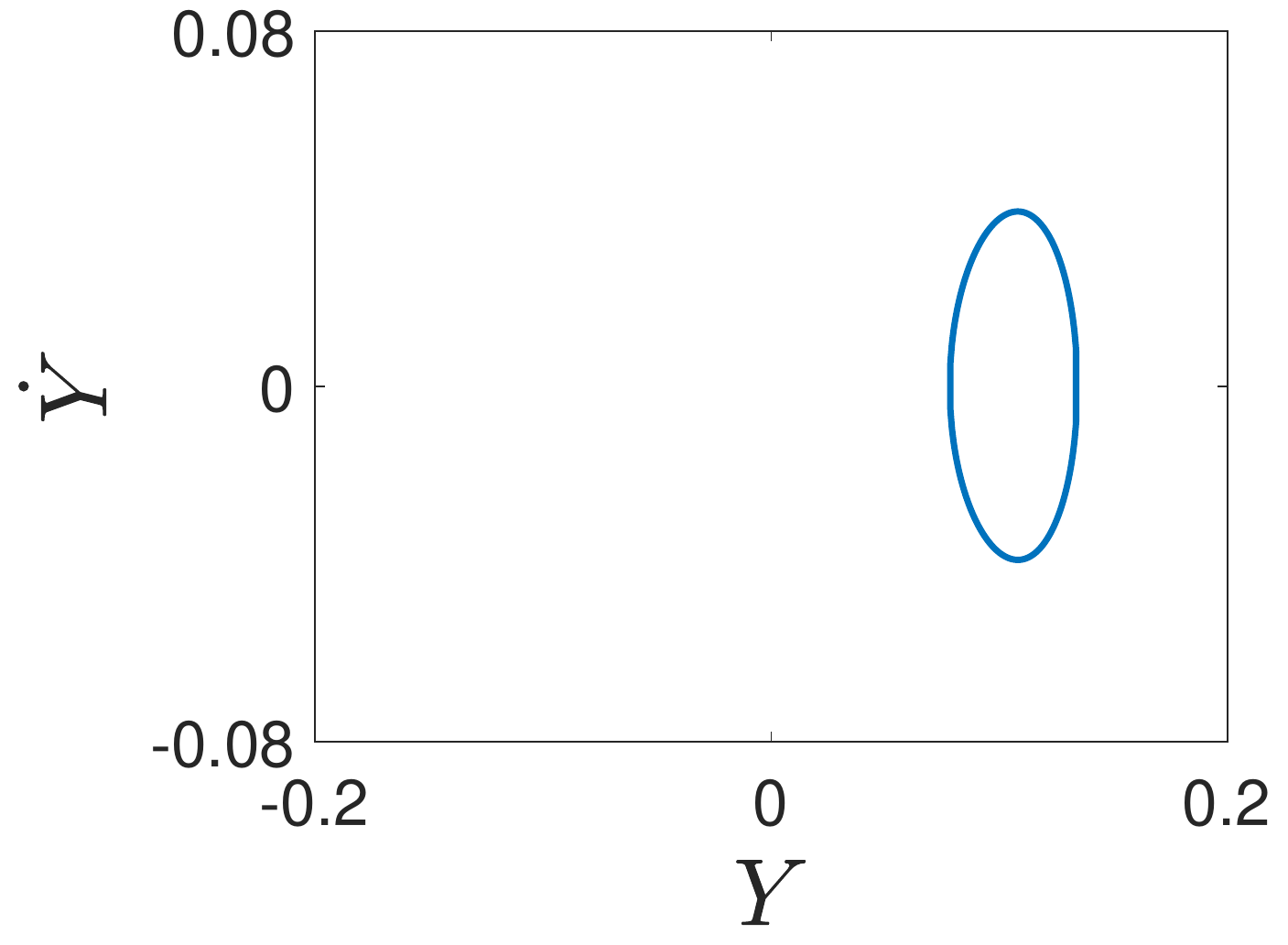} \\
    \includegraphics[width=0.31\textwidth]{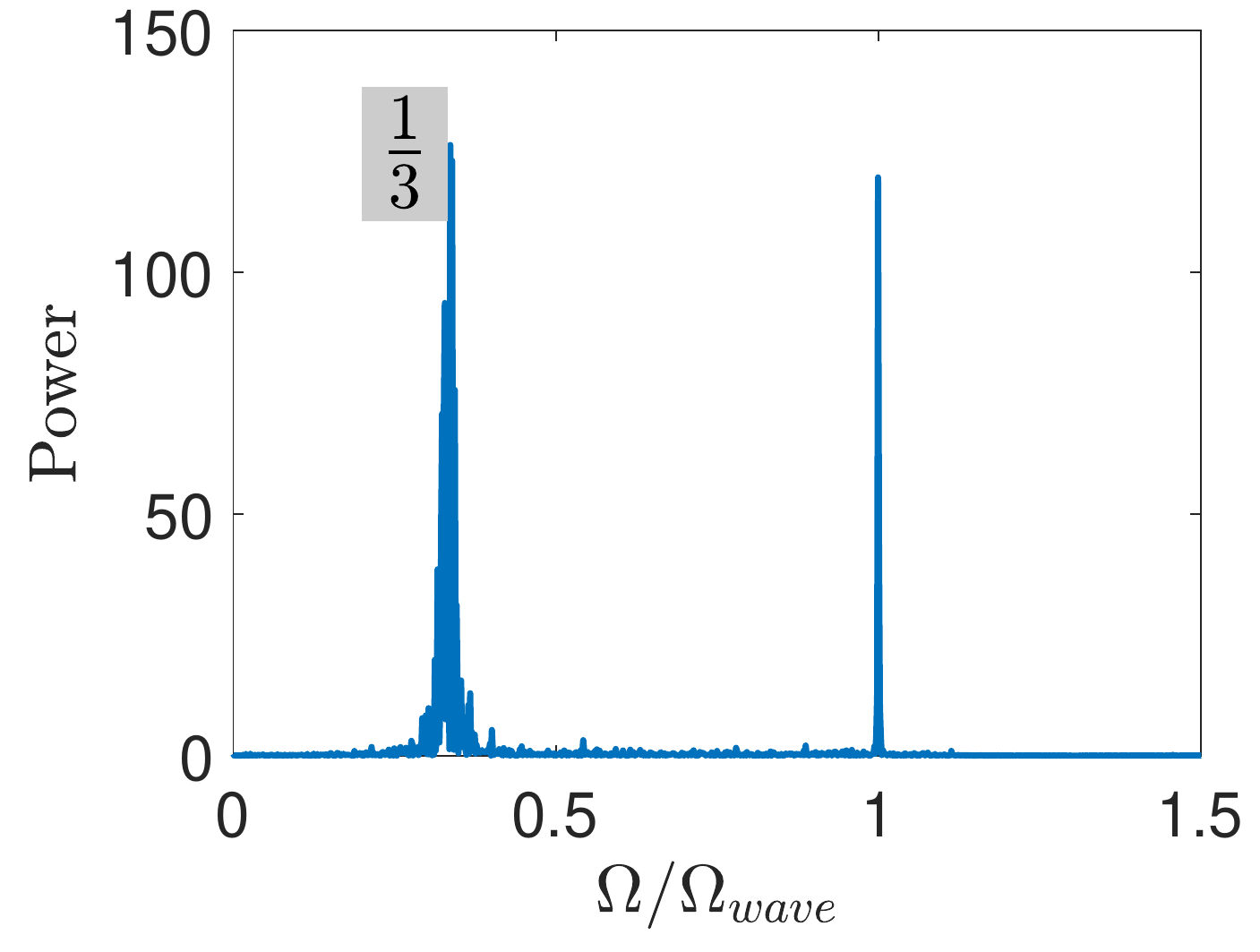} &
    \includegraphics[width=0.31\textwidth]{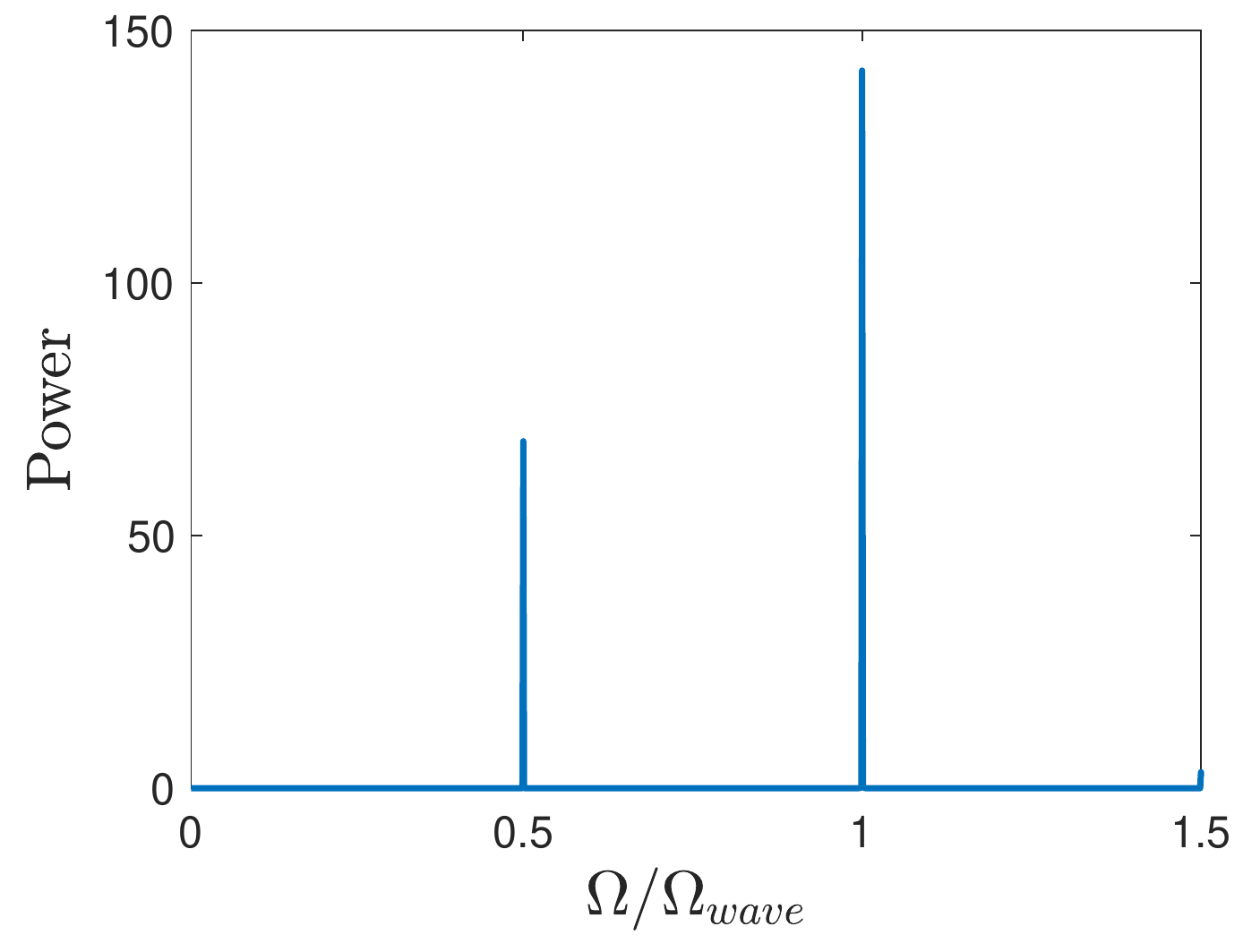} &
    \includegraphics[width=0.31\textwidth]{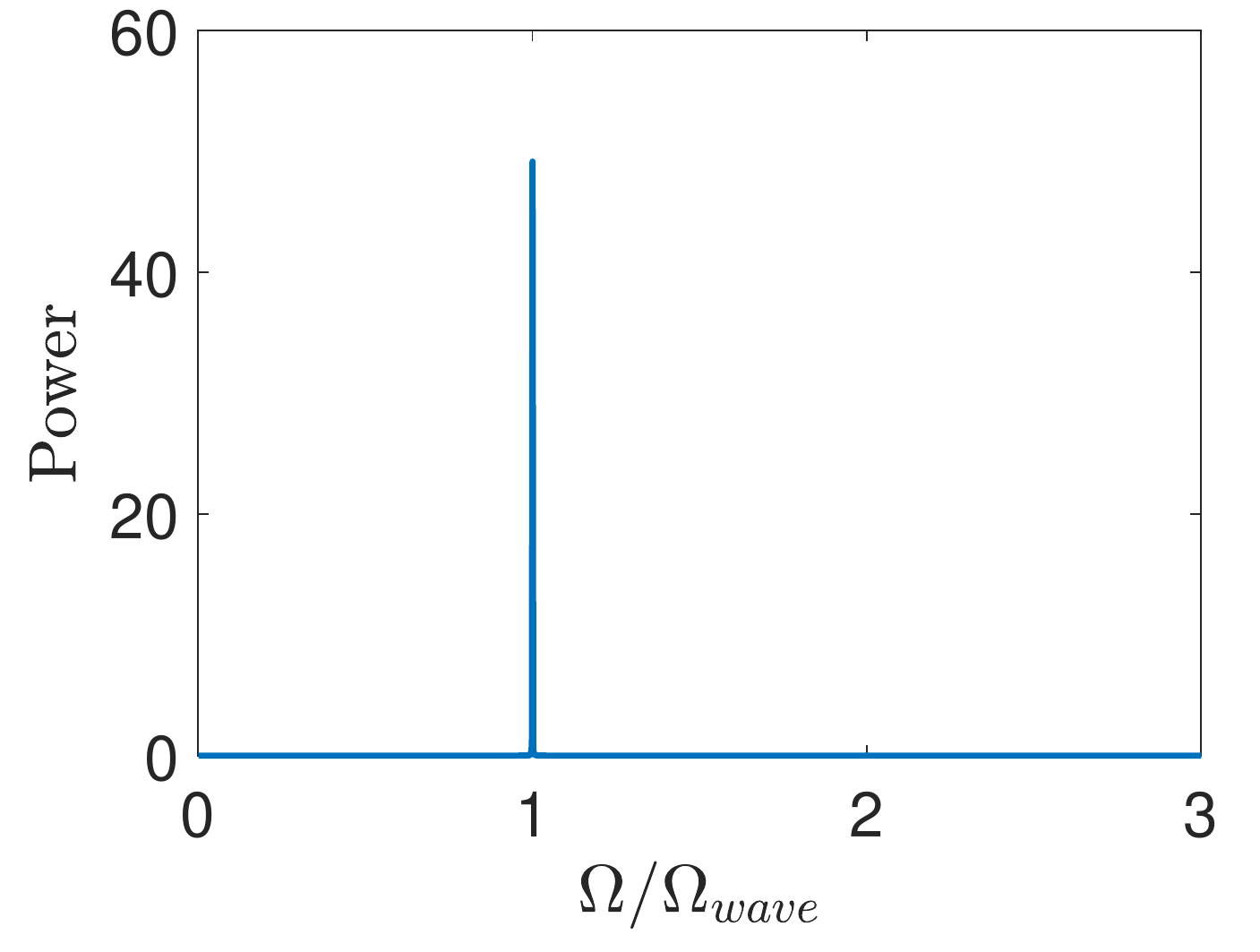}
      \end{tabular}
   \caption{Phase portraits and FFT spectra showing the instability route of the intra-well branch. $(a):\Omega_{wave}=1$, $(b):\Omega_{wave}=1.2$, and $(c):\Omega_{wave}=1.4$. }
  \label{fig: FFTs_Local}
\end{figure}

The superimposed stroboscopic bifurcation map in Figure \ref{fig:JacobianStability} obtained by a numerical integration of the original equations of motion, Equation (\ref{eq:Commins}), reveals a much more complex behaviour than what can be seen by relying on the Jacobian-based stability analysis. In particular, The bifurcation map reveals regions of aperiodic motions that extend over a wide range of frequencies. While we notice good agreement between the analytical solution and the numerical solution on the $B_r$ branch in the higher range of frequencies down to about $\Omega=1.2$. The Jacobian-based stability analysis does not reveal that the intra-well solutions on the branch, $B_r$, actually undergo a cascade of period-doubling, $pd$, bifurcations starting near $\Omega=1.2$. These bifurcations ultimately lead to a window of chaos, $CH$, which extends down to about $\Omega=0.9$. This window ultimately disappears in a boundary crisis. Phase portraits and the Fast Fourier Transform (FFT) of the period-doubled and chaotic solutions are shown in Figure \ref{fig: FFTs_Local} clearly demonstrating the period-doubling route to chaos.

\begin{figure}  [h!]
\centering
  \begin{tabular}{@{}c|c|c@{}}
  $(a)$ &
  $(b)$&
  $(c)$\\
    \includegraphics[width=0.31\textwidth]{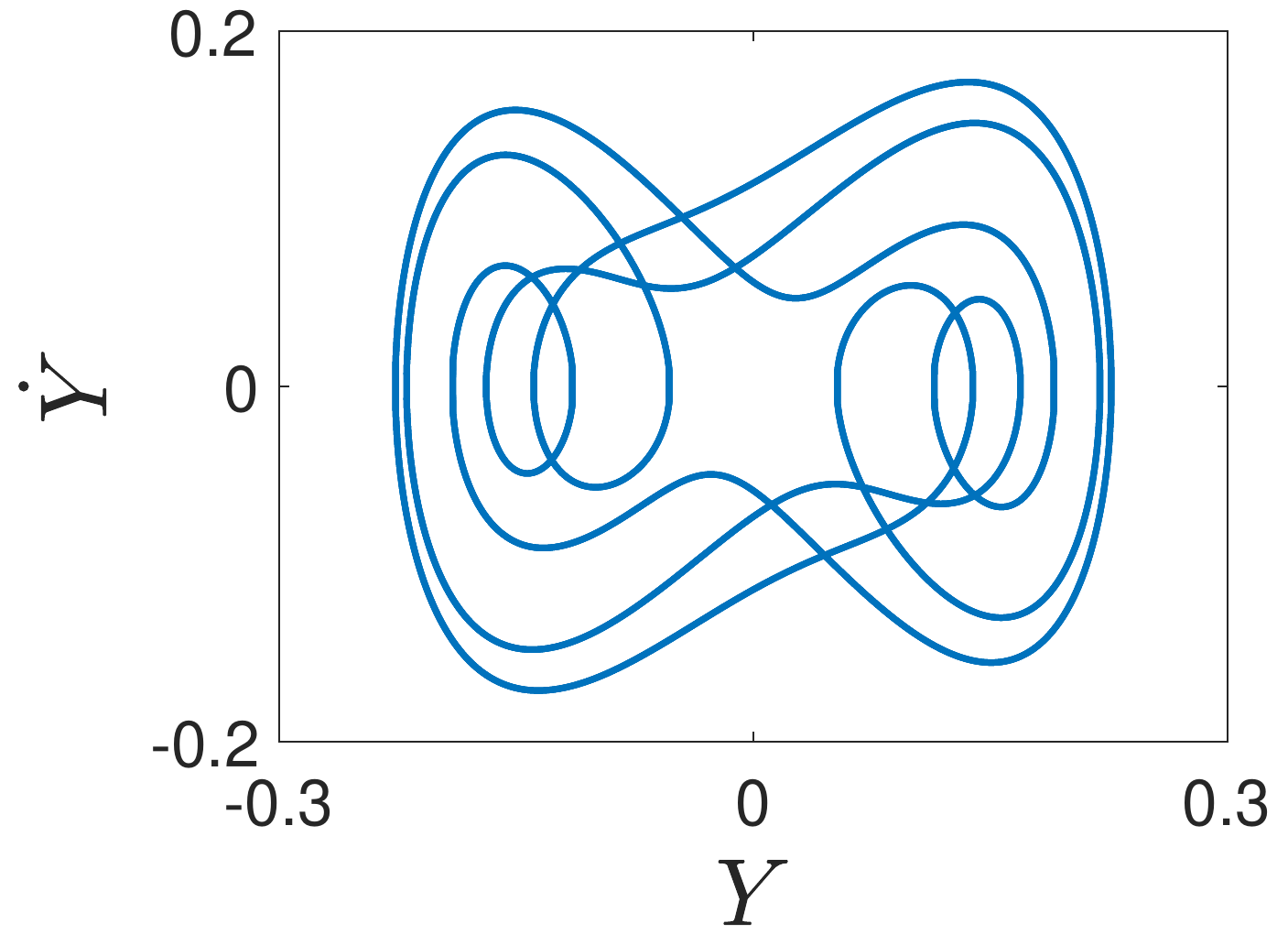} & 
    \includegraphics[width=0.31\textwidth]{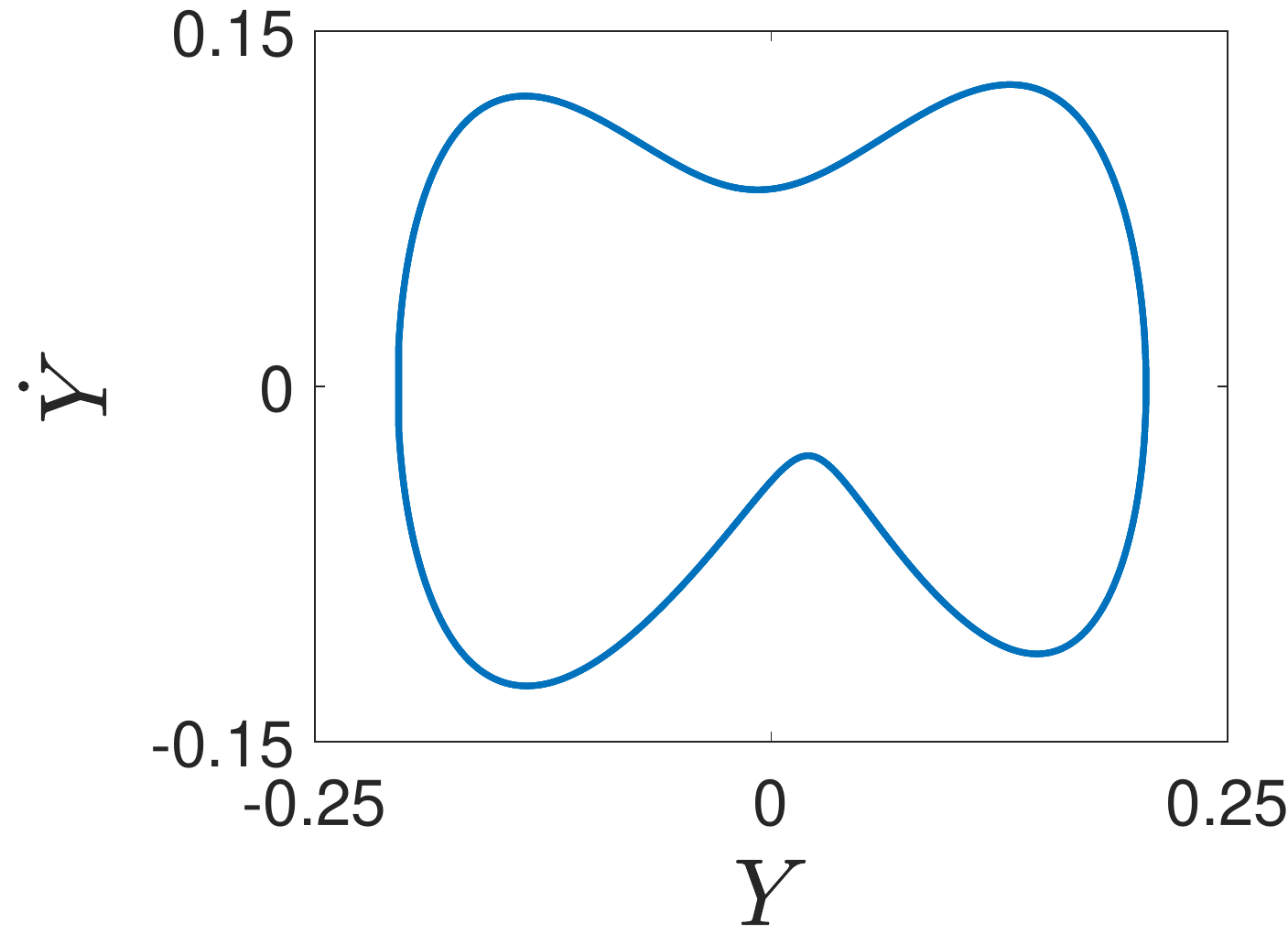}&
    \includegraphics[width=0.31\textwidth]{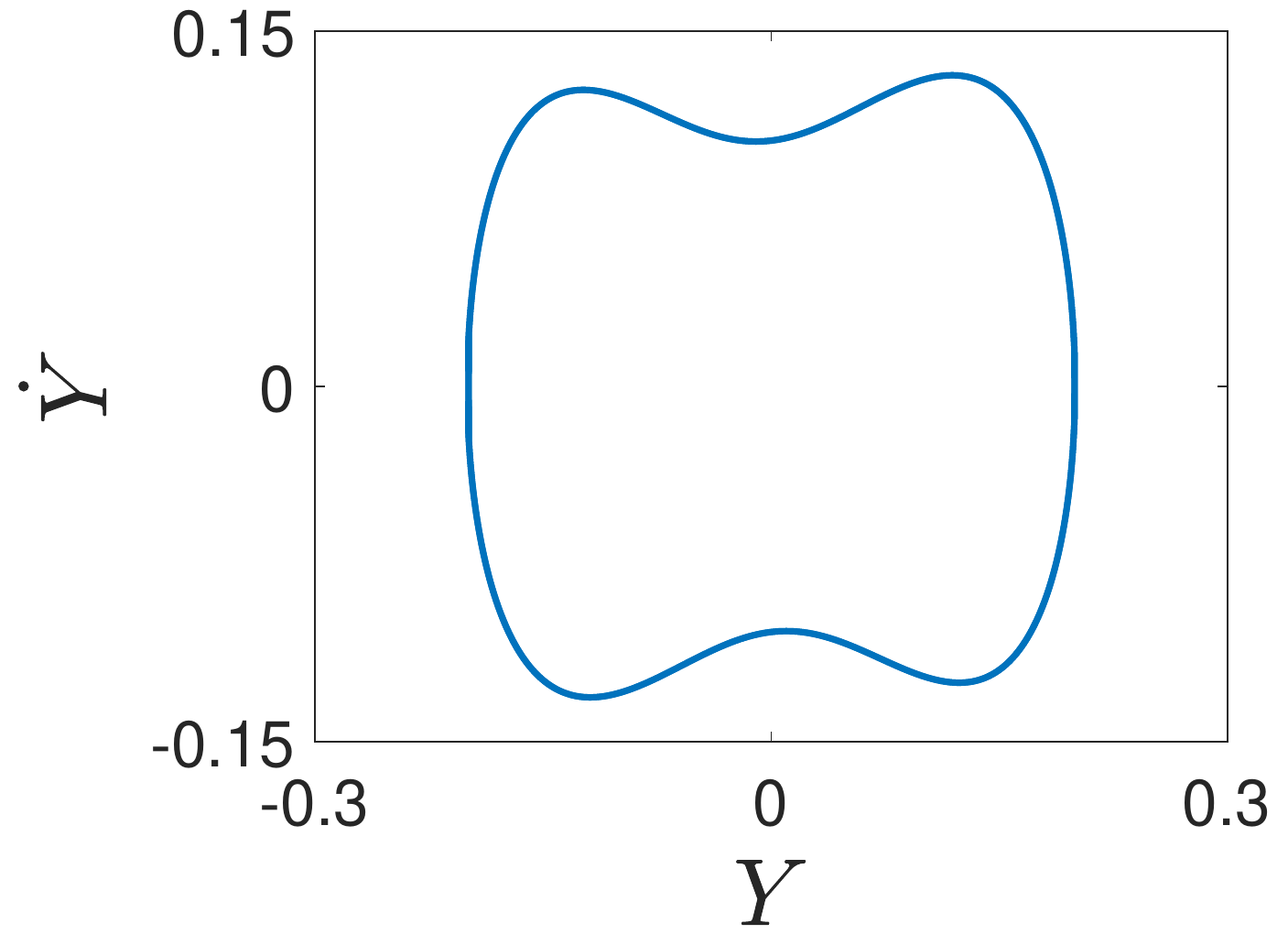} \\
    \includegraphics[width=0.31\textwidth]{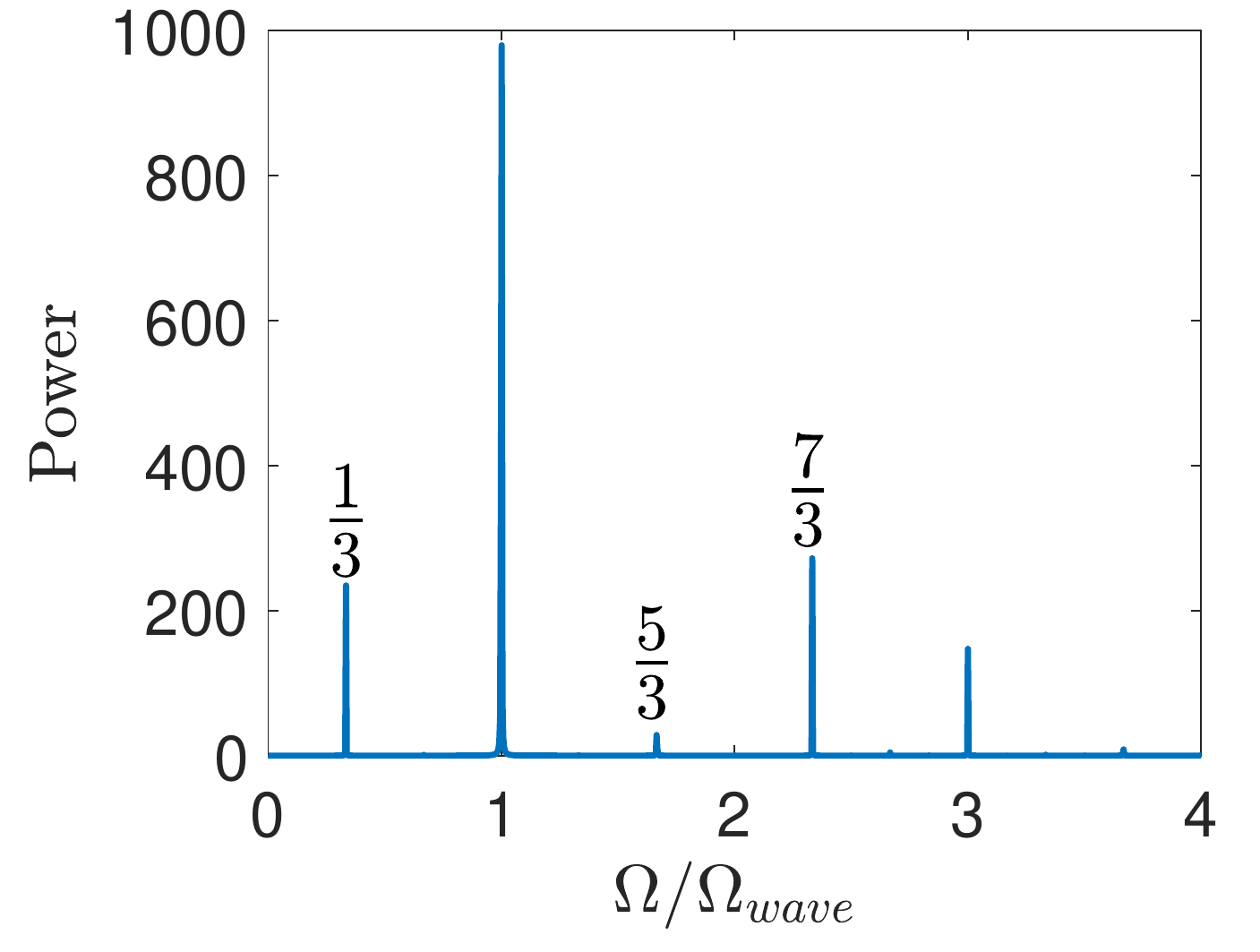} &
    \includegraphics[width=0.31\textwidth]{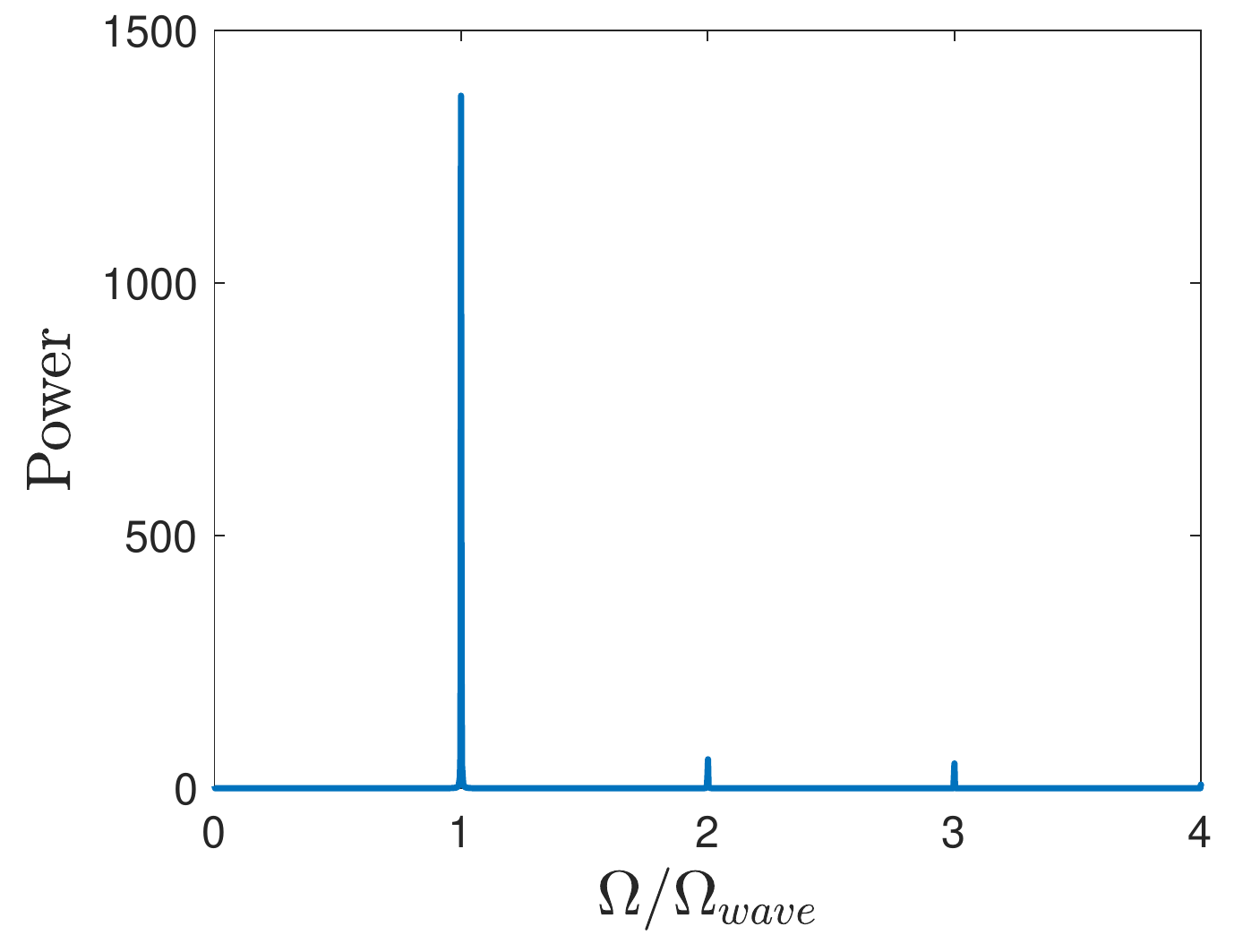} &
    \includegraphics[width=0.31\textwidth]{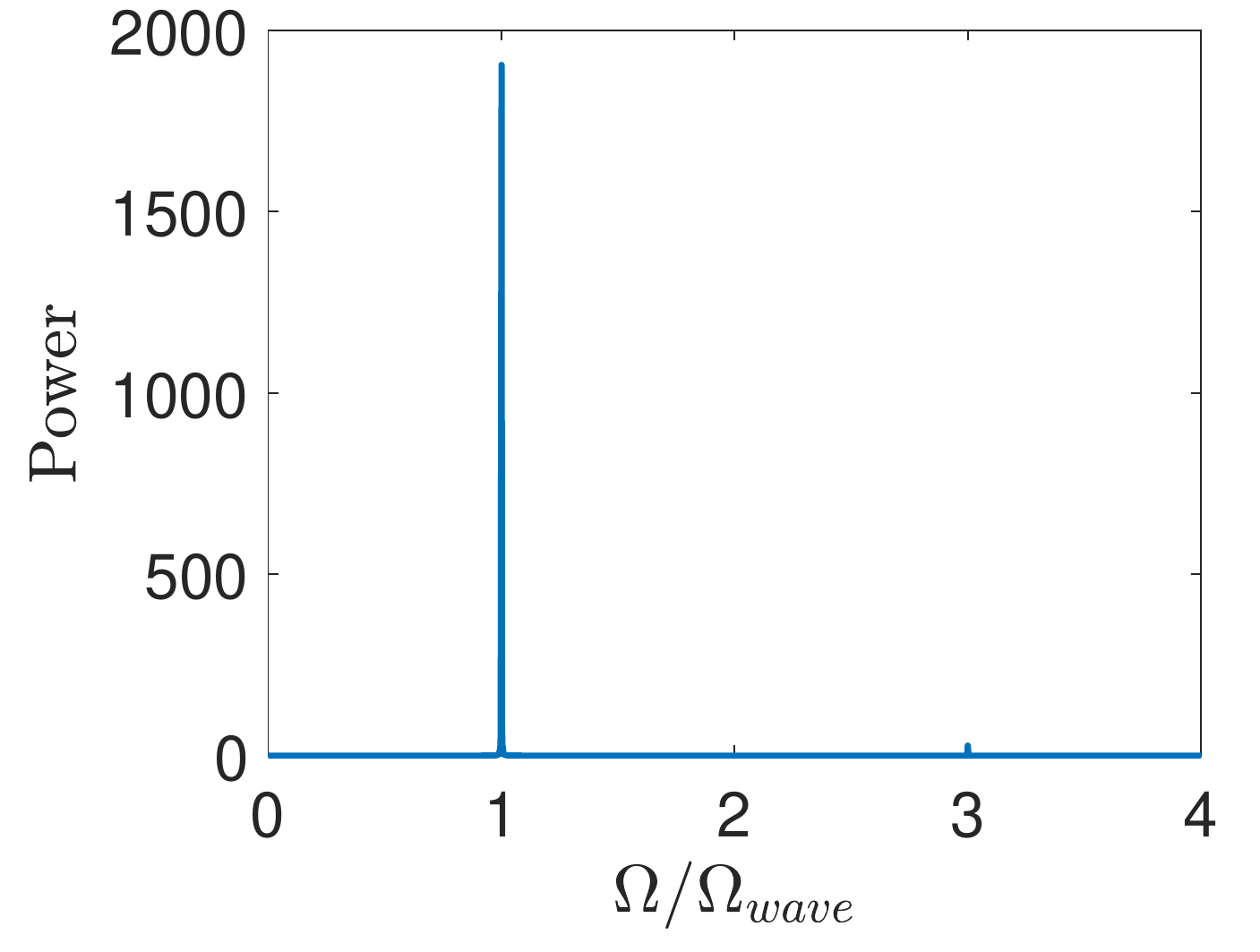}
      \end{tabular}
   \caption{Phase portraits and FFT spectra showing the instability route of the inter-well branch. $(a):\Omega_{wave}=0.42$, $(b):\Omega_{wave}=0.62$, and $(c):\Omega_{wave}=0.8$. }
  \label{fig: FFTs_Global}
\end{figure}

Furthermore, as shown in Figure \ref{fig:JacobianStability}, the period one solutions on branch $B_L$,  also lose stability via symmetry breaking bifurcations, $SB_1$ and $SB_2$ resulting in asymmetric periodic orbits as shown in the phase portraits and FFTs of Figure \ref{fig: FFTs_Global} . In the region between $SB_1$ and $SB_2$, the desired large orbit period one periodic solution is unique and, unlike the rest of the frequency bandwidth considered, does not coexist with other less desirable solutions. We coin this desired bandwidth as the effective bandwidth of the PWA. 

\begin{figure} [h!]
\centering
\includegraphics[width=0.75\textwidth]{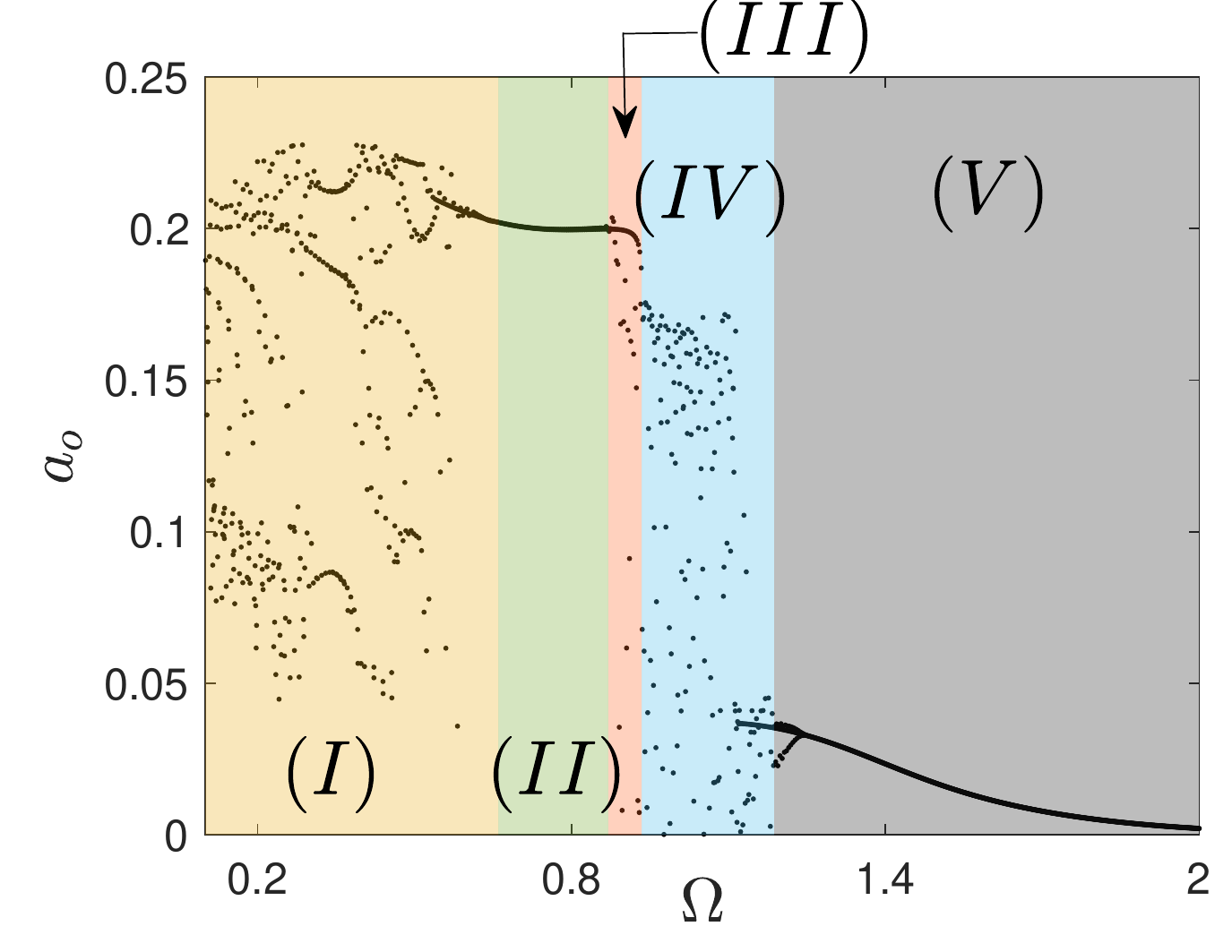}
\vspace{0.0cm}
\caption{Stroboscopic diagram for the Bi-stable PWA showing regions with different types of motion. Simulation performed at nondimensional wave
amplitude $A_{wave}/R$ = 0.1. }
\label{fig:motion_discription_NUM}
\end{figure}

Based on the numerical analysis, the response bandwidth of the absorber can be divided into different regions as shown in Figure \ref{fig:motion_discription_NUM}. The first region ($I$) occurs in the low frequency range and contains inter-well chaotic motions coexisting with small magnitude intra-well motions. The second region ($II$)  is the effective bandwidth, which contains the unique large-magnitude period one inter-well motions. The third region ($III$) contains chaotic motions coexisting with asymmetric periodic motions. Region four ($IV$) contains chaotic motions only. Finally, region five ($V$) contains unique small-amplitude period one intra-well periodic motions.

\begin{figure} [h!]
\centering
\includegraphics[width=0.75\textwidth]{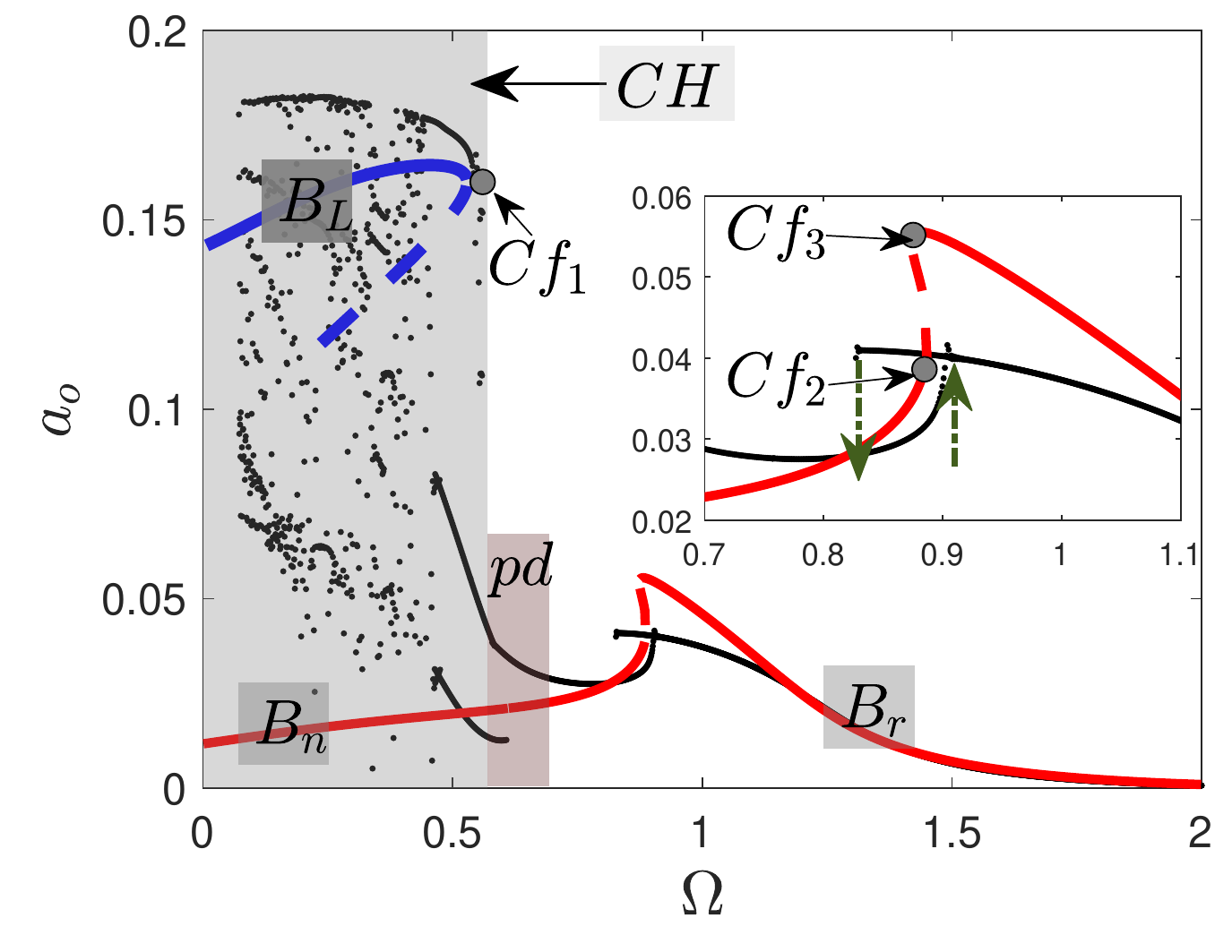}
\vspace{0.0cm}
\caption{Stroboscopic and analytical bifurcation diagrams for the bi-stable PWA under regular wave excitation at nondimensional wave
amplitude of 0.034. Solid lines: stable solution. Dashed lines: unstable solutions.}
\label{fig:JacobianStability2}
\end{figure}

At lower wave amplitudes, the bi-stable PWA reveals a different frequency response behavior as shown in Figure \ref{fig:JacobianStability2}. The most important feature is the birth of a new cyclic-fold bifurcation point denoted by $Cf_3$ and located at the tip of the resonant branch, $B_r$. In addition, we notice that the entire $B_r$ branch remains stable and the $pd$ bifurcation is shifted towards the nonresonant intra-well branch, $B_n$. This arrangement creates a region where two stable branches coexist, $B_r$ and $B_n$, which leads to jumps between them as depicted in the enlarged window in Figure \ref{fig:JacobianStability2}. In addition, the stroboscopic bifurcation map points shows a chaotic region that spans the entire frequency spectrum to the left of $Cf_1$. This Implies that the large-amplitude branch $B_L$ breaks symmetry at $Cf_1$, and no periodic inter-well motion can be realized.

Understanding the behavior of the absorber requires characterizing these regions as function of the frequency and amplitude of the incident waves. This can be realized by approximating the loci of the different bifurcations:  $Cf_1$, $Cf_2$, $Cf_3$ $pd$, $SB_1$ and $SB_2$. The first three bifurcations can be approximated using the Jacobian stability analysis, while the last three need the implementation of the Floquet theory.  In what follows, we obtain analytical approximation of all of these bifurcations as function of the design parameters of the PWA.

\subsection{Cyclic-fold bifurcation $Cf_1$ on the inter-well branch}
The cyclic-fold point $C_{f1}$ satisfies the relation $d\Omega/d a_o= 0$ on the resonant branch, $B_L$. Thus, its location can be obtained by differentiating Equations (\ref{eq:amplitude modulation-interwell}-\ref{eq:phase modulation-interwell}) with respect to $a_o$ and setting $d\Omega/d a_o$ to zero to arrive at the following 8$^{th}$ order polynomial:
\begin{equation} \label{eq: CF1}
    \begin{aligned}
        &\frac{45 \gamma^4}{4096 \omega_n^4}a_b^8 +\frac{9 \gamma^3}{128 \omega_n^2}a_b^6 + \left( \frac{99 \gamma^2}{64} - \frac{3 \gamma^2 \Omega_b^2}{128 \omega_n^2} + \frac{3 \gamma^2 \delta_1 \xi_1}{64 \omega_n} \right)a_b^4 -3\gamma \left(\Omega_b^2+\omega_n^2 -\delta_1 \omega_n \xi_1 \right)a_b^2\\
        &+\left(\Omega_b^2+\omega_n^2 -\delta_1 \omega_n \xi_1 \right)^2+\left(\delta_2 \omega_n + \delta_1 \omega_n\bar{\xi}_1\} \right)^2=0.
    \end{aligned}
\end{equation}
where $(\Omega_b, a_b)$ represent the wave frequency and response amplitude at which the bifurcation occurs. For any set of design parameters, the polynomial in Equation (\ref{eq: CF1}) can be solved to obtain the locus of the $Cf_1$ bifurcation in the parameter space of wave amplitude and frequency.
\subsection{Cyclic-fold bifurcations $Cf_2$ and $Cf_3$ on the intra-well branch}
These cyclic fold bifurcations also satisfy the relation $d \Omega/ d a_o = 0$ but on the intra-well branch, $B_n$; i.e. Equations (\ref{eq:amplitude modulation-intra-well}-\ref{eq:phase modulation-intra-well}). Differentiating Equations  (\ref{eq:amplitude modulation-intra-well}-\ref{eq:phase modulation-intra-well}) with respect to $a_o$ and setting $d \Omega/ d a_o = 0$, we obtain
\begin{equation} \label{eq: intra-well CF loci}
    a_b^2=\frac{-2\left(\Omega_b - \omega_o-\frac{\delta_1 \xi_0}{2}\right) \pm \sqrt{\left(\Omega_b - \omega_o-\frac{\delta_1 \xi_0}{2}\right)^2-3\left(\frac{\delta_1\bar{\xi}_0}{2}  + \frac{\delta_2}{2}\right)^2}}{3\left(\frac{5\eta^2}{12\omega_o^3} - \frac{3\gamma}{8\omega_o}\right)}.
\end{equation}
Here, $a_b$ and $\Omega_b$ are, respectively, the amplitude of oscillations and the wave frequency at which the cyclic-fold bifurcations, $Cf_2$ and $Cf_3$, occurs. Upon solving Equation (\ref{eq: intra-well CF loci}), we can obtain the loci of the $Cf_2$ and $Cf_3$ bifurcations in the parameter space of wave amplitude and frequency.

\subsection{Period doubling bifurcation} \label{sec: PD bifurcation}
One of the main dynamical characteristics of the intra-well oscillations is the appearance of a period-doubling bifurcation on the $B_r$ branch as the wave frequency is decreased toward lower values. To estimate this bifurcation point as function of the wave frequency and amplitude, we perturb the periodic orbit  ${Y}(t^*)$ by introducing the perturbation $\mathcal{p(t^*)}$ such that
\begin{equation} \label{eq:PD perturbation}
    \Tilde{Y}(t^*)=Y(t^*) + \mathcal{p(t^*)}.
\end{equation}
Here, $\Tilde{Y}(t^*)$ is a stable periodic orbit if $\mathcal{p(t^*) \to 0}$ as $t^* \to \infty$; otherwise it is unstable. Substituting Equation (\ref{eq:PD perturbation}) into Equation (\ref{eq:Commins intra-well}) and retaining only linear terms in $\mathcal{p(t^*)}$, we obtain, after simplifications, the following scaled equations:
 \begin{equation} \label{eq:Hill's}
 \begin{split}
      \mathcal{p}''+\epsilon \delta_1 \int_{0}^{t^*} \overline{h}(t^*-\tau) \mathcal{p}'(\tau) d\tau +\epsilon \delta_2 \mathcal{p}' &+ [ G_0+\epsilon G_1 \cos(\Omega t^* -\psi)+\epsilon^2 G_2 \cos(2\Omega t^* -2\psi) \\ 
      & +\epsilon^2G_3 \cos(3\Omega t^* -3 \psi) +\epsilon^2 G_4 \cos(4\Omega t^* -4\psi) ] \mathcal{p} =0,
 \end{split}
 \end{equation}
 where Equation (\ref{eq:Hill's}) represents a Mathieu’s-type differential equation with four parametric excitation terms. Each of these terms produces a principle subharmonic parametric instability when their frequency is half the natural frequency $\sqrt{G_0}$ of the perturbation  \cite{meyers2011mathematics, kovacic2018mathieu}. The interested reader can refer to Appendix \ref{Apndx: Parametric constants} for the full expressions of the parametric terms constants, $G_i$. 
 
 In order to examine the dynamics of the perturbation $\mathcal{p(t^*)}$, we derive an asymptotic approximation for the evolution of $\mathcal{p(t^*)}$ using the method of multiple scales assuming the following first-order expansion:
 \begin{equation} \label{eq: P expansion}
     \mathcal{p}(t^*,\epsilon)=\mathcal{p}_0(t^*)+\epsilon \mathcal{p}_1(t^*)+ O(\epsilon^2).
 \end{equation}
 Using the time scales and time derivatives defined in Equation (\ref{eq:time scales_intra}), and substituting Equation (\ref{eq: P expansion}) into Equation (\ref{eq:Hill's}), we obtain the following differential equations at the different orders of $\epsilon$:\\
\underline{$O(\epsilon^0)$}:\\
\begin{equation} \label{eq: order0 PD}
    D_0^2 \mathcal{p}_0 + G_0 \mathcal{p}_0=0,
\end{equation}
which admits the following homogeneous solution
\begin{equation} \label{eq: order0 soln PD}
    \mathcal{p}_0=Q(T_1) e^{i\sqrt{G_0}T_0}+\overline{Q}(T_1) e^{-i\sqrt{G_0}T_0},
\end{equation}
and \\
\underline{$O(\epsilon^1)$}:
\begin{equation} \label{eq: order1  PD}
%\begin{split}
        D_0^2 \mathcal{p}_1 +G_0 \mathcal{p}_1 = -2D_1 D_0 \mathcal{p}_0 - \delta_1 \int_0^{T_0} \overline{h}(T_0-\tau)D_0 \mathcal{p}_0 d\tau -\delta_2 D_0 \mathcal{p}_0 -G_1 \mathcal{p}_0 \cos(\Omega t^* - \psi).
%\end{split}
\end{equation}
Here, $Q(T_1)$ and its complex conjugate $\overline{Q}(T_1)$ are unknowns that can be expressed in the following polar form:
\begin{equation}
\begin{split}
    Q(T_1) &= \frac{q(T_1)}{2} e^{i\nu (T_1)},\\
    \overline{Q}(T_1) &= \frac{q(T_1)}{2} e^{-i\nu (T_1)},
\end{split}
\end{equation}
where $q$ and $\nu$ are, respectively, the amplitude and phase of the perturbation $\mathcal{p}$. Since we are interested in obtaining the point at which the first period-doubling bifurcation occurs, we seek to approximate the solution when $\Omega$  is near $2 \sqrt{G_0} $. Thus, we express the proximity of the excitation frequency to twice the natural frequency by introducing
\begin{equation} \label{eq: Hills detuning}
    \Omega=2 \sqrt{G_0} + \epsilon \sigma.
\end{equation}
Upon substituting (\ref{eq: order0 soln PD}) and (\ref{eq: Hills detuning}) into Equation (\ref{eq: order1  PD}), then eliminating the secular terms, we obtain the following equation which governs the locus of the period-doubling bifurcation as function of the design parameters of the absorber:
\begin{equation} \label{eq: Hill's Soln}
%\begin{split}
    (\sqrt{G_0} \Omega - 2G_0-\delta_1 \sqrt{G_0} \xi_2) ^2 
    + ( \delta_1 \sqrt{G_0} \bar{\xi}_2 + \delta_2 \sqrt{G_0} )^2 = \frac{G^2_1}{4},
 %   \end{split}
\end{equation}
where
\begin{equation} \label{eq:C_3 -pd}
    \begin{split}
        \xi_2&=\left( \frac{\lambda_1 \sqrt{G_0}}{\mu^2+G_0} + \frac{2\lambda_2 \mu^2 G_0}{4\mu^4+G_0^2} -\frac{\lambda_3 G_0^{\frac{2}{3}}}{4\mu^4+G_0^2} \right),\\
        \bar{\xi}_2&=\left( \frac{\lambda_1 \mu}{\mu^2+G_0} + \frac{\lambda_2(2\mu^3-\mu G_0)}{4\mu^4+G_0^2}-\frac{\lambda_3(2\mu^3+\mu G_0)}{4\mu^4+G_0^2} \right).
    \end{split}
\end{equation}

\subsection{Symmetry-break bifurcations of the inter-well branch}
The other dynamical features that we have a particular interest in estimating are the points of symmetry-break bifurcation on the symmetric inter-well solution branch. To this end, we examine the stability of the approximate expansion given in Equation (\ref{eq: interwell oscillations}) by introducing an infinitesimal perturbation $\mathcal{p}(t^*)$ to the periodic solution $Y(t^*)$ as
\begin{equation} \label{eq:CF perturbation}
    \tilde{Y}(t^*) = Y(t^*) + \mathcal{p}(t^*).
\end{equation}
Upon substituting Equation (\ref{eq:CF perturbation}) into Equation (\ref{eq:ComminsWN}), we obtain the following  Mathieu’s-type differential equation:
\begin{equation} \label{eq: Hill's CF}
          \mathcal{p}''(t^*) + \delta_1 \int_{0}^{t^*} h(t^*-\tau) \mathcal{p}'(\tau) d\tau + \delta_2 \mathcal{p}'(t^*)
          + \left(K_0 + \sum_{n=1}^N K_n \cos(2n \Omega t^*) \right) \mathcal{p}(t^*) =0.
 \end{equation}
Because of the presence of the term, $K_n \cos(2n \Omega t^*)$ in Equation (\ref{eq: Hill's CF}),  $\mathcal{p}(t^*)$ must admit solutions that have even frequencies, which breaks the symmetry of the original solution, $Y(t^*)$. Since the parametric terms are periodic with period $T=\pi/\Omega$, it follows by the virtue of Floquet theory that $\mathcal{p}(t^*)$ satisfies the following equation:
 \begin{equation}
     \mathcal{p}(t^*+T)=\lambda \mathcal{p}(t^*),
 \end{equation}
 where $\lambda$ is an eigenvalue called the Floquet multiplier. This multiplier is an eigenvalue of a matrix $C$ associated with the solution $\Phi(t^*)$ of the Equation (\ref{eq: Hill's CF}), which satisfies $\Phi(t^*+T)=\Phi(t^*)C$.  The matrix $C$, also known as the monodromy matrix, can be thought of as a transformation that maps $\Phi(0)$ to $\Phi(T)$. Specifying the initial condition $\Phi(0)=I$ yields
\begin{equation}
    C=\Phi(T).
\end{equation} 
 In order to assess the stability of $\mathcal{p}(t^*)$, we solve for $\Phi(t^*)$ starting at $t^*=0$ and ending at  $t^*=T$. We then examine the eigenvalues of the matrix $C=\Phi(T)$. The set of differential equations governing the fundamental matrix solution can be constructed using the following equation: 
\begin{equation} \label{eq:fundamental matrix solution}
\frac{d}{dt^*} \Phi(t^*)=
   \left[\includegraphics[width=0.25\textwidth ,valign=c]{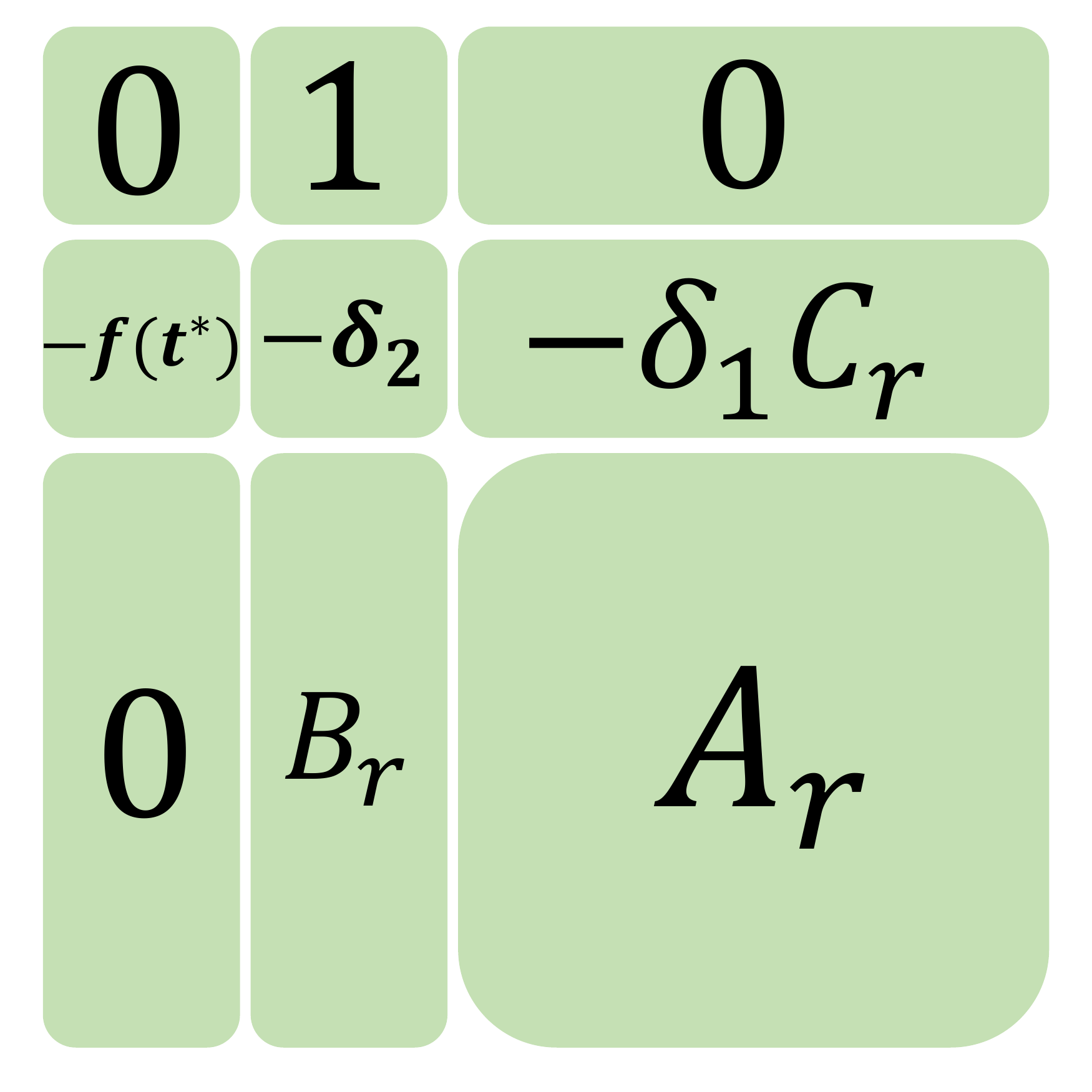}\right]
   \Phi(t^*),
\end{equation}
where $A_r$, $B_r$ and $C_r$ are the realized state-space accounting for radiation damping, and  $f(t^*)$ is the parametric excitation term given as:
\begin{equation}
    f(t^*)=K_0 + \sum_{n=1}^N  K_n \cos(2n \Omega t^*).
\end{equation}
The reader can refer to Appendix \ref{Apndx: Parametric constants} for the constants $K_i$. Equation (\ref{eq:fundamental matrix solution}) results in a set of 25 linear differential equations subjected to the initial conditions $\Phi(0)=I$. Here, $I$ is an identity matrix of dimension 5. 

Integrating Equation (\ref{eq:fundamental matrix solution}) numerically in [0, $T$] using the initial conditions of $\Phi(0)=I$, finding the Floquet multipliers, $\lambda$, of the resulting numerical matrix, then inspecting their location with respect the unit circle, we can find the loci of the symmetry breaking bifurcations in the wave amplitude versus frequency parameter space. In particular, a symmetry breaking bifurcation occurs when one of the Floquet multipliers exits the unit circle through $\lambda=1$.

\subsection{Bifurcation diagram} \label{sec:RES}
\begin{figure}  [h!]
\centering
  \begin{tabular}{@{}c@{}}
  $(a)$\\
    \includegraphics[width=.75\textwidth]{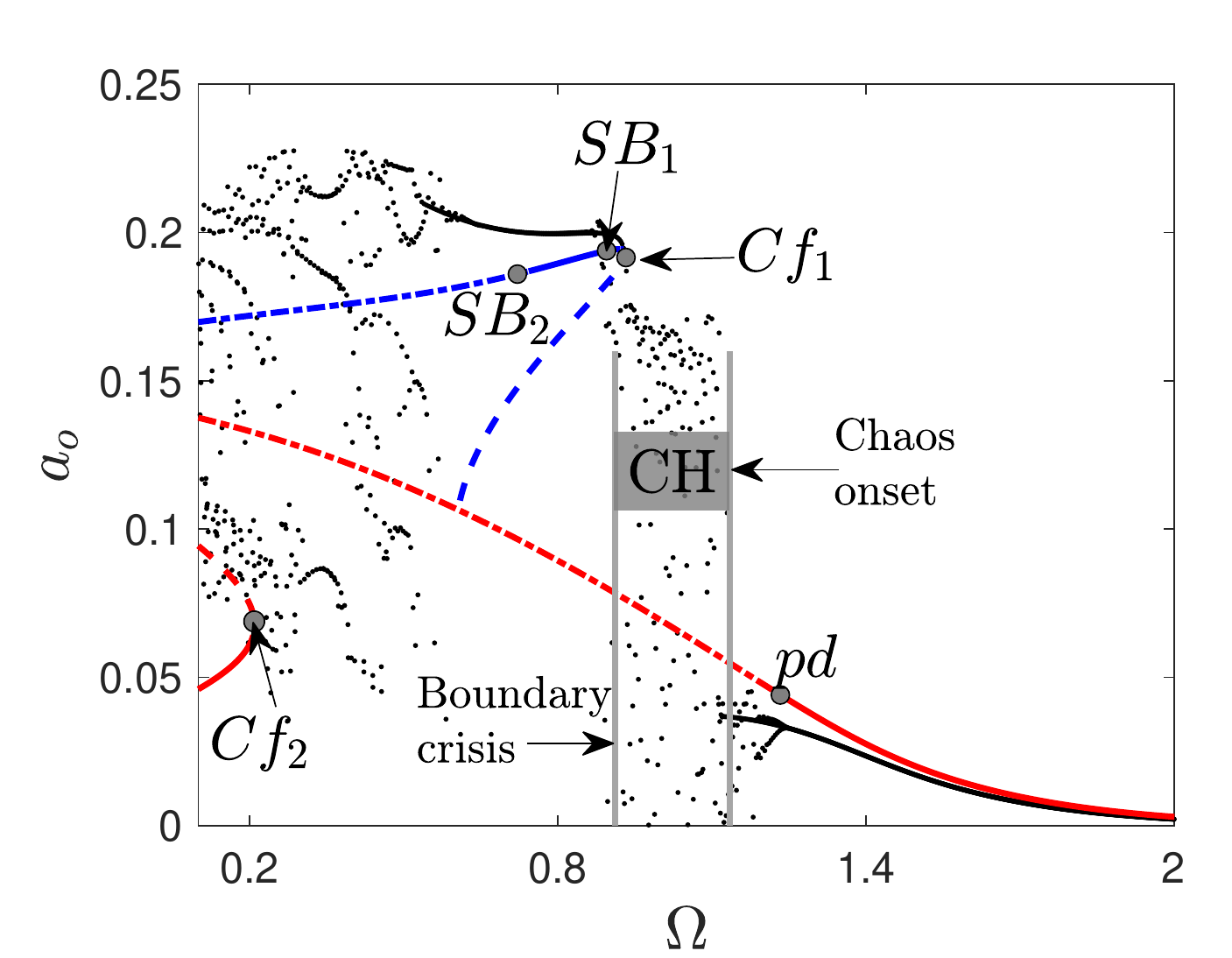} 
      \end{tabular}
  \begin{tabular}{@{}cc@{}}
  $(b)$ &
  $(c)$ \\
    \includegraphics[width=.5\textwidth]{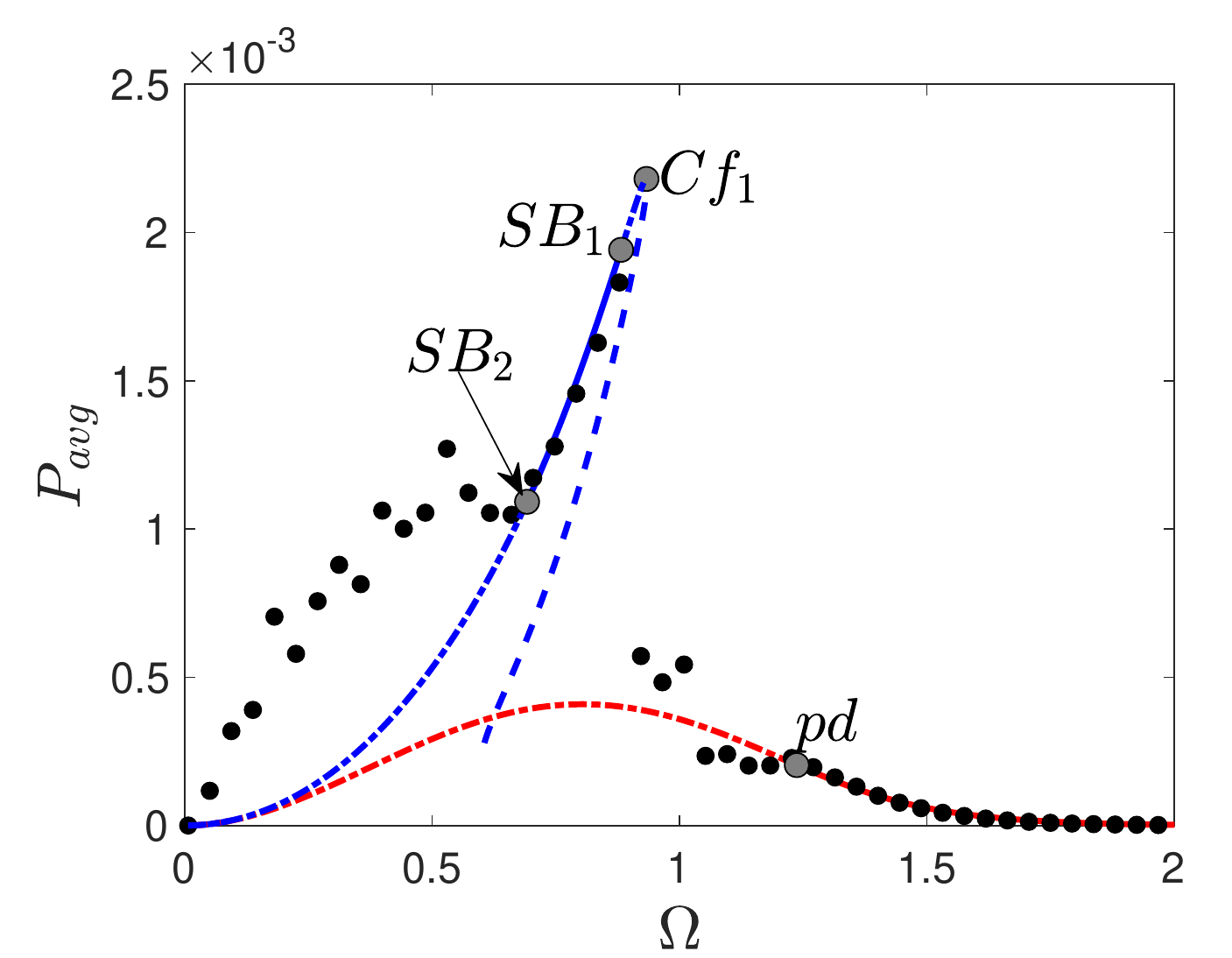} &
    \includegraphics[width=.5\textwidth]{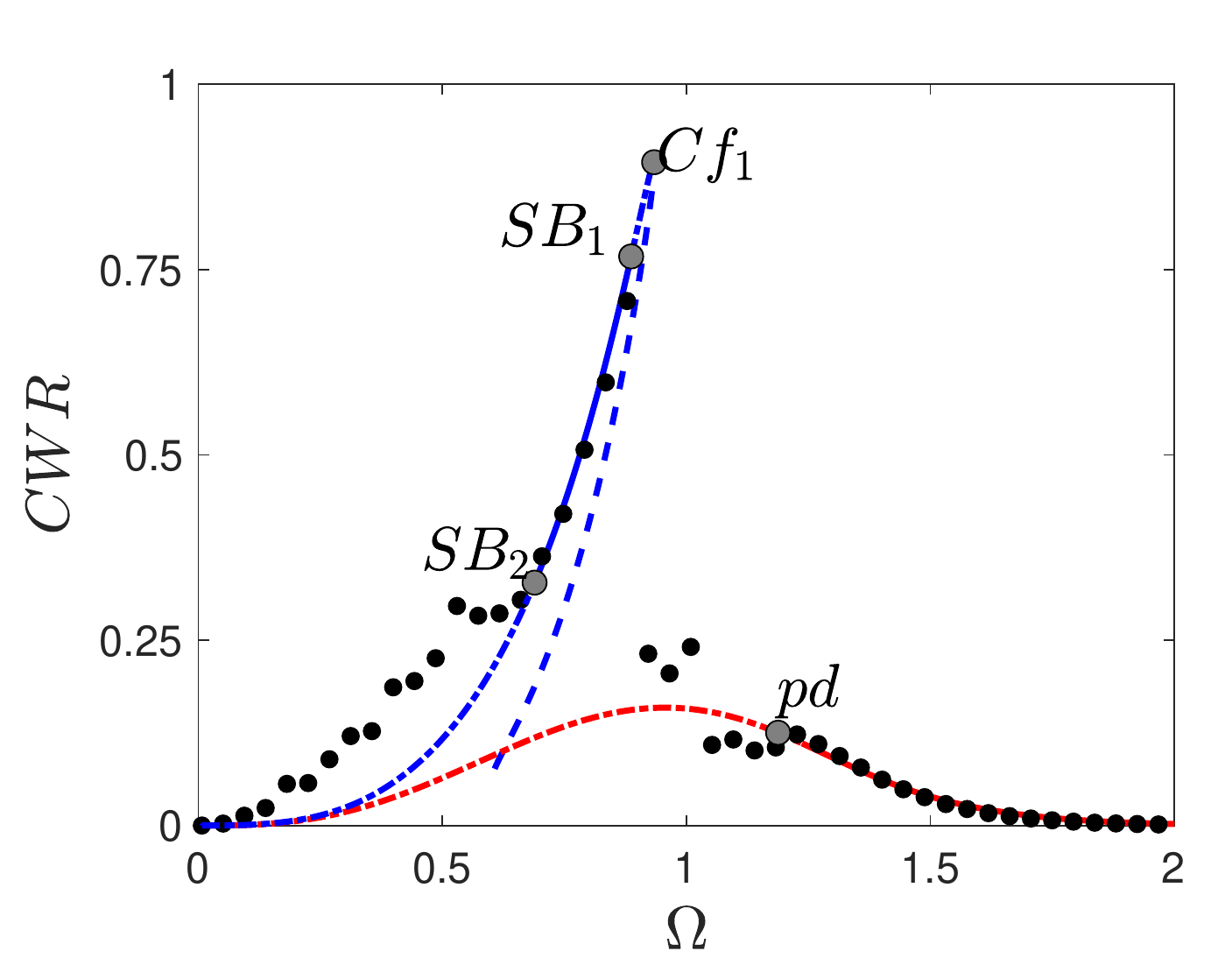}
      \end{tabular}
   \caption{Bi-stable PWA response under regular wave excitation at nondimensional wave amplitude $A_{wave}/R=0.1$. (a): Analytical and stroboscopic bifurcation map. (b): Nondimensional averaged power. (c): CWR. (Solid lines represent stable solutions. Dashed lines represent unstable solutions based on the Jacobian matrix criterion. Red dashed-dotted lines represent unstable solutions based on period-doubling instability analysis. Blue dashed-dotted lines represent unstable solutions based on Floquet analysis. Circles represent the numerical results). }
  \label{fig: Regular waves POT2 Results}
\end{figure}

The bifurcation diagram based on the aforedescribed stability analysis is shown in Figure \ref{fig: Regular waves POT2 Results} (a). It is evident that the bifurcations $SB_1$, $SB_2$, and $pd$ are in excellent agreement with the points where the stroboscopic numerical bifurcation maps undergo qualitative changes. The average power and the CWR curves shown in Figures \ref{fig: Regular waves POT2 Results} (b-c) also illustrate good agreement with the analytical solution in the range when the response exhibits unique period one periodic motion. It is evident that maximum average power and CWR are realized within the effective bandwidth of the absorber and that outside that bandwidth the average power drops whether the response is of the chaotic or intra-well type.

\section{The Effective Bandwidth} \label{sec: maps}
In this section, we use the stability analysis furnished in the previous sections to define an effective bandwidth for the bi-stable wave energy absorber by marking the boundaries at which the PWA switches to different types of motion in the wave amplitude versus frequency $(A_{wave}/R, \Omega)$ parameters space. We also draw a clearer picture of how the shape of the bi-stable potential influences the effective bandwidth.

\begin{figure} [h!]
\centering
\includegraphics[width=.75\textwidth]{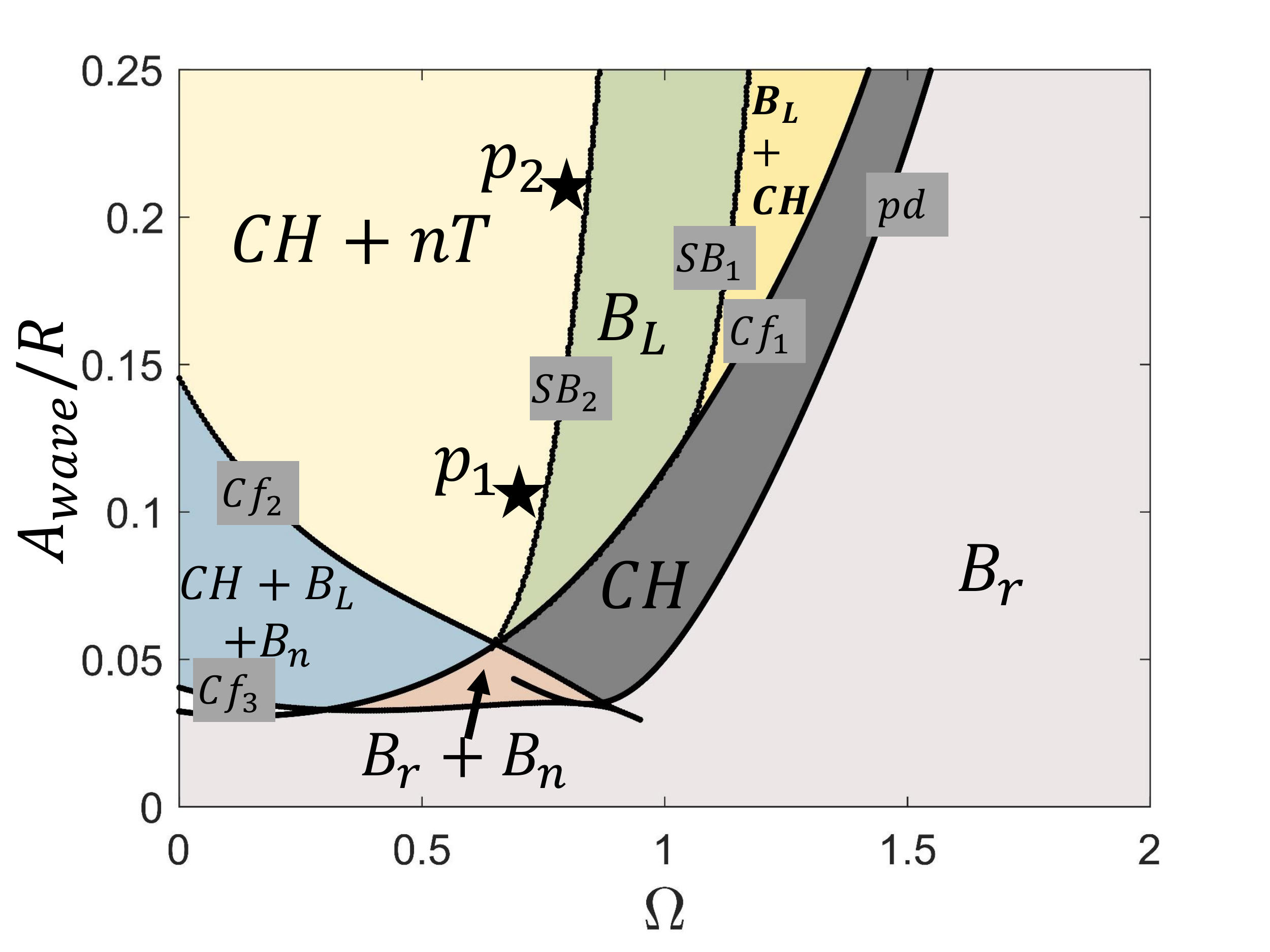}
\vspace{0.0cm}
\caption{A map demarcating regions of quantitatively different PWA responses.}
\label{fig:map1}
\end{figure}

\begin{figure} [h!]
\centering
\includegraphics[width=.5\textwidth]{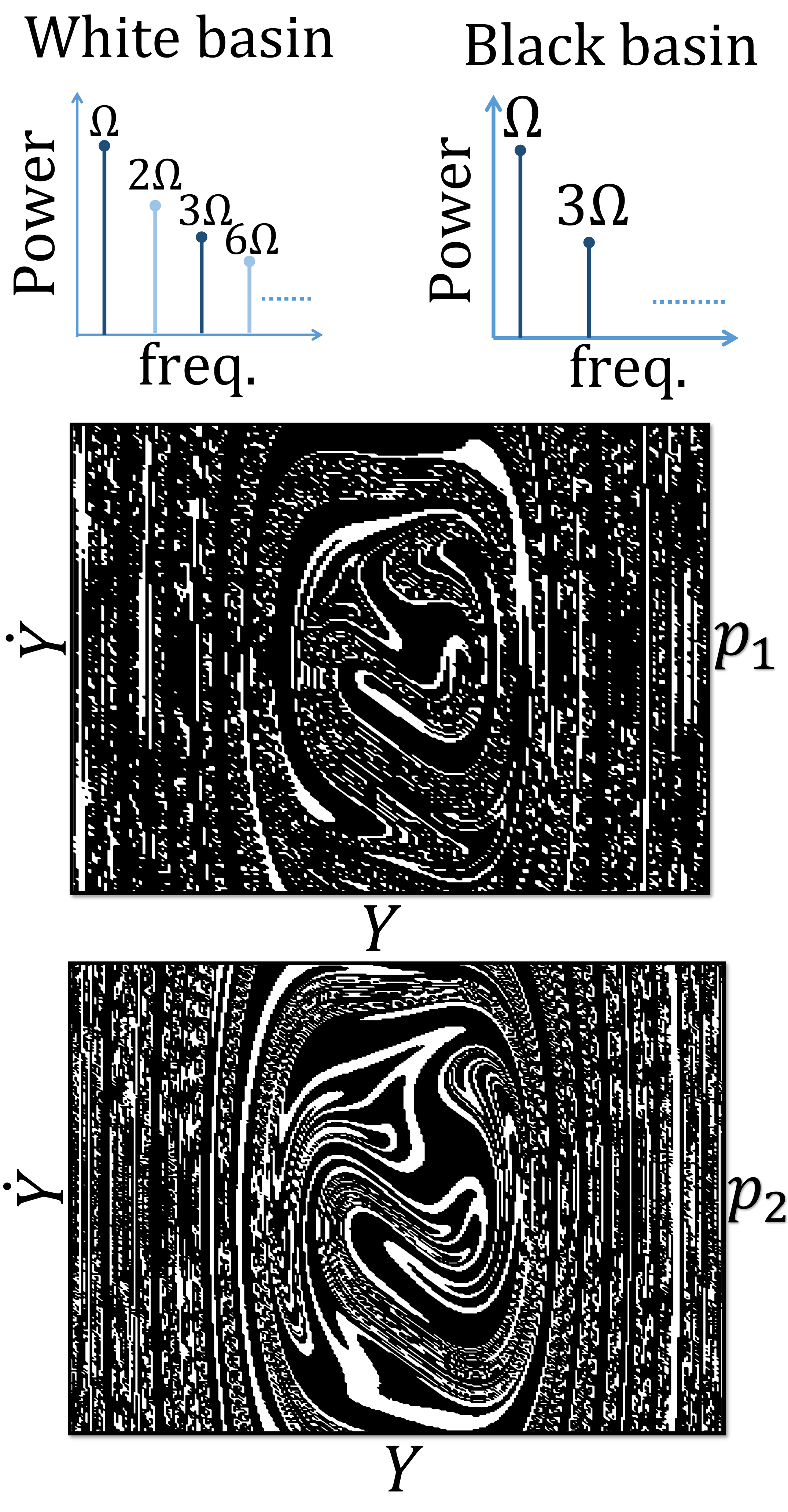}
\vspace{0.0cm}
\caption{Basins of attraction showing different coexisting orbits. }
\label{fig:basins}
\end{figure}

Using the loci of the different bifurcations, we create the design map shown in Figure \ref{fig:map1} which characterizes the type of motion realized for every combination of wave amplitude and frequency. The largest region on the figure is the one denoted by $B_r$, which represents $(A_{wave}/R, \Omega)$ combinations that incite small-amplitude inter-well motions. We can see that this type of response occurs for any frequency when the wave amplitude is small and for any wave amplitude when the wave frequency is large. In terms of size, the second region on the figure is that which results in coexisting chaotic, $CH$, and $n$-period, $nT$, periodic solutions. This region occurs slightly to the left of $\Omega=1$ and increases in bandwidth as  $A_{wave}/R$ is increased. The basin of attraction of the different coexisting solutions in this region are shown in Figure \ref{fig:basins} for the two points indicated on the map. We can see that there are three different basins: the black one represents symmetric period one motions, the white one represents asymmetric period one motions, and the one where the black and white colors blend together represents chaotic motions.

The most important region on the map is that denoted by $B_L$, which corresponds to the combination of wave parameters leading to a unique large-amplitude inter-well motion. We see that this region exists around $\Omega=1$ for wave amplitudes larger than $A_{wave}/R=0.05$. The size of the effective bandwidth of the absorber increases as the wave amplitude is increased up to a value of $A_{wave}/R=0.125$ beyond which the size of the bandwidth remains almost constant. The $B_L$ region is bordered from the right by the region, $B_L+CH$, where large-amplitude periodic motions coexist with chaotic motions, and from below by the region, $CH$, where the chaotic attractor is unique. The region denoted by $CH+B_L+B_n$ corresponds to wave parameters that results in a chaotic attractor coexisting with two types of periodic orbits: $B_L$ and $B_n$ each with competing basins of attraction.

\subsection{Influence of the shape of the potential energy function}
\begin{figure} [h!]
\centering
\includegraphics[width=0.75\textwidth]{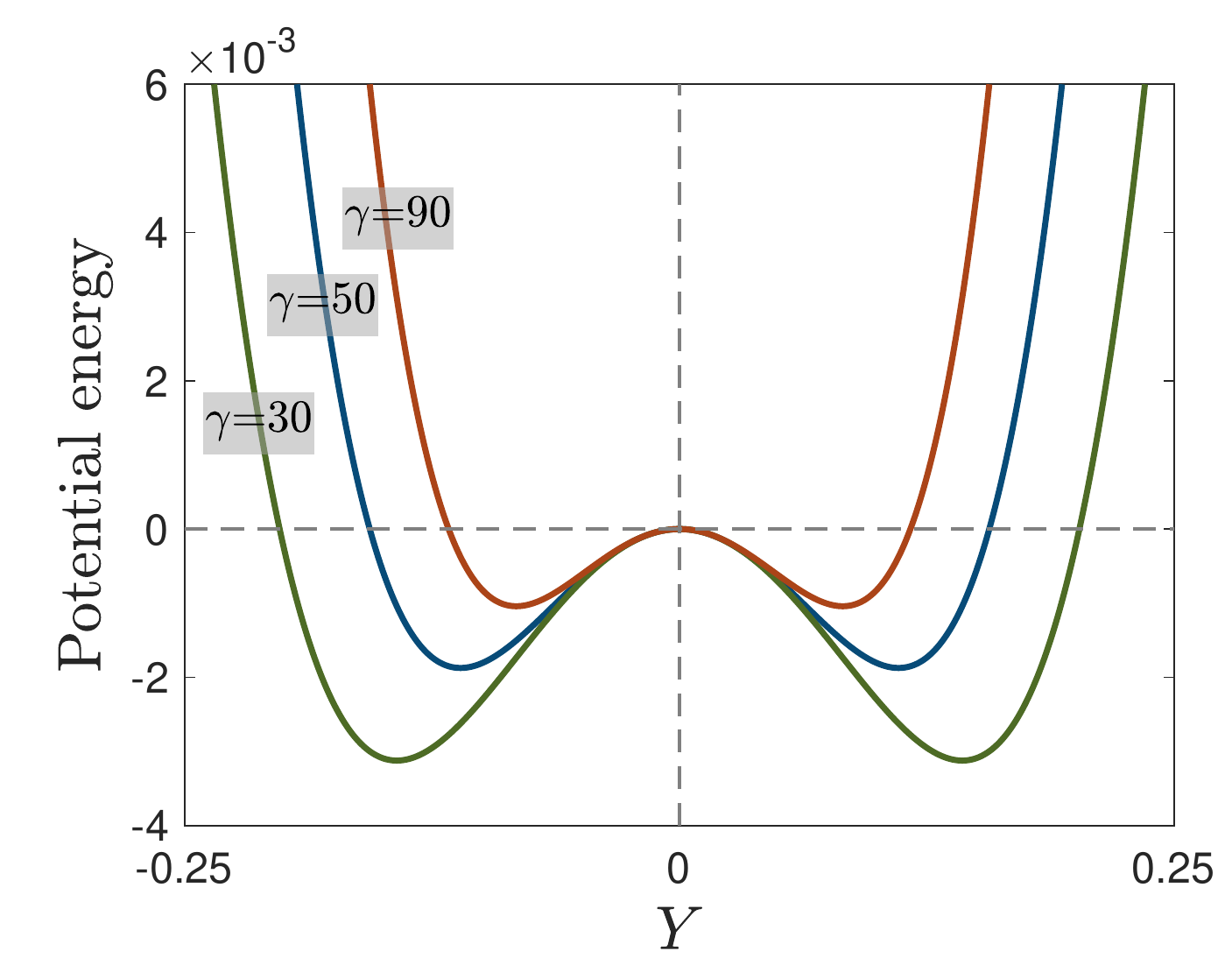}
\vspace{0.0cm}
\caption{Potential energy function of the bi-stable wave energy absorber at different values of $\gamma$ and $\omega_n=0.78$. (1): $\gamma=30$, (2): $\gamma=50$, (3): $\gamma=90$.}
\label{fig:POT123}
\end{figure}

In this section, we use the understanding developed in the previous section to examine the influence of the potential function on the effective bandwidth of the bi-stable PWA. In particular, we want to understand how the depth of the potential wells and their separation influences the size of the effective bandwidth and the other regions in the bifurcation map. To this end, we consider the three different potential energy functions shown in Figure \ref{fig:POT123}, which were obtained using $\gamma=$30, 50, and 90. The potential energy function associated with $\gamma=$30 is deeper with a larger separation between the stable equilibria, while the potential energy function associated with $\gamma=$90 has shallower potential wells and smaller separation. 

\begin{figure}
     \centering
     \begin{subfigure}[b]{0.5\textwidth}
         \centering
         \includegraphics[width=\textwidth]{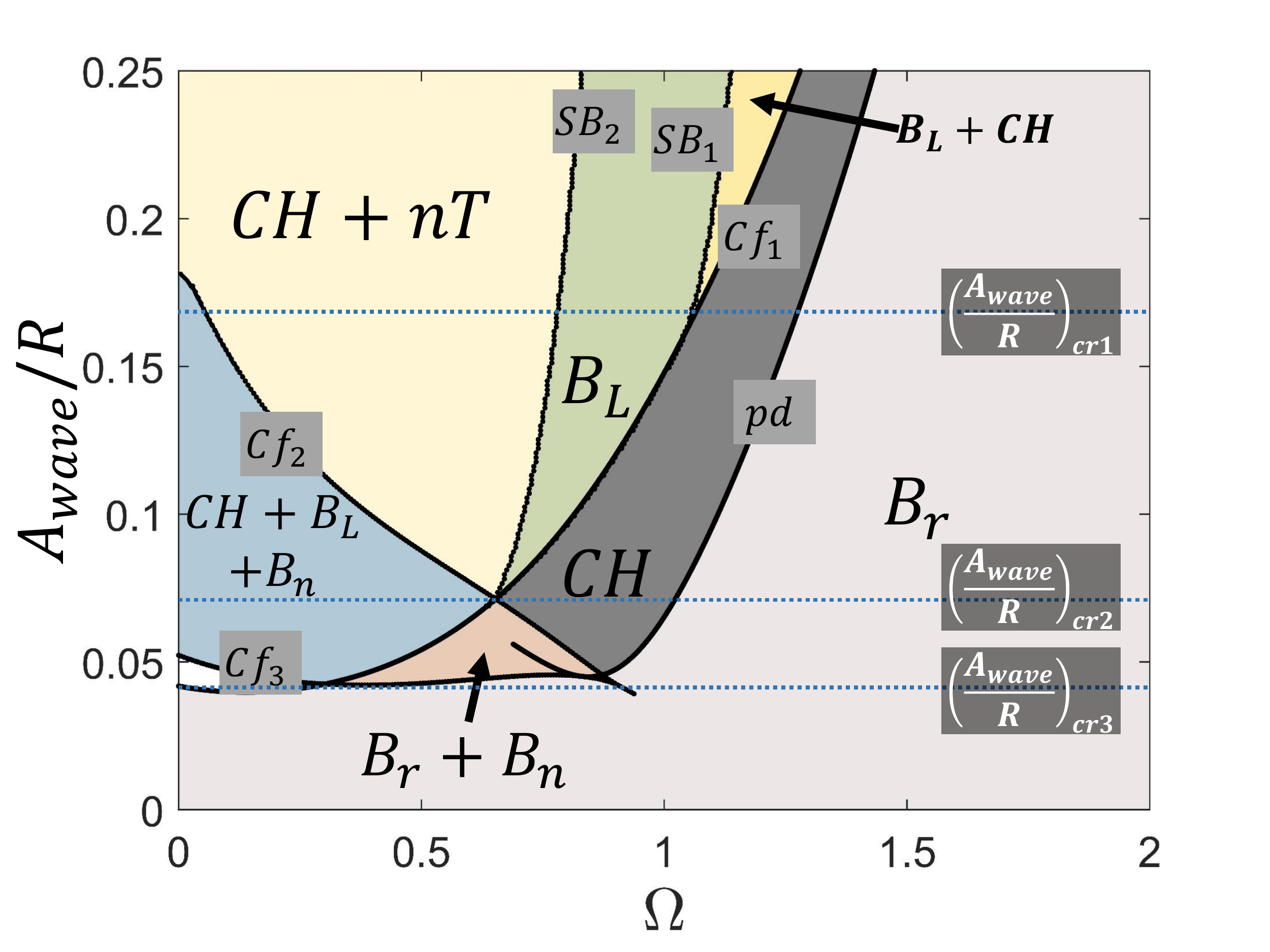}
         \caption{$\gamma=30$}
     \end{subfigure}
     \hfill
     \begin{subfigure}[b]{0.5\textwidth}
         \centering
         \includegraphics[width=\textwidth]{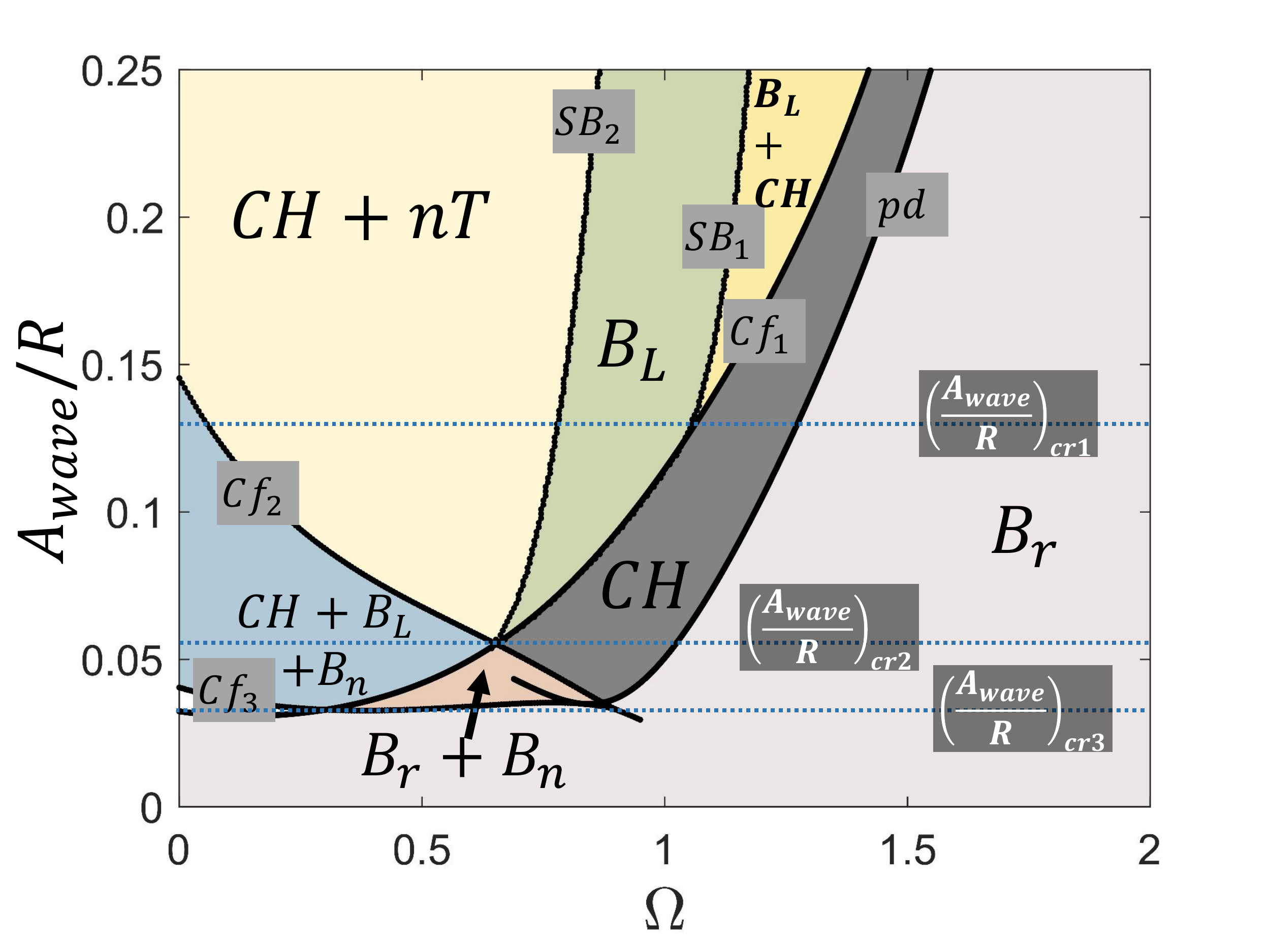}
         \caption{$\gamma=50$}
     \end{subfigure}
     \hfill
     \begin{subfigure}[b]{0.5\textwidth}
         \centering
         \includegraphics[width=\textwidth]{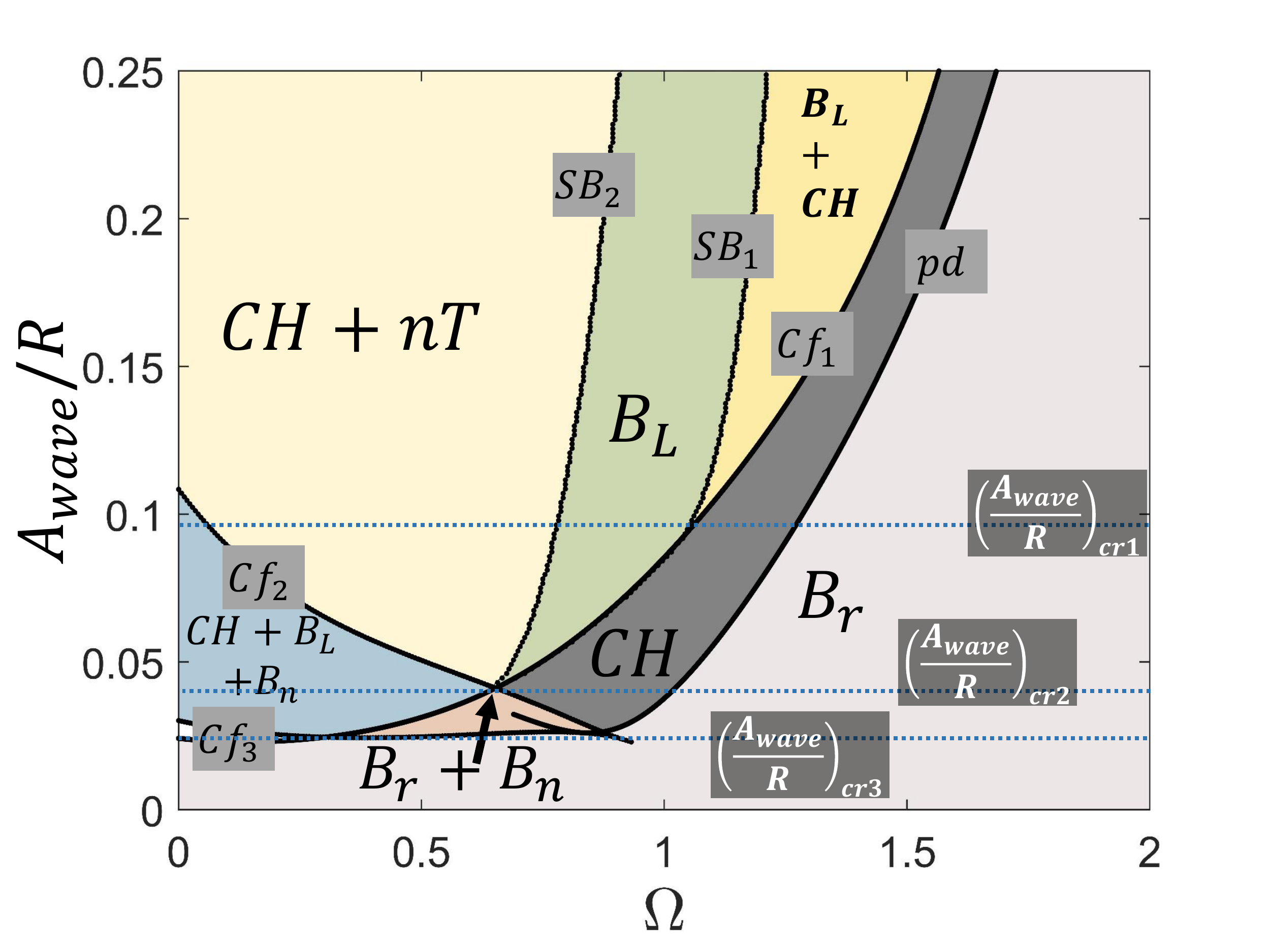}
         \caption{$\gamma=90$}
     \end{subfigure}
        \caption{Bifurcation maps in the wave amplitude versus frequency parameter space for the different potential functions shown in Figure \ref{fig:POT123}.}
        \label{fig:three graphs}
\end{figure}

The bifurcation maps showing the effective bandwidth associated with each of the potential energy functions considered are depicted in Figure \ref{fig:three graphs}. In order to make a quantitative comparison, three critical wave amplitudes were marked on the figures, and labeled as $(A_{wave}/R)_{cr1}$, $(A_{wave}/R)_{cr2}$, and $(A_{wave}/R)_{cr3}$. The first critical amplitude occurs at the intersection between the $Cf_1$ and $SB_1$ curves and can be used to define the wave amplitude at which the effective bandwidth of the absorber approaches its maximum size. We notice that this critical level of wave amplitude increases as the potential wells become deeper; that is, larger excitation levels become necessary to attain the unique periodic large-amplitude motions when increasing the depth and the separation distance between the potential energy wells. 

\begin{figure}  [h!]
\centering
  \begin{tabular}{@{}c@{}}
    \includegraphics[width=.5\textwidth,valign=c]{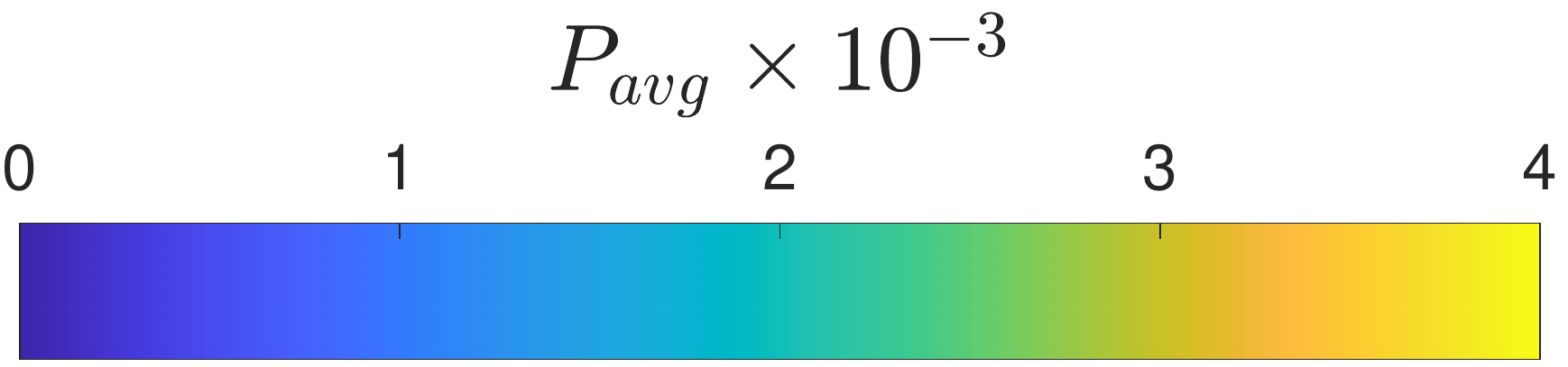}
      \end{tabular}
  \begin{tabular}{@{}cc@{}}
  $(a)$:  $\gamma=30$ &
  $(b)$:  $\gamma=50$\\
    \includegraphics[width=.475\textwidth,valign=c]{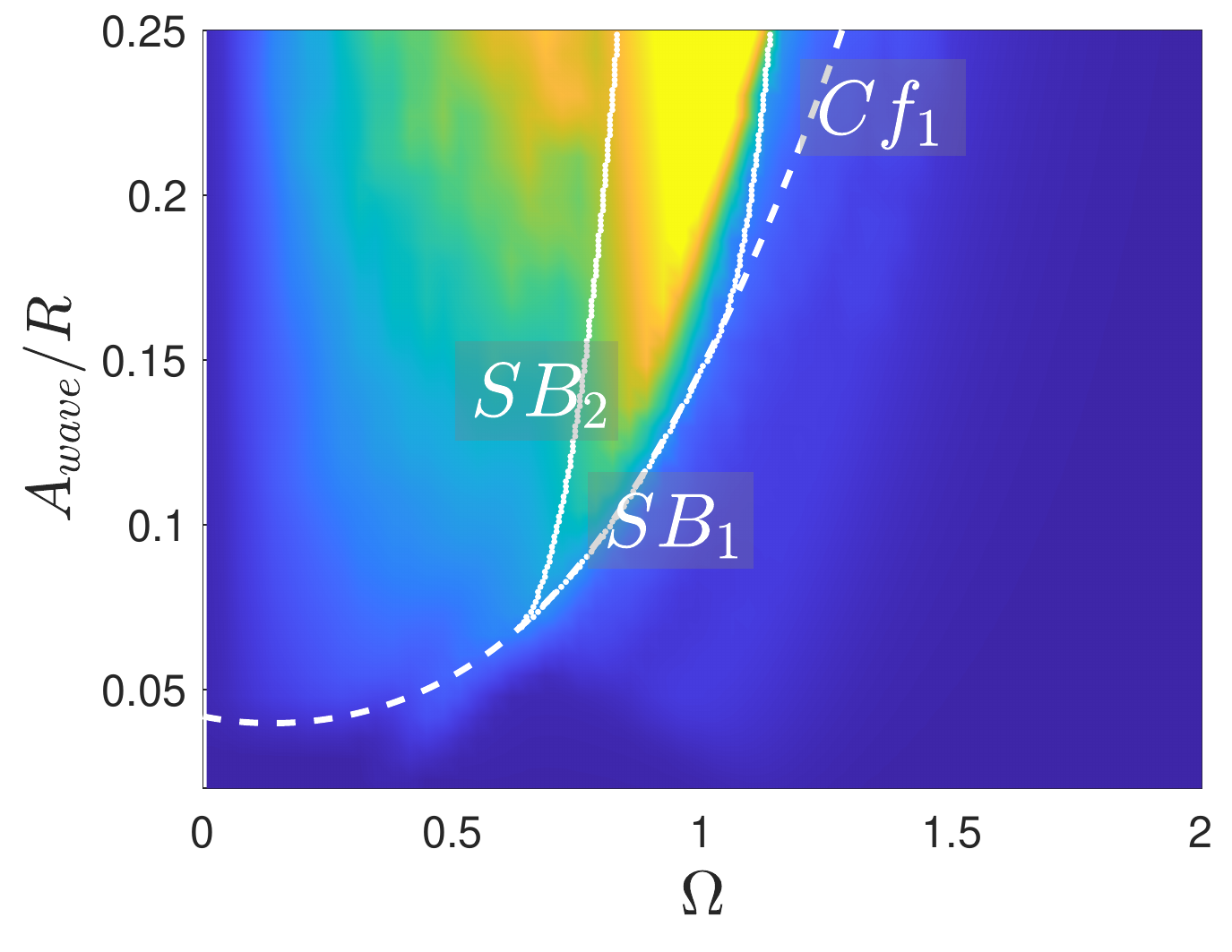} & 
    \includegraphics[width=.5\textwidth,valign=c]{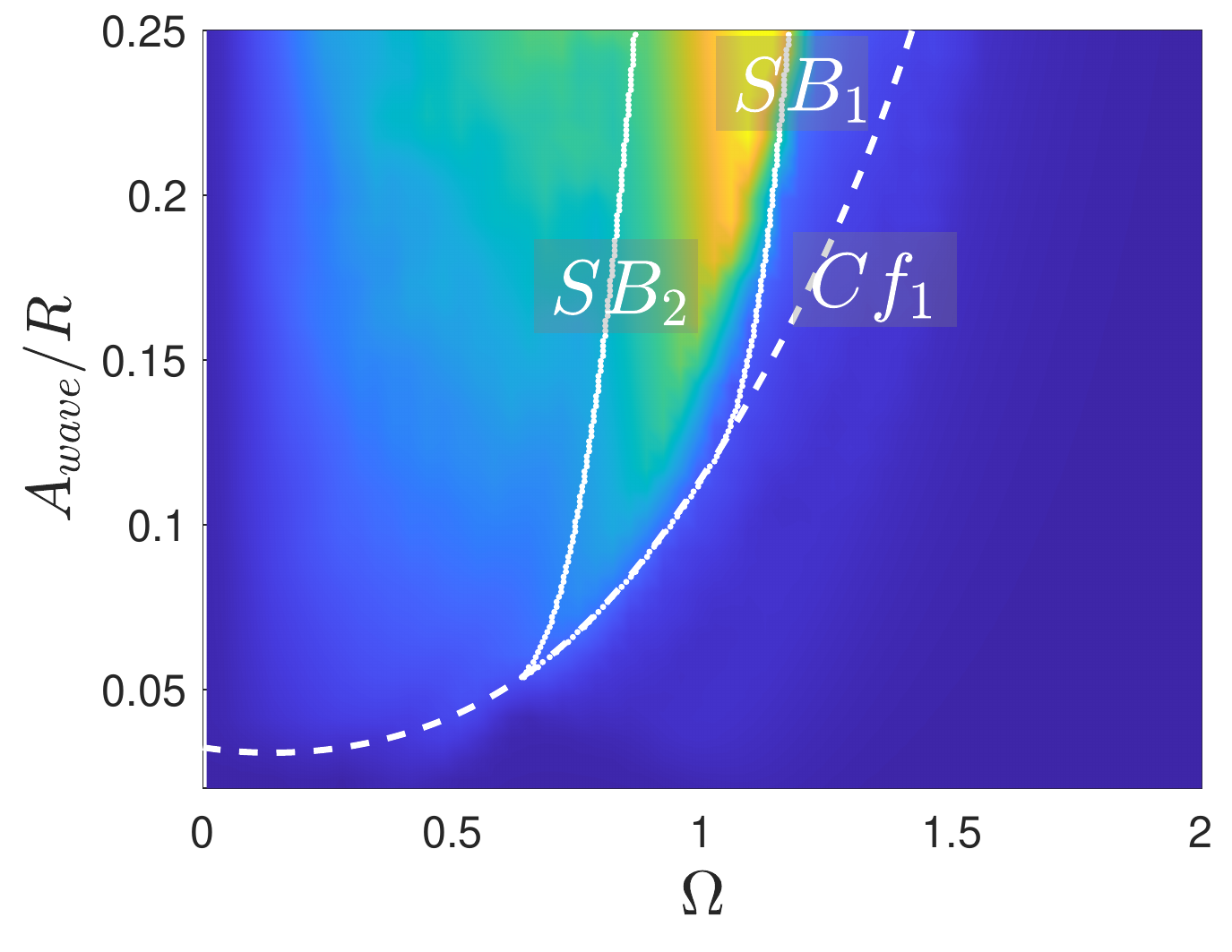} \\
      \end{tabular}
  \begin{tabular}{@{}c@{}}
  $(c)$:  $\gamma=90$\\
    \includegraphics[width=.5\textwidth,valign=c]{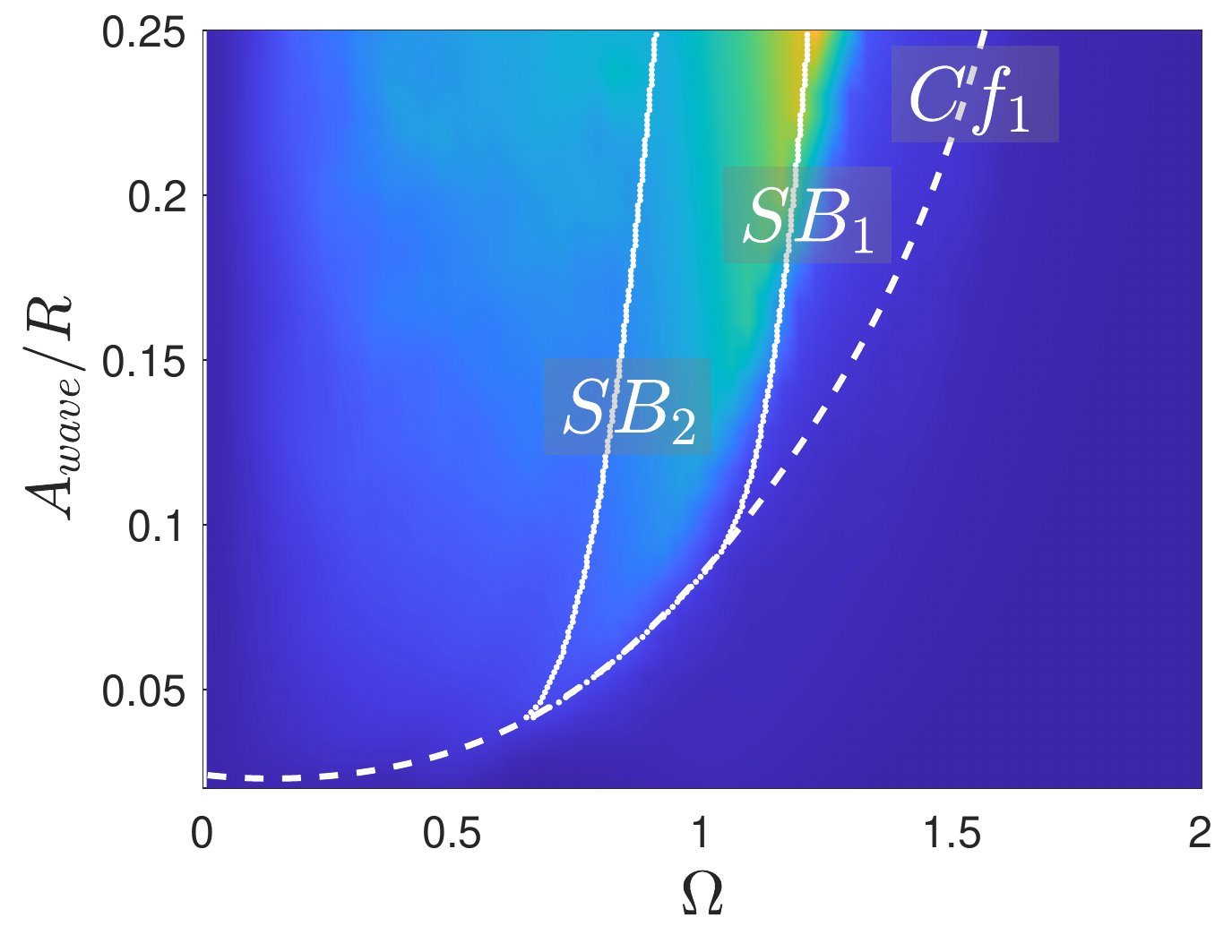} 
      \end{tabular}
  \caption{Comparison of the average generated power in the wave amplitude - frequency parameter space for the different potential functions shown in Figure \ref{fig:POT123}. The numerical simulations were performed on Equation (\ref{eq:ComminsWN}) at initial conditions $(Y_0, \dot{Y}_0)$ of $(0, 0)$.} 
  \label{fig: PWR MAPS 123}
\end{figure}

It is interesting to note that, for any wave amplitude above $(A_{wave}/R)_{cr1}$, the size of the effective bandwidth remains almost constant regardless of the shape of the potential energy function. However, while the effective bandwidth remains unchanged, the power levels within the effective bandwidth change considerably with the shape of the potential energy function as shown in Figure \ref{fig: PWR MAPS 123}. The bi-stable PWA with the deeper potential energy wells produces higher average power levels within the effective bandwidth (the region bounded between $SB_1$ and $SB_2$).

This critical level of excitation marked by $(A_{wave}/R)_{cr2}$ occurs at the intersection between the $Cf_1$ and $pd$ lines. It represents the minimum value of wave amplitude necessary to generate unique inter-well motions.  It is evident that this amplitude decreases as the potential wells become shallower. Finally, the amplitude level $(A_{wave}/R)_{cr3}$, represents the wave amplitude below which no bifurcations occur as the wave frequency is varied. At such a low level of wave amplitude, the response of the PWA resembles the bell-shaped response of the traditional linear PWA.

%\pagebreak
%\pagebreak
%\pagebreak
%\pagebreak
%\pagebreak

\section{Conclusion} \label{sec:conclusion}
This paper presented a theoretical analytical analysis of the response of bi-stable PWAs to harmonic wave excitations. To this end, approximate asymptotic solutions of the governing equations of motion were derived by implementing the method of multiple scales. A stability analysis of the attained solutions revealed the presence of key bifurcations that can be used to define an effective bandwidth of the generator. This effective bandwidth is characterized by the presence of a unique large-orbit inter-well motion for a set of wave amplitudes and frequencies. This effective bandwidth occurs slightly below the resonant frequency of the absorber and exists only above a certain threshold in the wave amplitude. This threshold increases as the depth of the potential well is increased. The size of the effective bandwidth increases as the wave amplitude is increased up to a certain threshold above which the effective bandwidth remains almost constant even when the wave amplitude is substantially increased. The size of the effective bandwidth is observed to be insensitive to variations in the depth of the potential well of the absorber. However, the power levels within the effective bandwidth change considerably with the shape of the potential energy function. In particular, a bi-stable PWA with deeper potential energy wells produces higher average power levels within its effective bandwidth. It is our belief that this comprehensive analytical treatment is key to designing effective bi-stable PWAs for known wave conditions and provide backbone results for future studies addressing more realistic regular non-harmonic wave excitations.

\section*{Funding}
This research was funded by Abu Dhabi Education and Knowledge Council (ADEK) under grant number AARE2019-161: Exploiting Bi-stability to Develop a Novel Broadband Point Wave Energy Absorber.
\section*{Conflict of interest}
The authors declare that they have no conflict of interest.

\section*{Data availability}
The data that support the findings of this study will be made available upon reasonable request.

\break
\appendix
\begin{appendices}
\numberwithin{equation}{section}
\makeatletter 
\newcommand{\section@cntformat}{Appendix \thesection:\ }
\makeatother
%\pagebreak
\section{Eigensystem realization algorithm} \label{apndx: ERA}
Consider the following single-input single-output discrete-time dynamical system:
\begin{equation}
\begin{split}
    \textbf{x}_{k+1} &=\textbf{A}\textbf{x}_k + \textbf{B}u_k,\\
    h_{k} &= \textbf{C}\textbf{x}_k + \textbf{D} u_k,
\end{split}
\end{equation}
and a discrete-time scalar input $u$:
\begin{equation}
u_{k}^{\delta}  \equiv u^{\delta} (k \Delta t)=
  \begin{cases}
  			1, & \text{if $k=0$}\\
            0, & \text{if $1, \dots \infty$}
  \end{cases}
\end{equation}

According to linear system theory, the discrete-time impulse response data $h_{k}^{\delta}$ could be expressed as:

\begin{equation} \label{impusledata2}
    h_{k}^{\delta} \equiv h^{\delta}(k \Delta t)=\textbf{C} \textbf{A}^{k} \textbf{B}, \hspace{0.5 cm} \left(k=0,1,\dots, \infty\right)
\end{equation}
In our analysis, the discrete-time impulse response data are obtained through substituting the radiation damping coefficients $B(\omega_i)$ into Equation (\ref{eq:ogilvie}). The next step is to proceed by constructing the generalized \textit{Hankel} matrix  $\textbf{H}_{r \times s} (g)$ for $g=0,1$  which consists of $r$ rows and $s$ columns, generated by stacking time-shifted impulse response data in the following order:
\begin{equation} \label{Hankel1}
\textbf{H}_{r \times s} (g) =
  \begin{pmatrix}
    h_{g}^{\delta} & h_{g+1}^{\delta} & \dots & h_{g+s-1}^{\delta} \\
    h_{g+1}^{\delta} & h_{g+2}^{\delta} & \dots & h_{g+s}^{\delta} \\
    \vdots & \vdots & \ddots & \vdots \\
    h_{g+r-1}^{\delta} & h_{g+r}^{\delta} & \dots & h_{g+r+s-2}^{\delta}
  \end{pmatrix}
\end{equation}
Using Equation (\ref{impusledata2}), we can express the generalized \textit{Hankel} matrix in Equation (\ref{Hankel1}) in terms of the realized state-space matrices $\textbf{A}_r$, $\textbf{B}_r$ and $\textbf{C}_r$ as:
\begin{equation} \label{Hankel2}
\textbf{H}_{r \times s} (g) =
  \begin{pmatrix}
    \textbf{C}_r \textbf{A}_r^{g} \textbf{B}_r & \textbf{C}_r \textbf{A}_r^{g+1} \textbf{B}_r & \dots & \textbf{C}_r \textbf{A}_r^{g+s-1} \textbf{B}_r \\
    \textbf{C}_r \textbf{A}_r^{g+1} \textbf{B}_r & \textbf{C}_r \textbf{A}_r^{g+2} \textbf{B}_r & \dots & \textbf{C}_r \textbf{A}_r^{g+s} \textbf{B}_r \\
    \vdots & \vdots & \ddots & \vdots \\
    \textbf{C}_r \textbf{A}_r^{g+r-1} \textbf{B}_r & \textbf{C}_r \textbf{A}_r^{g+r} \textbf{B}_r & \dots & \textbf{C}_r \textbf{A}_r^{g+r+s-2} \textbf{B}_r
  \end{pmatrix}
\end{equation}
which could be reduced as:
\begin{equation} \label{ReducedHankel}
   \textbf{H}_{r \times s} (g) = \mathcal{O} \textbf{A}_r^{g} \mathcal{C}
\end{equation}
where
    \begin{align*}
        \mathcal{O} &=\left(\textbf{C}_r \hspace{0.5 cm} \textbf{C}_r \textbf{A}_r \hspace{0.5 cm} \dots \hspace{0.5 cm} \textbf{C}_r \textbf{A}_r^{r-1}\right)^{T}\\ 
        \mathcal{C} &=\left(\textbf{B}_r \hspace{0.5 cm} \textbf{A}_r \textbf{B}_r \hspace{0.5 cm} \dots \hspace{0.5 cm} \textbf{A}_r^{r-1} \textbf{B}_r\right)^{T}
    \end{align*}
are respectively the generalized observability and controllability matrices, with observability and controllability indices of $r$ and $s$. Upon taking the singular value decomposition (SVD) for the first \textit{Hankel} matrix $\textbf{H}_{r \times s} (0)$ we get the following definition:
\begin{equation}
\begin{aligned}
\textbf{H}_{r \times s} (0) &= \textbf{U} \Sigma \textbf{V}^{T}\\
&=   \begin{pmatrix}
    \Tilde{\textbf{U}} & \textbf{U}_{t}\\
  \end{pmatrix}
  \begin{pmatrix}
    \Tilde{\Sigma} & 0\\
    0 & \Sigma_{t} \\
  \end{pmatrix}
    \begin{pmatrix}
     \Tilde{\textbf{V}}^T\\
     \Tilde{\textbf{V}}_{t}^T  \\
  \end{pmatrix} \\
  &\approx \Tilde{\textbf{U}} \Tilde{\Sigma} \Tilde{\textbf{V}}^T
\end{aligned}
\end{equation}
where,
    \begin{align*}
        \Tilde{\textbf{U}}^T \Tilde{\textbf{U}} &=\textbf{I}\\
        \Tilde{\textbf{V}}^T \Tilde{\textbf{V}}& =\textbf{I}\\
        \Tilde{\Sigma} &=
  \begin{pmatrix}
    \sigma_1 &          &        &   \\
             & \sigma_2 &        &   \\
             &          & \ddots &   \\
             &          &        & \sigma_N\\
  \end{pmatrix}
    \end{align*}
The diagonal matrix $\Tilde{\Sigma}$ which is constructed from the first $N \times N$ block of $\Sigma$ contains the dominant singular values $\sigma_i$ in the following order $(\sigma_1 \geq \sigma_2 \geq \dots \geq \sigma_N \geq 0)$. While $\Sigma_t$ contains the small truncated singular values. This truncation step is very vital for reducing the order of the realized state-space model; such that the realized dynamics matrix $A_r$ has size $N$. Also, vectors in $\Tilde{\textbf{U}}$ and $\Tilde{\textbf{V}}^T$ contain the dominant modes associated with the singular values retained in $\Tilde{\Sigma}$. As a result, the product $\Tilde{\textbf{U}} \Tilde{\Sigma} \Tilde{\textbf{V}}^T$  is considered to be a faithful representation of the original \textit{Hankel} matrix $\textbf{H}$ for the smallest size of $\Tilde{\Sigma}$.\\
It follows from Equation (\ref{ReducedHankel}) that
\begin{equation}
\begin{aligned}
        \textbf{H}_{r \times s} (0) &= \Tilde{\textbf{U}} \Tilde{\Sigma} \Tilde{\textbf{V}}^T\\
        &= \left(\Tilde{\textbf{U}} \Tilde{\Sigma}^{\frac{1}{2}} \right) \left( \Tilde{\Sigma}^{\frac{1}{2}} \Tilde{\textbf{V}}^T \right) = \mathcal{O} \mathcal{C}
\end{aligned}
\end{equation}
Using the above balanced decomposition of $\textbf{H}_{r \times s} (0)$ we can write:
\begin{equation*}
    \mathcal{O} = \Tilde{\textbf{U}} \Tilde{\Sigma}^{\frac{1}{2}} \hspace{1 cm} and \hspace{1 cm} \mathcal{C} = \Tilde{\Sigma}^{\frac{1}{2}} \Tilde{\textbf{V}}^T
\end{equation*}
Also from Equation (\ref{ReducedHankel}) we can express the second \textit{Hankel}  matrix $\textbf{H}_{r \times s} (1)$ as:
\begin{equation}
\begin{aligned}
    \textbf{H}_{r \times s} (1) &= \mathcal{O} \textbf{A}_r \mathcal{C}  \\
      &=  \left(\Tilde{\textbf{U}} \Tilde{\Sigma}^{\frac{1}{2}}\right) \textbf{A}_r \left(\Tilde{\Sigma}^{\frac{1}{2}} \Tilde{\textbf{V}}^T\right)
\end{aligned}
\end{equation}
Using the properties of $\textbf{U}$ and $\textbf{V}$ we can write:
\begin{equation}
    \Tilde{\Sigma}^{\frac{1}{2}} \textbf{A}_r \Tilde{\Sigma}^{\frac{1}{2}} = \Tilde{\textbf{U}}^T   \textbf{H}_{r \times s} (1) \Tilde{\textbf{V}}
\end{equation}
It follows that the matrix $\textbf{A}_r$ could be obtained through:
\begin{equation}
    \textbf{A}_r=\Tilde{\Sigma}^{-\frac{1}{2}} \Tilde{\textbf{U}}^T  \textbf{H}_{r \times s} (1) \Tilde{\textbf{V}} \Tilde{\Sigma}^{-\frac{1}{2}}
\end{equation}
Let:
\begin{equation*}
\textbf{E}^T_1 =
      \begin{pmatrix}
    1 & 0 & \dots & 0   \\
  \end{pmatrix} \hspace{1 cm} \textbf{E}^T_2=
      \begin{pmatrix}
    1 & 0 & \dots & 0   \\
  \end{pmatrix} 
\end{equation*}
where, $\textbf{E}^T_1$ and $\textbf{E}^T_2$ are respectively $1 \times r$ and $1 \times s$ vectors. We use that along with Equation (\ref{Hankel2}) to write the following balanced expression for $h_k^\delta$:
\begin{equation}
    \begin{aligned}
        h_k^\delta &= \textbf{E}^T_1 \textbf{H}_{r \times s} (g) \textbf{E}^T_2 \\
        &= \textbf{E}^T_1 (\mathcal{O} \textbf{A}_r^{g} \mathcal{C}) \textbf{E}_2 \\
        &= (\textbf{E}^T_1 \Tilde{\textbf{U}} \Tilde{\Sigma}^{\frac{1}{2}}) (\Tilde{\Sigma}^{-\frac{1}{2}} \Tilde{\textbf{U}}^{T} \textbf{H}_{r \times s} (1) \Tilde{\textbf{V}} \Tilde{\Sigma}^{-\frac{1}{2}})^{g} (\Tilde{\Sigma}^{\frac{1}{2}} \Tilde{\textbf{V}}^{T} \textbf{E}_2) \\
        &\equiv \textbf{C}_r \textbf{A}_r^{g} \textbf{B}_r
    \end{aligned}
\end{equation}
We use the decomposed expression above to obtain the reduced input and output matrices $\textbf{B}_r$ and $\textbf{C}_r$ as: 
\begin{equation}
    \textbf{B}_r=\Tilde{\Sigma}^{\frac{1}{2}} \Tilde{\textbf{V}}^{T} \textbf{E}_2
\end{equation}
\begin{equation}
    \textbf{C}_r=\textbf{E}^T_1 \Tilde{\textbf{U}} \Tilde{\Sigma}^{\frac{1}{2}}
\end{equation}
\section{Parametric terms constants} \label{Apndx: Parametric constants}
\begin{align}
    G_0&=\omega_o^2-\frac{\eta^2}{\omega_o^2}a_o^2+\frac{3\gamma}{2}a_o^2+\frac{3 \gamma \eta^2}{4 \omega_o^4}a_o^4 + \frac{\gamma \eta^2}{24 \omega_o^4}a_o^4\\
    G_1&= 2\eta a_o - \frac{5 \gamma \eta}{2 \omega_o^2} a_o^3 \\
    G_2&= \frac{\eta^2}{3 \omega_o^2}a_o^2 + \frac{3 \gamma}{2} a_o^2 - \frac{\gamma \eta^2}{2\omega_o^4}a_o^4\\
    G_3 &= \frac{\gamma \eta}{2 \omega_o^2} a_o^3 \\
    G_4 &= \frac{\gamma \eta^2}{24 \omega_o^4} a_o^4\\
K_0&=-\omega_n^2 + 3\gamma \left(\frac{R_1}{2}+\frac{R_3^2}{2}+\frac{R_5^2}{2}\right)\\
K_2&=3\gamma \left(\frac{R_1^2}{2}+R_1 R_3+R_3 R_5 \right)\\
K_4&=3\gamma \left(R_1 R_3+R_1 R_5 \right)\\
K_6&=3\gamma \left(\frac{R_3^2}{2}+R_1 R_5\right)\\
K_8&=3\gamma \left(R_3 R_5 \right)\\
K_{10}&=3\gamma \left(\frac{R_5^2}{2}\right)
\end{align}
where:
\begin{align}
R_1 &= a_o\\
R_3 &= \frac{\gamma}{32 \Omega^2} a_o^3 + \frac{3 \gamma^2}{1024 \Omega^4} a_o^5\\
R_5 &=\frac{\gamma^2}{1024 \Omega^4} a_o^5
\end{align}

\end{appendices}
\bibliographystyle{unsrt}
%\newpage
\bibliography{main}

\begin{thebibliography}{10}

\bibitem{al2019point}
Elie Al~Shami, Ran Zhang, and Xu~Wang.
\newblock Point absorber wave energy harvesters: A review of recent
  developments.
\newblock {\em Energies}, 12(1):47, 2019.

\bibitem{Falnes2012}
Johannes Falnes and J{\o}rgen Hals.
\newblock Heaving buoys, point absorbers and arrays.
\newblock {\em Philosophical Transactions of the Royal Society A: Mathematical,
  Physical and Engineering Sciences}, 370(1959):246--277, January 2012.

\bibitem{drew2009review}
Benjamin Drew, Andrew~R Plummer, and M~Necip Sahinkaya.
\newblock A review of wave energy converter technology.
\newblock {\em Proceedings of the Institution of Mechanical Engineers, Part A:
  Journal of Power and Energy}, 223(8):887--902, 2009.

\bibitem{younesian2017multi}
Davood Younesian and Mohammad-Reza Alam.
\newblock Multi-stable mechanisms for high-efficiency and broadband ocean wave
  energy harvesting.
\newblock {\em Applied energy}, 197:292--302, 2017.

\bibitem{schubert2020performance}
Benjamin~W Schubert, William~SP Robertson, Benjamin~S Cazzolato, Mergen~H
  Ghayesh, and Nataliia~Y Sergiienko.
\newblock Performance enhancement of submerged wave energy device using
  bistability.
\newblock {\em Ocean Engineering}, 213:107816, 2020.

\bibitem{daqaq2014role}
Mohammed~F Daqaq, Ravindra Masana, Alper Erturk, and D~Dane~Quinn.
\newblock On the role of nonlinearities in vibratory energy harvesting: a
  critical review and discussion.
\newblock {\em Applied Mechanics Reviews}, 66(4), 2014.

\bibitem{xi2021high}
Ru~Xi, Haicheng Zhang, Huai Zhao, Ramnarayan Mondal, et~al.
\newblock High-performance and robust bistable point absorber wave energy
  converter.
\newblock {\em Ocean Engineering}, 229:108767, 2021.

\bibitem{xiao2017comparative}
Xiaolong Xiao, Longfei Xiao, and Tao Peng.
\newblock Comparative study on power capture performance of oscillating-body
  wave energy converters with three novel power take-off systems.
\newblock {\em Renewable Energy}, 103:94--105, 2017.

\bibitem{zhang2019efficiency}
Haicheng Zhang, Ru~Xi, Daolin Xu, Kai Wang, Qijia Shi, Huai Zhao, and Bo~Wu.
\newblock Efficiency enhancement of a point wave energy converter with a
  magnetic bistable mechanism.
\newblock {\em Energy}, 181:1152--1165, 2019.

\bibitem{zhang2016oscillating}
Xian-tao Zhang, Jian-min Yang, and Long-fei Xiao.
\newblock An oscillating wave energy converter with nonlinear snap-through
  power-take-off systems in regular waves.
\newblock {\em China Ocean Engineering}, 30(4):565--580, 2016.

\bibitem{zhang2018application}
Xiantao Zhang, Xinliang Tian, Longfei Xiao, Xin Li, and Lifen Chen.
\newblock Application of an adaptive bistable power capture mechanism to a
  point absorber wave energy converter.
\newblock {\em Applied Energy}, 228:450--467, 2018.

\bibitem{zhang2019mechanism}
Xiantao Zhang, XinLiang Tian, Longfei Xiao, Xin Li, and Wenyue Lu.
\newblock Mechanism and sensitivity for broadband energy harvesting of an
  adaptive bistable point absorber wave energy converter.
\newblock {\em Energy}, 188:115984, 2019.

\bibitem{song2020performance}
Yang Song, Xiaoxian Guo, Hongchao Wang, Xinliang Tian, Handi Wei, and Xiantao
  Zhang.
\newblock Performance analysis of an adaptive bistable point absorber wave
  energy converter under white noise wave excitation.
\newblock {\em IEEE Transactions on Sustainable Energy}, 2020.

\bibitem{hulme1982wave}
A~Hulme.
\newblock The wave forces acting on a floating hemisphere undergoing forced
  periodic oscillations.
\newblock {\em Journal of Fluid Mechanics}, 121:443--463, 1982.

\bibitem{ogilvie1969rational}
T~Francis Ogilvie and Ernest~O Tuck.
\newblock A rational strip theory of ship motions: part i.
\newblock Technical report, University of Michigan, 1969.

\bibitem{haskind2010exciting}
MD~Haskind.
\newblock The exciting forces and wetting of ships in waves.
\newblock {\em report}, 2010.

\bibitem{newman1962exciting}
John~Nicholas Newman.
\newblock The exciting forces on fixed bodies in waves.
\newblock {\em Journal of ship research}, 6(04):10--17, 1962.

\bibitem{nayfeh2008perturbation}
Ali~H Nayfeh.
\newblock {\em Perturbation methods}.
\newblock John Wiley \& Sons, 2008.

\bibitem{brunton2019data}
Steven~L Brunton and J~Nathan Kutz.
\newblock {\em Data-driven science and engineering: Machine learning, dynamical
  systems, and control}.
\newblock Cambridge University Press, 2019.

\bibitem{wang2017modelling}
Liguo Wang.
\newblock {\em Modelling and advanced control of fully coupled wave energy
  converters subject to constraints: the wave-to-wire approach}.
\newblock PhD thesis, Acta Universitatis Upsaliensis, 2017.

\bibitem{meyers2011mathematics}
Robert~A Meyers.
\newblock {\em Mathematics of complexity and dynamical systems}.
\newblock Springer Science \& Business Media, 2011.

\bibitem{kovacic2018mathieu}
Ivana Kovacic, Richard Rand, and Si~Mohamed~Sah.
\newblock Mathieu's equation and its generalizations: overview of stability
  charts and their features.
\newblock {\em Applied Mechanics Reviews}, 70(2), 2018.

\end{thebibliography}
\end{document}